\documentclass[10pt,letter]{amsart}
\usepackage{amssymb,amsbsy,amsmath,amsfonts,amssymb}
\usepackage{latexsym,euscript,exscale}
\title[Turbulence, amalgamation, and generic automorphisms]{Turbulence, amalgamation, and generic automorphisms
of homogeneous structures}
\author {Alexander S. Kechris and Christian Rosendal}
\date {October 2004}
\newcommand {\ca} {{2^\N}}

\newcommand {\N}{\mathbb N}

\newcommand {\Z}{\mathbb Z}

\newcommand{\norm}[1]{\lVert#1\rVert}

\newcommand{\iso}{\cong}

\newcommand{\tom} {\emptyset}

\newcommand{\saa}{\Longrightarrow}

\newcommand{\til}{\rightarrow}

\newcommand{\Lim}[1]{\mathop{\longrightarrow}\limits_{#1}}

\newcommand {\del}{ \; \big| \;}

\newcommand {\ku} {\mathcal}

\newcommand{\ov}{\overline}
\newcommand{\inv}{^{-1}}

\newcommand{\conti}{2^{\aleph_0}}

\renewcommand {\pmb} {\boldsymbol}

\newtheorem{thm}{Theorem}[section]
\newtheorem{cor}[thm]{Corollary}
\newtheorem{lemme}[thm]{Lemma}
\newtheorem{prop} [thm] {Proposition}
\newtheorem{defi} [thm] {Definition}

\newtheorem{rem}[thm]{Remark}

\newcommand{\bbZ}{\mathbb Z}
\newcommand{\bbN}{\mathbb N}

\newcommand{\bbQ}{\mathbb Q}

\newcommand{\bfa}{\mathbf A}
\newcommand{\bfb}{\mathbf B}
\newcommand{\bfc}{\mathbf C}
\newcommand{\bfd}{\mathbf D}
\newcommand{\bfe}{\mathbf E}
\newcommand{\bff}{\mathbf F}
\newcommand{\bfg}{\mathbf G}
\newcommand{\bfh}{\mathbf H}

\newcommand{\bfj}{\mathbf J}
\newcommand{\bfk}{\mathbf K}
\newcommand{\bfl}{\mathbf L}
\newcommand{\bfm}{\mathbf M}
\newcommand{\bfn}{\mathbf N}

\newcommand{\bfp}{\mathbf P}
\newcommand{\bfq}{\mathbf Q}
\newcommand{\bfr}{\mathbf R}
\newcommand{\bfs}{\mathbf S}
\newcommand{\bft}{\mathbf T}
\newcommand{\bfu}{\mathbf U}

\newcommand{\bfuo}{\mathbf U_0}

\newcommand{\aaa }{\mathcal A}
\newcommand{\bbb }{\mathcal B}

\newcommand{\e}{\mathcal E}
\newcommand{\f}{\mathcal F}
\newcommand{\g}{\mathcal G}
\newcommand{\h}{\mathcal H}
\newcommand{\jjj}{\mathcal J}
\newcommand{\kkk}{\mathcal K}
\newcommand{\llll}{\mathcal L}

\newcommand{\n}{\mathcal N}
\newcommand{\mba}{\mathcal {MBA}}
\newcommand{\ooo}{\mathcal O}
\newcommand{\p}{\mathcal P}

\newcommand{\rrr}{\mathcal R}
\newcommand{\s}{\mathcal S}
\newcommand{\ttt}{\mathcal T}
\newcommand{\calu}{\mathcal U}

\def\fraisse{Fra\"\i ss\'e }
\def\abert{Ab\'ert }

\begin{document}
\maketitle

\begin{abstract} We study topological properties of conjugacy classes in
Polish groups, with emphasis on automorphism groups of homogeneous
countable structures. We first consider the existence of dense
conjugacy classes (the topological Rokhlin property). We then
characterize when an automorphism group admits a comeager
conjugacy class (answering a question of Truss) and apply this to
show that the homeomorphism group of the Cantor space has a
comeager conjugacy class (answering a question of
Akin-Hurley-Kennedy). Finally, we study Polish groups that admit
comeager conjugacy classes in any dimension (in which case the
groups are said to admit ample generics). We show that Polish
groups with ample generics have the small index property
(generalizing results of Hodges-Hodkinson-Lascar-Shelah) and
arbitrary homomorphisms from such groups into separable groups are
automatically continuous.  Moreover, in the case of oligomorphic
permutation groups, they have uncountable cofinality and the
Bergman property. These results in particular apply to
automorphism groups of many $\omega$-stable, $\aleph_0$-categorical
structures and of the random graph. In this connection, we also show that the infinite symmetric group $S_\infty$ has a unique non-trivial separable group topology. For several interesting groups we
also establish Serre's properties (FH) and (FA).
\end{abstract}

\section{Introduction}

\subsection{Polish groups} We study in this paper topological properties of
conjugacy classes in Polish groups.  There are two questions which
we are particularly interested in.  First, does a Polish group $G$
have a dense conjugacy class?  This is equivalent (see, e.g.,
Kechris \cite{kechris1} (8.47)) to the following generic ergodicity property
of $G$:  Every conjugacy invariant subset $A\subseteq G$ with the
Baire property (e.g., a Borel set) is either meager or comeager.
There is an extensive list of Polish groups that have dense
conjugacy classes, like, e.g., the automorphism group Aut$(X,\mu
)$ of a standard measure space $(X,\mu )$, i.e., a standard Borel
space $X$ with a non-atomic Borel probability measure $\mu$ (see,
e.g., Halmos \cite{halmos}), the unitary group $U(H)$ of separable
infinite-dimensional Hilbert space $H$ (see, e.g., Choksi-Nadkarni
\cite{choksi}) and the homeomorphism groups $H(X)$ of various compact
metric spaces $X$ with the uniform convergence topology, including
$X=[0,1]^\bbN$ (the Hilbert cube), $X=2^\bbN$ (the Cantor space),
$X=S^{2d}$ (even dimensional spheres), etc. (see Glasner-Weiss
\cite{glasner} and Akin-Hurley-Kennedy \cite{akin}).  In Glasner-Weiss \cite{glasner} groups
that have dense conjugacy classes are said to have the {\it
topological Rokhlin property}, motivated by the existence of dense
conjugacy classes in Aut$(X,\mu )$, which is usually seen as a
consequence of the well-known Rokhlin Lemma in ergodic theory.
Recently, Akin, Glasner and Weiss \cite{akin2} also found an example of a locally compact
Polish group with dense conjugacy class. Any such group must be
non-compact and it appears likely that it must also be totally
disconnected. (Karl Hofmann has shown that if a locally compact
group $\neq\{1\}$ has a dense conjugacy class, then it is not
pro-Lie and in particular, its quotient by the connected component
of the identity cannot be compact.)

The second question we consider is whether $G$ has a (necessarily unique) dense $G_\delta$ conjugacy class (which is well-known to
be equivalent to whether it has a dense non-meager class, see, e.g., Becker-Kechris \cite{becker}). Following Truss \cite{truss}, we call any
element of $G$ whose conjugacy class is dense $G_\delta$ a {\it generic element} of $G$. This is a much stronger property which
fails in very ``big'' groups such as Aut$(X,\mu )$ or $U(H)$ but
can often occur in automorphism groups Aut$(\bfk )$ of countable
structures $\bfk$. It was first studied in this context by Lascar
\cite{lascar} and Truss \cite{truss}. For example, Aut$(\bbN ,=)$ (i.e., {\it the
infinite symmetric group} $S_\infty$), Aut$(\bbQ ,<)$, Aut$(\bfr
)$, where $\bfr$ is the random graph, Aut$(\bfp )$, where $\bfp$
is the random poset, all have a dense $G_\delta$ conjugacy class
(see Truss \cite{truss}, Kuske-Truss \cite{kuske}).  For more on automorphism
groups having dense $G_\delta$ classes, see also the recent paper
Macpherson-Thomas \cite{macpherson2}.  But the question of whether certain groups have a dense $G_\delta$ conjugacy class arose as well in
topological dynamics, where Akin-Hurley-Kennedy \cite{akin} p. 104 posed the problem of the existence of a dense $G_\delta$ conjugacy class in ${\it H}(2^\bbN )$, i.e., the existence of a generic homeomorphism of the Cantor space.

\subsection{Model theory} Our goal here is to study these questions in the context
of automorphism groups, Aut$(\bfk )$, of countable structures
$\bfk$.  It is well-known that these groups are (up to topological
group isomorphism) exactly the closed subgroups of $S_\infty$.
However, quite often such groups can be densely embedded in other
Polish groups $G$ (i.e., there is an injective continuous
homomorphism from Aut$(\bfk )$ into $G$ with dense image) and
therefore establishing the existence of a dense conjugacy class of
Aut$(\bfk )$ implies the same for $G$.  So we can use this method
to give simple proofs that such conjugacy classes exist in several
interesting Polish groups. This can be viewed as another instance
of the idea of reducing questions about the structure of certain
Polish groups to those of automorphism groups Aut$(\bfk )$ of
countable structures $\bfk$, where one can employ methods of model
theory and combinatorics.  An earlier use of this methodology is
found in Kechris-Pestov-Todorcevic \cite{kechris4} in connection with the study
of extreme amenability and its relation with Ramsey theory.

It is also well-known that to every closed subgroup $G\leq
S_\infty$ one can associate the structure $\bfk_G=(\bbN
,\{R_{i,n}\})$, where $R_{i,n}\subseteq\bbN^n$ is the $i$th orbit
in some fixed enumeration of the orbits of $G$ on $\bbN^n$, which
is such that $G=$ Aut$(\bfk_G)$ and moreover $\bfk_G$ is {\it
ultrahomogeneous}, i.e., every isomorphism between finite
substructures of $\bfk_G$ extends to an automorphism of $\bfk_G$.
\fraisse has analyzed such structures in terms of their finite
``approximations."

To be more precise, let $\kkk$ be a class of finite structures in a
fixed (countable) signature $L$, which has the following
properties:

(i) (HP) $\kkk$ is {\it hereditary} (i.e., $\bfa\leq\bfb\in\kkk$
implies $\bfa\in\kkk$, where $\bfa\leq\bfb$ means that $\bfa$ can be
embedded into $\bfb$),

(ii) (JEP) $\kkk$ satisfies the {\it joint embedding property}
(i.e., if $\bfa ,\bfb\in\kkk$, there is $\bfc\in\kkk$ with $\bfa
,\bfb\leq\bfc$),

(iii) (AP) $\kkk$ satisfies the {\it amalgamation property} (i.e.,
if $f: \bfa\rightarrow\bfb ,g:\bfa\rightarrow\bfc$, with $\bfa
,\bfb ,\bfc\in\kkk$, are embeddings, there is $\bfd\in\kkk$ and
embeddings $r:\bfb\rightarrow\bfd , s:\bfc\rightarrow\bfd$ with
$r\circ f=s\circ g$).

(iv) $\kkk$ contains only countably many structures, up to
isomorphism, and contains structures of arbitrarily large (finite)
cardinality.

We call any $\kkk$ that satisfies (i)-(iv) a {\it \fraisse class}.
Examples include the class of trivial structures $(L=\emptyset )$,
graphs, linear orderings, Boolean algebras, metric spaces with
rational distances, etc. For any \fraisse class $\kkk$ one can
define, following \fraisse \cite{fraisse} (see also Hodges \cite{hodges1}), its
so-called {\it \fraisse limit} $\bfk =$ Flim$(\kkk )$, which is the
unique countably infinite structure satisfying:

(a) $\bfk$ is locally finite (i.e., finitely generated
substructures of $\bfk$ are finite),

(b) $\bfk$ is {\it ultrahomogeneous} (i.e., any isomorphism
between finite substructures of $\bfk$ extends to an automorphism
of $\bfk$),

(c) Age$(\bfk )=\kkk$, where Age$(\bfk )$ is the class of all finite
structures that can be embedded in $\bfk$.

A countably infinite structure $\bfk$ satisfying (a), (b) is
called a {\it \fraisse structure}.  If $\bfk$ is a \fraisse
structure, then Age$(\bfk )$ is a \fraisse class, therefore the
maps $\kkk\mapsto$ Flim$(\kkk )$, $\bfk\mapsto$ Age$(\bfk )$ provide a
canonical bijection between \fraisse classes and structures.
Examples of \fraisse structures include the trivial structure
$(\bbN , =),\ \bfr =$ the random graph, $(\bbQ ,<), \bfb_\infty =$
the countable atomless Boolean algebra, $\bfuo =$ the rational
Urysohn space (= the \fraisse limit of the class of finite metric
spaces with rational distances), etc.

For further reference, we note that condition (b) in the
definition of a \fraisse structure can be replaced by the
following equivalent condition, called the {\it extension
property}:  If $\bfa ,\bfb$ are finite, $\bfa , \bfb\leq\bfk$, and
$f:\bfa\rightarrow\bfk ,g:\bfa\rightarrow\bfb$ are embeddings,
then there is an embedding $h:\bfb\rightarrow\bfk$ with $h\circ
g=f$.

Thus every closed subgroup $G\leq S_\infty$ is of the form $G=$
Aut$(\bfk )$ for a \fraisse structure $\bfk$.  We study in this
paper the question of existence of dense or comeager conjugacy
classes in the Polish groups Aut$(\bfk )$, equipped with the
pointwise convergence topology, for \fraisse structures $\bfk$, in
terms of properties of $\bfk$, and also derive consequences for
other groups.

\subsection{Dense conjugacy classes}  Following Truss \cite{truss}, we associate to each \fraisse
class $\kkk$ the class $\kkk_p$ of all systems, $\s =\langle\bfa ,\psi
:\bfb\rightarrow \bfc\rangle$, where $\bfa ,\bfb ,\bfc\in\kkk
,\;\bfb ,\bfc\subseteq\bfa$ (i.e., $\bfb ,\bfc$ are {\it
substructures} of $\bfa$) and $\psi$ is an isomorphism of $\bfb$
and $\bfc$.  For such systems, we can define the notion of
embedding as follows:  An {\it embedding} of $\s$ into $\ttt
=\langle\bfd , \varphi :\bfe\rightarrow\bff\rangle$ is an
embedding $f:\bfa\rightarrow \bfd$ such that $f$ embeds $\bfb$
into $\bfe ,\bfc$ into $\bff$ and $f\circ\psi\subseteq\phi\circ
f$.  Using this concept of embedding it is clear then what it
means to say that $\kkk_p$ satisfies JEP or AP.

We now have:
\begin{thm} Let $\kkk$ be a \fraisse class with
\fraisse limit $\bfk =$ {\rm Flim}$(\kkk )$.  Then the following are
equivalent:

(i) There is a dense conjugacy class in {\rm Aut}$(\bfk )$.

(ii) $\kkk_p$ satisfies the {\rm JEP}.
\end{thm}
For example, it is easy to verify JEP of $\kkk_p$ for the following
classes $\kkk$:

(i) $\kkk =$ finite metric spaces with rational distances,

(ii) $\kkk =$ finite Boolean algebras,

(iii) $\kkk =$ finite measure Boolean algebras with rational
measure,

\noindent and this immediately gives that Aut$(\bfuo )={\rm Iso}
(\bfuo )$ (= the isometry group of the rational Urysohn space
$\bfuo$), Aut$(\bfb_\infty )$, Aut$(\bff ,\lambda )$ (the
automorphism group of the Boolean algebra generated by the
rational intervals of $[0,1]$ with Lebesgue measure) all have
dense conjugacy classes.  But Aut$(\bfuo )$ can be densely
embedded into Iso$(\bfu )=$ the isometry group of the Urysohn
space (see Kechris-Pestov-Todorcevic \cite{kechris4}), Aut$(\bfb_\infty )$ is
isomorphic to ${\it H}(2^\bbN )$ by Stone duality, and Aut$(\bff
,\lambda )$ can be densely embedded in Aut$(X,\mu )$, so we have
simple proofs of the following:
\begin{cor} 
The following Polish groups have dense conjugacy
classes:

(i) {\rm Iso}$(\bfu )$,

(ii) {\rm(Glasner-Weiss \cite{glasner}, Akin-Hurley-Kennedy \cite{akin})} ${\it H}(2^\bbN )$,

(iii) {\rm (Rokhlin)} {\rm Aut}$(X,\mu )$.
\end{cor}
Glasner and Pestov have also proved part (\emph{i}). We also
obtain a similar result for the (diagonal) conjugacy action of $G$
on $G^{\bbN}$ for all the above groups $G$.

We say that a Polish group $G$ has a {\it cyclically dense
conjugacy class} if there are $g,h \in G$ such that
$\{g^nhg^{-n}\}_{n\in \Z}$ is dense in $G$. In this case $G$ is
topologically 2-generated, i.e., has a dense 2-generated subgroup.
For example, $S_{\infty}$ has this property and so does the automorphism group of the random graph (see Macpherson \cite{macpherson1}). Using a version of Theorem 1, we can give simple proofs that the following groups admit cyclically dense conjugacy classes and therefore are topologically 2-generated: ${\it H}(2^{\bbN})$, ${\it H}(2^{\bbN},\sigma)$ (the group of the measure preserving
homeomorphisms of $2^{\bbN}$ with the usual product measure), Aut$(X,\mu)$, and Aut$(\bbN^{< \bbN})$ (the automorphism group of
the infinitely splitting rooted tree).

\subsection{Comeager conjugacy classes}
We now turn to the existence of a dense $G_\delta$
conjugacy class in Aut$(\bfk ),\bfk$ a \fraisse structure, i.e.,
the existence of a {\it generic automorphism} of $\bfk$. Truss
\cite{truss} showed that the existence of a subclass $\llll\subseteq\kkk_p$ such that $\llll$ is cofinal under embeddability and satisfies the AP, a
property which we will refer to as the {\it cofinal amalgamation
property} (CAP), together with the JEP for $\kkk_p$, is sufficient
for the existence of a generic automorphism. He also raised the
question of whether the existence of a generic automorphism is
equivalent to some combination of amalgamation and joint embedding
properties for $\kkk_p$.

Motivated by this problem, we have realized that the question of
generic automorphisms is closely related to Hjorth's concept of
turbulence (see, Hjorth \cite{hjorth} or Kechris \cite{kechris2}) for the conjugacy
action of the automorphism group Aut$(\bfk )$, a connection that
is surprising at first sight since turbulence is a phenomenon
usually thought of as incompatible with actions of closed
subgroups of the infinite symmetric group, like Aut$(\bfk )$.
Once however this connection is realized, it leads naturally to
the formulation of an appropriate amalgamation property for $\kkk_p$
which in combination with JEP is equivalent to the existence of a
generic automorphism.

We say that $\kkk_p$ satisfies the {\it weak amalgamation property}
(WAP) if for any $\s =\langle\bfa ,\psi
:\bfb\rightarrow\bfc\rangle$ in $\kkk_p$, there is $\ttt =\langle\bfd
,\varphi :\bfe\rightarrow\bff\rangle$ and an embedding
$e:\s\rightarrow\ttt$, such that for any embeddings $f:\ttt\rightarrow
\ttt_0 ,g:\ttt\rightarrow\ttt_1$, where $\ttt_0,\ttt_1\in\kkk_p$, there is
$\calu\in\kkk_p$ and embeddings $r:\ttt_0\rightarrow\calu ,s:
\ttt_1\rightarrow\calu$ with $r\circ f\circ e=s\circ g\circ e$.

We now have:
\begin{thm} 
Let $\kkk$ be a \fraisse class and $\bfk =$
{\rm Flim}$(\kkk )$ its \fraisse limit.  Then the following are
equivalent:

(i) $\bfk$ has a generic automorphism.

(ii) $\kkk_p$ satisfies {\rm JEP} and {\rm WAP}.
\end{thm}
After obtaining this result, we found out that Ivanov \cite{ivanov} had
already proved a similar theorem, in response to Truss' question,
in a somewhat different context, that of $\aleph_0$-categorical
structures (he also calls WAP the {\it almost amalgamation
property}).  Our approach however, through the idea of turbulence,
is different and we proceed to explain it in more detail.

Suppose a Polish group $G$ acts continuously on a Polish space
$X$.  Given $x\in X$, an open nbhd $U$ of $x$ and an open
symmetric nbhd $V$ of the identity of $G$, the $(U,V)$-{\it local
orbit} of $x, \ooo (x,U,V)$, is the set of all $y\in U$ for which
there is a finite sequence $g_0,g_1,\dots ,g_k\in V$ with $x_0=x,
g_i\cdot x_i=x_{i+1},x_{k+1}=y$ and $x_i\in U,\forall i$. A point
$x$ is {\it turbulent} if for every $U$, $V$ as above
Int($\overline{\ooo (x,U,V)})\neq\emptyset$.  This turns out to be
equivalent to saying that $x\in{\rm Int}(\overline{\ooo (x,U,V)})$,
see Kechris \cite{kechris2}. It is easy to see that this only depends on the
orbit $G\cdot x$ of $x$, so we can refer to {\it turbulent
orbits}.  This action is called {\it (generically) turbulent} if:

(i) Every orbit is meager.

(ii) There is $x\in X$ with dense, turbulent orbit.

(This is not quite the original definition of turbulence, as in
Hjorth \cite{hjorth}, but it is equivalent to it, see, e.g., Kechris \cite{kechris2}.)

Examples of turbulent actions, relevant to our context, include
the conjugacy actions of $U(H)$, see Kechris-Sofronidis \cite{kechris3}, and
Aut$(X,\mu )$, see Foreman-Weiss \cite{foreman}.

Now Hjorth \cite{hjorth} has shown that no closed subgroup of $S_\infty$
has a turbulent action and from this one has the following
corollary, which can be also easily proved directly.
\begin{prop} 
Let $G$ be a closed subgroup of $S_\infty$ and suppose $G$
acts continuously on the Polish space $X$.  Then the following are
equivalent for any $x\in X$:

(i) The orbit $G\cdot x$ is dense $G_\delta$.

(ii) $G\cdot x$ is dense and turbulent.
\end{prop}
Thus a dense $G_\delta$ orbit exists iff a dense turbulent orbit
exists.

We can now apply this to the conjugacy action of Aut$(\bfk )$, $\bfk$
a \fraisse structure, on itself by conjugacy and this leads to the
formulation of the WAP and our approach to the proof of Theorem 3.
In view of that result, it is interesting that, in the context of
this action, the existence of dense, turbulent orbits is
equivalently manifested as a combination of joint embedding and
amalgamation properties.

\subsection{Homeomorphisms of the Cantor space}
We next use these ideas to answer the question of
Akin-Hurley-Kennedy \cite{akin} about the existence of a generic
homeomorphism of the Cantor space $2^\bbN$, i.e., the existence of
a dense $G_\delta$ conjugacy class in ${\it H}(2^\bbN )$.  Since
${\it H}(2^\bbN )$ is isomorphic (as a topological group) to
Aut$(\bfb_\infty )$, where $\bfb_\infty$ is the countable atomless
Boolean algebra (i.e., the \fraisse limit of the class of finite
Boolean algebras), this follows from the following result.
\begin{thm} 
Let $\bbb\aaa $ be the class of finite Boolean algebras. Then
$\bbb\aaa _p$ has the {\rm CAP}. So $\bfb_\infty$ has a generic
automorphism and there is a generic homeomorphism of the Cantor
space. 
\end{thm}
Using well-known results, we can also see that there is a generic
element of the group of order-preserving homeomorphisms of the interval
$[0,1]$.

\subsection{Ample generics}
Finally, we discuss the concept of ample generics in
Polish groups, which is a tool that has been used before (see,
e.g., Hodges et al. \cite{hodges2}) in the structure theory of automorphism
groups. We say that a Polish group $G$ has {\it ample generic}
elements if for each finite $n$ there is a comeager orbit for the
(diagonal) conjugacy action of $G$ on $G^n:g\cdot (g_1,\dots
,g_n)=(gg_1g^{-1},\dots ,gg_ng^{-1})$. Obviously this is a
stronger property than just having a comeager conjugacy class and
for example ${\rm Aut}(\bbQ, <)$ has the latter, but not the former
(see Kuske and Truss \cite{kuske} for a discussion of this).

There is now an extensive list of permutation groups known to have
ample generics. These include the automorphism groups of the
following structures: many $\omega$-stable, $\aleph_0$-categorical
structures (see Hodges et al. \cite{hodges2}), the random graph (Hrushovski
\cite{hrushovski}, see also Hodges et al. \cite{hodges2}), the free group on countably
many generators (Bryant and Evans \cite{bryant}) and arithmetically
saturated models of true arithmetic (Schmerl \cite{schmerl}). Moreover,
Herwig and Lascar \cite{herwig} have extended the result of Hrushovski to a
much larger class of structures in finite relational languages and
the isometry group of the rational Urysohn space $\bfu_0$ is now
also known to have ample generics (this follows from recent
results of Solecki \cite{solecki} and Vershik).

We will add another two groups to this list, which incidentally
are automorphism groups of structures that are not
$\aleph_0$-categorical, namely, the group of (Haar) measure
preserving homeomorphisms of the Cantor space, $H(2^\bbN,\sigma)$,
and the group of Lipschitz homeomorphisms of the Baire space. This
latter group is canonically isomorphic to $\rm Aut(\bbN^{<\bbN})$,
where $\bbN^{<\bbN}$ is seen as the infinitely splitting regular
rooted tree.

Let us first notice that in the same manner as for the existence
of a comeager conjugacy class, we are able to determine an
equivalent model-theoretic condition on a \fraisse class $\kkk$ for
when the automorphism group of its \fraisse limit $\bfk$ has
ample generic elements. This criterion is in fact a trivial
generalization of the one dimensional case. But we shall be more
interested in the consequences of the existence of ample generics.

Recall that a second countable topological group $G$ is said to have the {\it small
index property} if any subgroup of index $<2^{\aleph_0}$ is open.
Then we can show the following, which generalizes the case of
automorphism groups of $\omega$-stable, $\aleph_0$-categorical
structures due to Hodges et al. \cite{hodges2}.
\begin{thm} 
Let $G$ be a Polish group with ample generic
elements. Then $G$ has the small index property.
\end{thm}
In the case of $G$ being a closed subgroup of $S_\infty$, i.e.,
$G$ having a neighborhood basis at the identity consisting of
clopen subgroups, this essentially says that the topological
structure of the group is completely determined by its algebraic
structure.

We subsequently study the cofinality of Polish groups. Recall that
the {\it cofinality} of a group $G$ is the least cardinality of a
well ordered chain of proper subgroups whose union is $G$. Again
generalizing results of Hodges et al. \cite{hodges2}, we prove:
\begin{thm} 
Let $G$ be a Polish group with ample generic
elements. Then $G$ is not the union of countably many non-open
subgroups (or even cosets of subgroups).
\end{thm}
Note that if $G$, a closed subgroup of $S_\infty$, is
oligomorphic, i.e., has only finitely many orbits on each
$\bbN^n$, then, by a result of Cameron, any open subgroup of $G$
is contained in only finitely many subgroups of $G$, thus, if $G$
has ample generics, it has uncountable cofinality. The same holds
for connected Polish groups, and Polish groups with a finite
number of topological generators.

It turns out that the existence of generic elements of a Polish
group has implications for its actions on trees. So let us recall
some basic notions of the theory of group actions on trees (see
Serre \cite{serre}):

A group $G$ is said to act {\it without inversion} on a tree $T$
if $G$ acts on $T$ by automorphisms such that for no $g\in G$
there are two adjacent vertices $a,b\in T$ such that $g\cdot a=b$
and $g\cdot b=a$. The action is said to have a fixed point if
there is an $a\in T$ such that $g\cdot a=a$, for all $g\in G$. We
say that a group $G$ has {\it property} (FA) if whenever $G$ acts
without inversion on a tree, there is a fixed point.  When $G$ is
not countable this is known to be equivalent to the conjunction of
the following three properties (Serre \cite{serre}):

(i) $G$ is not a non-trivial free product with amalgamation,

(ii) $\bbZ$ is not a homomorphic image of $G$,

(iii) $G$ has uncountable cofinality.

Macpherson and Thomas \cite{macpherson2} recently showed that (i) follows if $G$
has a comeager conjugacy class and, as (ii) trivially also holds
in this case, we are left with verifying (iii).

Another way of proving property (FA) is through a slightly different study of the generation of Polish groups with ample generics. Obviously, any generating set for an uncountable group must be uncountable, but ideally we would still
like to understand the structure of the group by studying a set of generators. Let us say that a group $G$ has the {\it Bergman
property} if for each exhaustive sequence of subsets $W_0\subseteq W_1\subseteq W_2\subseteq \ldots\subseteq G$, there are $n$ and
$k$ such that $W_n^k=G$. Bergman \cite{bergman} showed this property for $S_\infty$ by methods very different from those employed here and
we will extend his result to a fairly large class of automorphism groups:
\begin{thm} 
Let $G$ be a closed oligomorphic subgroup of
$S_\infty$ with ample generic elements. Then $G$ has the Bergman
property.
\end{thm}
In particular, this result applies to many automorphism groups of $\omega$-stable, $\aleph_0$-categorical structures and the automorphism group of the random graph.

It is not hard to see, and has indeed been noticed independently by other authors (e.g., Cornulier \cite{cornulier}) that the  Bergman property also implies that any action of the group by isometries on a metric space has bounded orbits. In fact, this is actually an equivalent formulation of the Bergman property. But well-known results of geometric group theory (see B. Bekka, P. de la Harpe and A. Valette \cite{bekka}) state that if a group action by isometries on a real Hilbert space has a bounded orbit, then it has a fixed point. Similarly for an action by automorphisms without inversion on a tree. Thus Bergman groups automatically have property (FH) and (FA), where property (FH) is the statement that any isometric action on a real Hilbert space has a fixed point.
\begin{cor}Let $G$ be a closed oligomorphic subgroup of
$S_\infty$ with ample generic elements. Then $G$ has properties {\rm (FA)} and {\rm (FH)}.
\end{cor}
The phenomenon of automatic continuity is well known and has been
extensively studied, in particular in the context of Banach
algebras. In this category morphisms of course preserve much more
structure than homomorphisms of the underlying groups and
therefore automatic continuity is easier to obtain. But there are
also plenty of examples of this phenomenon for topological groups,
provided we add some definability constraints on the
homomorphisms. An example of this is the classical result of
Pettis, saying that any Baire measurable homomorphism from a
Polish group into a separable group is continuous. Surprisingly
though, when one assumes ample generics one can completely
eliminate any definability assumption and still obtain the same
result. In fact, one does not even need as much as separability
for the target group, but essentially need only rule out that its
topology is discrete. Let us recall that the {\it Souslin number}
of a topological space is the least cardinal $\kappa$ such that
there is no family of $\kappa$ many disjoint open sets. In analogy
with this, if $H$ is a topological group, we let the {\it uniform
Souslin number} of $H$ be the least cardinal $\kappa$ such that
there is no non-$\emptyset$ open set having $\kappa$ many disjoint
translates. Then we can prove the following:
\begin{thm}\label{autom} 
Suppose $G$ is a Polish group with ample
generic elements and $\pi:G\rightarrow H$ is a homomorphism into a
topological group with uniform Souslin number at most
$2^{\aleph_0}$ (in particular, if $H$ is separable). Then $\pi$ is
continuous.
\end{thm}
This in particular shows that any action of a Polish group with
ample generics by isometries on a Polish space or by
homeomorphisms on a compact metric space is actually a continuous
action. For such an action is essentially just a homomorphism into
the isometry group, respectively into the homeomorphism group.

Moreover, one also sees that in this case there is a unique Polish
group topology, and coupled with a result of Gaughan \cite{gaughan},
one can in the case of $S_\infty$ prove the following stronger fact.
\begin{thm}
There is exactly one non-trivial separable group topology
on $S_\infty$.
\end{thm}
In the literature one can find several results on automatic continuity of homomorphisms between topological groups when one puts restrictions on the target groups. For example, the small index property of a group $G$ can be seen to imply that any homomorphism from $G$ into $S_\infty$ is continuous. Also a classical theorem due to Van der Waerden (see
Hofmann and Morris \cite{hofmann}) states that if $\pi\colon G\to
H$ is a group homomorphism from an $n$-dimensional Lie group $G$,
whose Lie algebra agrees with its commutator algebra, into a compact
group $H$, then $\pi$ is continuous. The final result we should
mention is due to Dudley \cite{dudley}, which says that any
homomorphism from a complete metric group into a ``normed'' group with
the discrete topology is automatically continuous. We shall not go
into his definition of a normed group, other than saying that these
include the additive group of the integers and more general free
groups. The novelty of Theorem \ref{autom} lies in the fact that it places no restrictions on the target group (other than essentially ruling out the trivial case of the topology being discrete).

Extending the list of Polish groups with ample generics with
${\it H}(2^\bbN, \sigma)$ and Aut$(\bbN^{<\bbN})$, we finally show the
following, where a closed subgroup of $S_\infty$ has the {\it
strong small index property} if any subgroup of index
$<2^{\aleph_0}$ is sandwiched between the pointwise and setwise
stabilizer of a finite set.
\begin{thm} 
Let $G$ be either ${\it H}(2^\bbN, \sigma)$, the
group of measure preserving homeomorphisms of the Cantor space, or
{\rm Aut}$(\bbN^{<\bbN})$, the group of Lipschitz homeomorphisms
of the Baire space. Then

(i)  $G$ has ample generic elements.

(ii) $G$ has the strong small index property.

(ii) $G$ has uncountable cofinality.

(iv) $G$ has the Bergman property and thus properties {\rm (FH)} and {\rm (FA)}.
\end{thm}
The strong small index property for {\rm Aut}$(\bbN^{<\bbN})$ was
previously proved by R\"ognvaldur M\"oller \cite{moller}.

{\it Acknowledgments}.  The research of ASK was partially
supported by NSF Grants DMS 9987437 and DMS 0455285, the Centre de Recerca
Matem\`{a}tica, Bellatera, and a Guggenheim Fellowship. We would
like to thank Eli Glasner, Karl H. Hofmann, Simon Thomas, Vladimir
Pestov, S\l awek Solecki and the anonymous referee for many useful comments.

\section{Automorphisms with dense conjugacy classes}\label{dense}

Let $\kkk$ be a \fraisse class.  We let $\kkk_p$ be the class of all
systems of the form $\s =\langle\bfa,\psi
:\bfb\rightarrow\bfc\rangle$, where $\bfa ,\bfb ,\bfc\in\kkk ,\;\bfb
,\bfc\subseteq\bfa$ and $\psi$ is an isomorphism of $\bfb$ and
$\bfc$.  An embedding of one system $\s =\langle \bfa ,\psi
:\bfb\rightarrow\bfc\rangle$ into another $\ttt =\langle\bfd , \phi
:\bfe\rightarrow\bff\rangle$ is an embedding
$f:\bfa\rightarrow\bfd$ such that $f$ embeds $\bfb$ into $\bfe
,\bfc$ into $\bff$ and moreover $f\circ\psi \subseteq\phi\circ f$.

So we can define JEP and AP for $\kkk_p$ as well.

\begin{thm}\label{criterion for density} Let $\kkk$ be a \fraisse class with
\fraisse limit $\bfk$.  Then the following are equivalent:

(i) There is a dense conjugacy class in {\rm Aut}$(\bfk )$.

(ii) $\kkk_p$ satisfies the {\rm JEP}.\end{thm}

\begin{proof} (i) $\Rightarrow$ (ii):  Fix some element $f\in$ Aut$(\bfk )$
having a dense conjugacy class in Aut$(\bfk )$ and suppose $\s
=\langle \bfa ,\psi :\bfb\rightarrow\bfc\rangle ,\ttt =\langle\bfd
,\phi :\bfe \rightarrow\bff\rangle$ are two systems in $\kkk_p$.

Replacing $\s$ and $\ttt$ by isomorphic copies we can of course
assume that $\bfa ,\bfd\subseteq\bfk$.  Now by ultrahomogeneity of
$\bfk$ we know that both $\psi$ and $\phi$ have extensions in
Aut$(\bfk )$, so, by the density of $f$'s conjugacy class, there
are $g_1,g_2\in$ Aut$(\bfk )$ such that $\psi\subseteq
g^{-1}_1fg_1$ and $\phi\subseteq g^{-1}_2fg_2$. Let $\bfh =\langle
g''_1\bfa\cup g''_2\bfd\rangle ,\bfm =\langle g''_1 \bfb\cup
g''_2\bfe\rangle ,\bfn =\langle g''_1\bfc\cup g''_2\bff\rangle$
and let $\chi =f\restriction\bfm$.  Then it is easily seen that
$g_1 \restriction_\bfa$ and $g_2\restriction_\bfd$ embed $\s$ and
$\ttt$ into $\langle\bfh ,\chi :\bfm\rightarrow\bfn\rangle$.

(ii) $\Rightarrow$ (i):  A basis for the open subsets of Aut$(\bfk
)$ consists of sets of the form:
\[
[\psi ]=[\psi :\bfb\rightarrow\bfc ]=\{f\in {\rm Aut}(\bfk
):f\supseteq \psi\},
\]
where $\psi :\bfb\rightarrow\bfc$ is an isomorphism of finite
substructures of $\bfk$.  We refer to such $\psi
:\bfb\rightarrow\bfc$ as {\it conditions}. We wish to construct an
element $f\in$ Aut$(\bfk )$ such that for any condition $\psi
:\bfb\rightarrow\bfc$ there is a $g\in$ Aut$(\bfk )$ such that
$g\psi g^{-1}\subseteq f$.  So let
\[
D(\psi :\bfb\rightarrow\bfc )=\{f\in{\rm Aut}(\bfk ):\exists g\in
{\rm Aut}(\bfk )(g\psi g^{-1}\subseteq f)\}.
\]
$D(\psi :\bfb\rightarrow\bfc )$ is clearly open, but we shall see
that it is also dense. If $[\phi : \bfe\rightarrow\bff ]$ is a basic
open set, there is, by the JEP of $\kkk_p$ and the extension
property of $\bfk$, some isomorphism $\chi:\bfh \rightarrow\bfl$ of
finite substructures of $\bfk$ such that $\phi\subseteq \chi$, and
some $g\in$ Aut$(\bfk )$ embedding $\psi :\bfb\rightarrow\bfc$ into
$\chi :\bfh\rightarrow\bfl$.  Then any $f\in$ Aut$(\bfk )$ extending
$\chi$ witnesses that $D(\psi :\bfb\rightarrow\bfc )\cap [\phi
:\bfe\rightarrow \bfd ]\neq\emptyset$.  Therefore any $f\in$
Aut$(\bfk )$ that is in the intersection of all $D(\psi
:\bfb\rightarrow\bfc )$ has dense conjugacy class.\end{proof}

{\bf Remark}.  One can also give a more direct proof of Theorem
1.1 by using the following standard fact:  If $G$ is a Polish
group which acts continuously on a Polish space $X$ and $\bbb$ is a
countable basis of non-$\emptyset$ open sets for the topology of
$X$, then there is a dense orbit for this action iff $\forall
U,V\in\bbb (G\cdot U\cap V\neq\emptyset )$. (The direction
$\Rightarrow$ is obvious.  For the direction $\Leftarrow$, let for
$V\in\bbb ,D_V=\{x\in X:G\cdot x\cap V\neq\emptyset\}$.  This is
clearly open, dense, and so $\bigcap_{V\in\bbb}D_V\neq\emptyset$.
Any $x\in \bigcap_{V\in\bbb}D_V$ has dense orbit.)

We can simply apply this to the conjugacy action of Aut$(\bfk )$
on itself and the basis consisting of the sets $[\psi ]$ as above.
\medskip

We will proceed to some applications of Theorem \ref{criterion for density}.

\begin{thm}\label{dense urysohn} Let $\bfuo$ be the rational Urysohn space. Then
{\rm Aut}$(\bfuo )$ has a dense conjugacy class. Thus, as {\rm
Aut}$(\bfuo )$ can be continuously, densely embedded into {\rm
Iso}$(\bfu )$, the isometry group of the Urysohn space, {\rm
Iso}$(\bfu )$ also has a dense conjugacy class.\end{thm}

\begin{proof} Let $\kkk$ be the class of finite metric spaces with rational
distances; we will see that $\kkk_p$ has the JEP.

Given $\s =\langle\bfa,\psi :\bfb\rightarrow\bfc\rangle$ and $\ttt
=\langle\bfd , \phi :\bfe\rightarrow\bff\rangle\in\kkk_p$,  we let
$\bfh =\bfa\sqcup\bfb$ be the disjoint union of the two metric
spaces $\bfa$ and $\bfb$, where we have put some distance $k>$
diam$(\bfa )+$ diam$(\bfb )$ between any two points $x\in\bfa$ and
$y\in\bfb$.  The triangle inequality is clearly satisfied. So now we
can let $\bfm$ and $\bfn$ be the corresponding unions of $\bfb
,\bfe$ and $\bfc ,\bff$, and finally we let $\chi =\psi
\cup\varphi$.  Then clearly both $\s$ and $\ttt$ embed into $\langle
\bfh ,\chi :\bfm\rightarrow\bfn\rangle$.\end{proof}

Glasner and Pestov informed us that they had also proved that
Iso$(\bfu )$ has a dense conjugacy class.

We next consider the class $\kkk =\mba_\bbQ$ of all finite boolean
algebras with an additive probability measure taking positive
rational values on each nonzero element.  So an element of $\kkk$ is
on the form $(\bfa ,\mu :\bfa\rightarrow [0,1]\cap\bbQ )$, but, in
order to get into first order model theory, we will view $\mu$ as
a collection of unary predicates on $\bfa ,\{M_r\}_{r\in
[0,1]\cap\bbQ}$, where $M_r(x) \Leftrightarrow\mu (x)=r$.

\begin{prop}\label{amalgamation for boolean algebras} $\kkk$ is a \fraisse class.\end{prop}

\begin{proof}
Obviously $\kkk$ has the HP.  Now as all structures in $\kkk$
have a common substructure (the trivial boolean algebra \{0,1\}) it
is enough to verify AP, from which JEP will follow.

So suppose $f:(\bfa ,\mu )\rightarrow (\bfb ,\nu )$ and $g:(\bfa
,\mu )\rightarrow (\bfc ,\rho )$ are embeddings.  Let $a_1,\dots
,a_n$ be the atoms of $\bfa$, $b_1,\dots , b_m$ the atoms of
$\bfb$, $c_1,\dots ,c_k$ the atoms of $\bfc$ and let
$\Gamma_1\sqcup\dots\sqcup\Gamma_n$ partition $\{1,\dots ,m\}
,\Lambda_1 \sqcup\dots\sqcup\Lambda_n$ partition $\{1,\dots ,k\}$
such that $f(a_l ) =\bigvee_{i\in\Gamma_l}b_i$ and $g(a_l
)=\bigvee_{j\in\Lambda_l}c_j$. Let $\bfd$ be the boolean algebra
with formal atoms $b_i\otimes c_j$ for
$(i,j)\in\Gamma_l\times\Lambda_l ,l\leq n$, and let the embeddings
$e:\bfb\rightarrow\bfd , h:\bfc\rightarrow\bfd$ be defined by
\begin{align*}
e(b_i)&=\bigvee_{j\in\Lambda_l}b_i\otimes c_j,\text{ where }i\in
\Gamma_l ,\\
h(c_j)&=\bigvee_{i\in\Gamma_l}b_i\otimes c_j, \text{ where
}j\in\Lambda_l .
\end{align*}
This is the usual amalgamation of boolean algebras and we now only
have to check that we can define an appropriate measure $\sigma$
on $\bfd$ to finish the proof.

Notice first that $\mu (a_l )=\sum_{i\in\Gamma_l}\nu (b_i)=
\sum_{j\in\Lambda_l}\rho (c_j)$, so the appropriate measure is
$\sigma (b_i\otimes c_j)=\tfrac{\nu (b_i)\rho (c_j)}{\mu (a_l )}$,
for $(i,j)\in\Gamma_l\times\Lambda_l$.  Then $e$ and $h$ are indeed
embeddings, as, e.g.,
\begin{equation}
\sigma (e(b_i))=\sigma (\bigvee_{j\in\Gamma_l} b_i\otimes
c_j)=\sum_{j\in\Gamma_l}\tfrac{\nu (b_i)\rho (c_j)} {\mu (a_l )}=\nu
(b_i).
\end{equation}
\end{proof}

We can now identify the \fraisse limit of the class $\kkk$.  It
follows from \fraisse's construction that the limit must be the
countable atomless boolean algebra equipped with some finitely
additive probability measure. Let $\bff$ be the boolean algebra
generated by the rational intervals in the measure algebra of
$([0,1],\lambda )$, where $\lambda$ is Lebesgue measure, together
with the restriction of Lebesgue measure.

\begin{prop}\label{measure} $(\bff ,\lambda )$ is the \fraisse
limit of $\kkk$.\end{prop}

\begin{proof}  The age of $(\bff ,\lambda )$ is clearly $\kkk$, so we need
only check that $(\bff ,\lambda )$ has the extension property.  So
suppose $(\bfa ,\mu )\subseteq (\bfb ,\nu )\in\kkk$ and
$f:\bfa\rightarrow\bff$ is an embedding of $(\bfa ,\mu )$ into
$(\bff ,\lambda )$.  Then we can  extend $f$ to an embedding
$\tilde f$ of $(\bfb ,\nu )$ into $(\bff ,\lambda )$ as follows:

An atom $a$ of $\bfa$ corresponds by $f$ to a finite disjoint union
of rational intervals in [0,1] and is also the join of atoms
$b_1,\dots ,b_k$ in $\bfb$ with rational measure.  So by
appropriately decomposing these rational intervals into finitely
many pieces we can find images for each of $b_1,\dots ,b_k$ of the
same measure.\end{proof}

If we let $\kkk =\mba_{\bbQ_2}$ be the class of finite boolean
algebras with a measure that takes positive dyadic rational values
on each non-0 element, then one can show (though it is more complicated than the preceding proof) that $\kkk$ is a \fraisse class
with limit (clop$(2^\bbN ), \sigma$), where $\sigma$ is the usual product measure on $2^\bbN$.  Moreover, ${\rm Aut}({\rm clop}(2^\bbN
),\sigma)$ is (isomorphic to) the group of measure preserving homeomorphisms of $2^\bbN$ (with the uniform convergence
topology), ${\it H}(2^\bbN, \sigma)$.

We now have:

\begin{thm}\label{density in measure spaces} There is a dense conjugacy class in {\rm
Aut}$(\bff , \lambda )$.\end{thm}

\begin{proof}  We verify that $\kkk_p$ has the JEP.  So suppose $\langle\bfa,
\mu, \psi :\bfb\rightarrow\bfc\rangle$ and $\langle\bfd ,\nu ,\phi
: \bfe\rightarrow\bfh\rangle$ are given, where $\bfa\supseteq\bfb,
\bfc$ and $\bfd\supseteq\bfe ,\bfh$ are finite boolean algebras
with rational valued probability measures $\mu ,\nu$ and $\psi$
and $\phi$ are isomorphisms preserving the measure.  We amalgamate
$\langle\bfa ,\mu \rangle$ and $\langle\bfd ,\nu\rangle$ over the
trivial boolean algebra as in the proof of Proposition \ref{amalgamation for boolean algebras}, so
our atoms in the new algebra $\bfj$ are $a_i\otimes d_j$, where
$a_1,\dots ,a_n$ and $d_1,\dots , d_m$ are the atoms at $\bfa$ and
$\bfd$, respectively, and $i\leq n, j\leq m$.  Find partitions
$$
\Gamma_1\sqcup\dots\sqcup\Gamma_k =\Lambda_1\sqcup\dots
\sqcup\Lambda_k=\{1,\dots ,n\}
$$
and
$$
\Delta_1\sqcup\dots\sqcup\Delta_l
=\Theta_1\sqcup\dots\sqcup\Theta_l =\{1,\dots ,m\}
$$
such that $\bigvee_{\Gamma_e}a_i$, $\bigvee_{\Lambda_e}a_i$,
$\bigvee_{\Delta_e}d_j$, $\bigvee_{\Theta_e}d_j$ are the atoms of
$\bfb ,\bfc ,\bfe , \bfh$ respectively and
$$
\psi (\bigvee_{\Gamma_e} a_i )=\bigvee_{\Lambda_e}a_i,\;\;\;\;\;\;
\phi (\bigvee_{\Delta_e}d_j )= \bigvee_{\Theta_e}d_j.
$$
Then we let $\bfm$ be the subalgebra of $\bfj$ generated by the
atoms $\bigvee_{\Gamma_e\times\Delta_f}a_i\otimes d_j,\; \bfn$ be
the subalgebra generated by the atoms $\bigvee_{\Lambda_e
\times\Theta_f}a_i\otimes d_j$ and $\chi :\bfm\rightarrow\bfn$ the
isomorphism given by
$$
\chi \Big(\bigvee_{\Gamma_e\times\Delta_f}a_i\otimes
d_j\Big)=\bigvee_{\Lambda_e \times\Theta_f}a_i\otimes d_j .
$$
Then $\langle\bfj ,\chi :\bfm\rightarrow\bfn\rangle$ clearly
amalgamates $\langle\bfa ,\psi :\bfb\rightarrow\bfc\rangle$ and
$\langle\bfd ,\phi :\bfe\rightarrow\bfh\rangle$ over the trivial
boolean algebra, so the only thing we need to check is that $\chi$
preserves the measure $\sigma$ on $\bfj$ given by $\sigma
(a_i\otimes d_j)=\mu (a_i)\cdot\nu (d_j)$:
\begin{align*}
&\sigma \Big(\bigvee_{\Gamma_e\times\Delta_f}a_i\otimes d_j\Big)=\sum_{i\in
\Gamma_e}\sum_{j\in\Delta_f}\mu (a_i)\nu (d_j)\\
&=\Big(\sum_{i\in\Gamma_e}\mu (a_i)\Big )\cdot\Big
(\sum_{j\in \Delta_f}\nu (d_j)\Big )=\Big
(\sum_{i\in\Lambda_e}\mu (a_i)\Big )
\cdot\Big (\sum_{j\in\Theta_f}\nu (d_j)\Big )\\
&=\sum_{i\in \Lambda_e}\sum_{j\in\Theta_f}\mu (a_i)\nu
(d_j)=\sigma\Big (\bigvee_{\Lambda_e\times\Theta_f}a_i\otimes
d_j\Big)
\end{align*}
\end{proof}
By using $\mba_{\bbQ_2}$ instead of $\mba_{\bbQ}$, we obtain the
same result for ${\it H}(2^\mathbb{N}, \sigma)$.

Let now $\kkk =\bbb\aaa $ be the class of all finite Boolean algebras. It
is well-known that $\kkk$ is a \fraisse class (the argument is
essentially that of Proposition \ref{amalgamation for boolean algebras}, forgetting about the
measures), whose \fraisse limit is the countable atomless Boolean
algebra $\bfb_\infty$.  As in the proof of \ref{density in measure spaces}, we can easily
verify the JEP for $\kkk_p$, so we have:

\begin{thm} There is a dense conjugacy class in ${\rm
Aut}(\bfb_\infty )$.\end{thm}

From this we immediately deduce:

\begin{cor}\label{dense cantor}{\rm (Glasner-Weiss \cite{glasner}, Akin-Hurley-Kennedy
\cite{akin})} In the uniform topology there is a dense conjugacy class on
${\it H}(2^\bbN )$ (the group of homeomorphisms of the Cantor
space).\end{cor}

Note that Corollary \ref{dense cantor} has as a consequence the known fact that the aperiodic homeomorphisms form a dense $G_\delta$ in $H(\ca)$. This is because for each $n>0$, $\{h\in H(\ca)\colon  \exists x\; h^n(x)=x\}$ is closed and conjugacy invariant, so it has empty interior. Thus, ${\rm APER}=\{h\in H(\ca)\colon \forall x\; \forall n>0\; h^n(x)\neq x\}$ is dense $G_\delta$.

Consider now the group ${\rm Aut}(I,\lambda )$ of measure preserving
automorphisms of the unit interval $I$ with Lebesgue measure
$\lambda$.  By a theorem of Sikorski (see Kechris \cite{kechris1}) this is
the same as the group of automorphisms of the measure algebra of
$(I,\lambda )$.

Every element $\varphi\in$ Aut$(\bff ,\lambda )$ induces a unique
automorphism $\varphi^*$ of the measure algebra, as $\bff$ is
dense in the latter.  And therefore $\varphi^*$ can be seen as an
element of Aut$(I, \lambda )$.

It is not hard to see that the mapping $\varphi\mapsto\varphi^*$
is a continuous injective homomorphism from Aut$(\bff ,\lambda )$
into Aut$(I, \lambda )$.  Moreover, the image of Aut$(\bff
,\lambda )$ is dense in Aut$(I,\lambda )$.  To see this, use for
example the following fact (see Halmos \cite{halmos}): the set of $f\in$
Aut$(I,\lambda )$, such that for some $n$, $f$ is just shuffling
the dyadic intervals at length $2^n ,[\tfrac{i}
{2^n},\tfrac{i+1}{2^n}]$ (linearly on each interval), is dense in
Aut$(I, \lambda )$.  Obviously any such $f$ is the image of some
element of Aut$(\bff ,\lambda )$.

\begin{cor} {\rm Aut}$(I,\lambda )$ has a dense
conjugacy class.\end{cor}

This is of course a weaker version of the conjugacy lemma of
ergodic theory which asserts that the conjugacy class of any
aperiodic transformation is dense in Aut$(I,\lambda )$.  But the
above proof may be of some interest as it avoids the use of
Rokhlin's Lemma.

Finally, let $\kkk$ be the \fraisse class of all finite linear
orderings. Trivially $\kkk_p$ has the JEP, so the automorphism group
of $(\bbQ ,<)$ has a dense conjugacy class.  But Aut$(\bbQ ,<)$
can be densely and continuously embedded into ${\it H}_+([0,1])$, the
group of order preserving homeomorphisms of the unit interval with
the uniform topology.  So we have the following:

\begin{cor} {\rm (Glasner-Weiss \cite{glasner})}  There is a dense
conjugacy class in ${\it H}_+([0,1])$.\end{cor}

Obviously there cannot be any dense conjugacy class in ${\it H}
([0,1])$, as it has a proper clopen normal subgroup, ${\it H}_+([0,1])$.

In the same manner one can easily check that for $\kkk =$ the class
of finite graphs, finite hypergraphs, finite posets, etc., $\kkk_p$
has the JEP.  So the automorphism group of the \fraisse limit,
e.g., the random graph, has a dense conjugacy class.
\medskip

{\bf Remark}.  For an example of a \fraisse class $\kkk$ for which
$\kkk_p$ does not have the JEP, consider the class of finite
equivalence relations with at most two equivalence classes.  Its
\fraisse limit is $(\bbN ,E)$ where $nEm\Leftrightarrow$
parity$(n)=$ parity$(m)$.  So any $f\in$ Aut$(\bbN ,E)$ either
fixes each class setwise or permutes the two classes.  Therefore
the group of $f$ setwise fixing the classes is a clopen normal
subgroup (of index 2) and there cannot be any dense conjugacy
class in Aut$(\bbN ,E)$.  So the JEP does not hold for $\kkk_p$ and
the counterexample is of course two equivalence relations with two
classes each and two automorphisms, one of which fixes the two
classes and the other switches them.

\medskip
In certain situations one can obtain more precise information
concerning dense conjugacy classes, which also has further
interesting consequences.

Suppose $G$ is a Polish group. We say that $G$ has {\it cyclically
dense conjugacy classes} if there are $g,h \in G$ such that
$\{g^nhg^{-n}\}_{n\in \Z}$ is dense in $G$.

Notice that the set
$$
D=\{(g,h)\in G^2\colon \{g^nhg^{-n}\}_{n\in \Z} \textrm{ is dense in } G\}
$$
is $G_\delta$. Also, if some section $D_g=\{h\in G\colon \{g^nhg^{-n}\}_{n\in \Z} \text{ is dense in } G\}$ is non-empty, then it is dense in $G$. Moreover, as $D_{fgf\inv}=fD_gf\inv$, the set of $g\in G$ for which $D_g\neq \tom$ is conjugacy invariant. So if there is some $g$ with a dense conjugacy class such that $D_g\neq \tom$, then there is a dense set of $f\in G$ for which $D_f$ is dense and hence $D$ itself is dense and thus comeager in $G^2$.

For each Polish group $G$, we denote by $n(G)$ the smallest number
of {\it topological generators} of $G$, i.e., the smallest $1\leq
n\leq \infty$ such that there is an $n$-generated dense subgroup
of $G$. Thus $n(G) =1$ iff $G$ is monothetic. It follows
immediately that if $G$ admits a cyclically dense conjugacy class,
then $n(G) \leq 2$.

Consider now a \fraisse structure $\bfk$ and suppose we can find
$g\in {\rm Aut}(\bfk)$ with the following property:

\medskip
(*) For any finite $\bfa,\bfb \subseteq \bfk$ and any isomorphisms
$\varphi : \bfc \rightarrow \bfd, \psi : \bfe\rightarrow\bff$,
where $\bfc, \bfd \subseteq\bfa ; \bfe,\bff\subseteq\bfb$, there
are $m,n \in \bbZ$ such that $g^m\varphi g^{-m}$, $g^n\psi g^{-n}$
have a common extension.

\medskip
Then, as in the proof of Theorem 1.1, we see that there is
$h\in{\rm Aut}(\bfk)$ with $\{g^nhg^{-n}\}$ dense in Aut$(\bfk)$,
i.e., $\bfk$ admits a cyclically dense conjugacy class, and
moreover Aut$(\bfk)$ is topologically 2-generated.

Macpherson \cite{macpherson1} 3.3 has shown that the automorphism group of the
random graph has cyclically dense conjugacy classes and so is
topologically 2-generated. We note below that
Aut$(\bfb_{\infty})$, and thus ${\it H}(2^{\bbN})$, as well as
${\it H}(2^{\bbN}, \sigma)$, Aut$(X,\mu)$, Aut$(\bbN^{<\bbN})$ admit
cyclically dense conjugacy classes, so in particular, they are
topologically 2-generated. That Aut$(X,\mu)$ is topologically
2-generated was earlier proved by different means in Grzaslewicz
\cite{grzaslewicz} and Prasad \cite{prasad}.

\begin{thm} Each of the groups ${\it H}(2^{\bbN})$,
${\it H}(2^{\bbN},\sigma)$, ${\rm Aut}(X,\mu)$, and ${\rm Aut}(\bbN^{<\bbN})$ has
cyclically dense conjugacy classes and  is topologically
2-generated. \end{thm}

\begin{proof} We will sketch the proof for ${\it H}(2^{\bbN})$ or rather its
isomorphic copy Aut$(\bfb_{\infty})$. The proofs in the cases of
${\it H}(2^{\bbN},\sigma)$ and Aut$(\bbN^{<\bbN})$ are similar. Since
${\it H}(2^{\bbN},\sigma)$ can be densely embedded in Aut$(X,\mu)$, the
result follows for this group as well.

It is clear (see the proof of \ref{density in measure spaces}) that it is enough to show the
following, in order to verify property (*) for $\bfk =
\bfb_{\infty}$: There is an automorphism $g$ such that given any
finite subalgebras $\bfa,\bfb$, there is $n\in \bbZ$ with $g^n
(\bfa), \bfb$ independent (i.e., any non-zero elements of $g^n
(\bfa), \bfb$ have non-zero join).

To find such a $g$ view $\bfb_{\infty}$ as the algebra of clopen
subsets of $2^{\bbZ}$ and take $g$ to be the Bernoulli shift on this
space.\end{proof}

It is not hard to show that also $U(\ell_2)$ has a cyclicaly dense conjugacy class.

Solecki has recently shown that the group of isometries of the
Urysohn space admits a cyclically dense conjugacy class and thus
it is topologically 2-generated.
\medskip

{\bf Remarks.} Prasad \cite{prasad} showed that there is a comeager set of
pairs $(g,h)$ in Aut$(X,\mu)^2$ generating a dense subgroup of
Aut$(X,\mu)$ (moreover, using the results of section \ref{generic freeness} one sees easily that generically the subgroup generated is free). This also follows from our earlier remarks, since
the $g$ that witnesses the cyclically dense conjugacy class is in
this case the shift, whose conjugacy class is dense (by the conjugacy lemma of ergodic theory).

\medskip
We can also consider the diagonal conjugacy action of Aut$(\bfk
)$ on Aut$(\bfk)^n$, for $n = 1,2, \dots, \bbN$. Note that there
is a dense diagonal conjugacy class in Aut$(\bfk)^\bbN$ iff for each $n =
1,2 \dots ,$ there is a dense diagonal conjugacy class in Aut$(\bfk)^n$.
This follows from the fact that in the latter case, the set of elements $(g_m)\in {\rm Aut}(\bfk)^\bbN$, whose diagonal conjugacy class is dense in Aut$(\bfk)^n$, is dense $G_\delta$ for all $n$, and hence we can pick a $(g_m)$ such that it holds for all $n$ and therefore also for $\N$.

\medskip
Consider a \fraisse class $\kkk$ with \fraisse limit $\bfk$ and let $G={\rm Aut}(\bfk)$. Define for each $n\geq 1$ the following sets (which are all $G_\delta$).
\begin{displaymath}\begin{split}
F_n&=\{(f_1,\ldots,f_n)\in G^n\colon \forall x\in \bfk\; x\text{'s orbit under } \langle f_1,\ldots,f_n\rangle \text{ is finite}\} \\
&=\{(f_1,\ldots,f_n)\in G^n\colon \langle f_1,\ldots,f_n\rangle \text{ is relatively compact in } G\}.\\
D_n&=\{(f_1,\ldots,f_n)\in G^n\colon \langle f_1,\ldots,f_n\rangle \text{ is dense in } G\}.\\
E_n&=\{(f_1,\ldots,f_n)\in G^n\colon \langle f_1,\ldots,f_n\rangle \text{ is non-discrete in } G\}.\\
H_n&=\{(f_1,\ldots,f_n)\in G^n\colon \langle f_1,\ldots,f_n\rangle \text{ freely generates a free subgroup of } G\}.
\end{split}\end{displaymath}
Since the sets $F_n$, $D_n$, $E_n$ and $H_n$ are $G_\delta$ sets invariant under the diagonal conjugacy action of $G$ on $G^n$, it follows that if $G^n$ has a dense diagonal conjugacy class then each of them is either comeager  or nowhere dense in $G^n$. We also see that $D_n$ is never dense unless $G=\{1\}$, for if $H<G$ is a proper open subgroup and $(g_1,\ldots,g_n)\in H^n$, then $\langle g_1,\ldots,g_n\rangle\subseteq H$ and hence is not dense. So $D_n\cap H^n=\emptyset$. On the other hand, if $\kkk$ satisfies the Hrushovski property (see Definition \ref{hrushovski property}), then $F_n$ is automatically comeager. Moreover, unless $G$ itself is compact, the sets $F_n$ and $D_n$ are of course disjoint and thus if also $\kkk$ has the Hrushovski property, $D_n$ is nowhere dense.

These results can also be combined with the comments in section \ref{generic freeness}, where we will study the sets $H_n$.

Suppose $\kkk$ is a \fraisse class, and $\bfk$ its \fraisse limit.
We would like to characterize as before when Aut$(\bfk)^n$ has a
dense diagonal conjugacy class. For this we introduce the class of
$n$-systems $\kkk^n_p$ for each $n\geq 1$. An $n$-{\it system} $\s
=\langle\bfa ,\psi_1:\bfb_1\rightarrow\bfc_1, \dots
,\psi_n:\bfb_n\rightarrow\bfc_n\rangle$ consists of finite
structures $\bfa ,\bfb_i,\bfc_i\in\kkk$ with
$\bfb_i,\bfc_i\subseteq\bfa$ and $\psi_i$ an isomorphism of
$\bfb_i$ and $\bfc_i$.  As before, an embedding of one $n$-system
$\s =\langle\bfa ,\psi_1:\bfb_1\rightarrow\bfc_1,\dots ,\psi_n:
\bfb_n\rightarrow\bfc_n\rangle$ into another $\ttt =\langle\bfd
,\phi_1: \bfe_1\rightarrow\bff_1,\dots
,\phi_n:\bfe_n\rightarrow\bff_n\rangle$ is a function
$f:\bfa\rightarrow\bfd$ embedding $\bfa$ into $\bfd ,\bfb_i$ into
$\bfe_i$ and $\bfc_i$ into $\bff_i$ such that $f\circ\psi_i
\subseteq \phi_i\circ f, i=1,\dots ,n$.  So we can talk about JEP,
WAP, etc., for $n$-systems as well.

We can prove exactly as before:

\begin{thm} Let $\kkk$ be a \fraisse class and $\bfk$
its \fraisse limit.  Then the following are equivalent:

(i) There is a dense diagonal conjugacy class in ${\rm Aut}(\bfk )^n$,

(ii) $\kkk^n_p$ has the {\rm JEP}.\end{thm}

From this we easily have the following:

\begin{thm} There is a dense diagonal conjugacy class in each of
$${\it H}(2^{\bbN})^{\bbN}, {\it H}(2^{\bbN},\sigma)^{\bbN}, {\rm
Aut}(X,\mu)^{\bbN}, {\rm Aut}(\bbN^{<\bbN})^{\bbN}, {\rm
Aut}(\bfuo )^{\bbN}, {\rm Iso}(\bfu)^{\bbN}.$$\end{thm}

\begin{defi} Let $\kkk$ be a \fraisse class and let $\kkk_p$ be the
corresponding class of systems. We say that $\kkk_p$ satisfies the
{\em cofinal joint embedding property} (CJEP) if for each
$\bfa\in \kkk$ there is $\bfa\leq\bfb\in\kkk$ such that for
$\jjj_\bfb=\langle\bfb,id\colon\bfb\rightarrow\bfb\rangle$, any
systems $\ttt_0,\ttt_1$ and embeddings $e\colon \jjj_\bfb\rightarrow
\ttt_0$, $f\colon \jjj_\bfb\rightarrow\ttt_1$, there is a system $\rrr$
and embeddings $g\colon\ttt_0\rightarrow \rrr$,
$h\colon\ttt_1\rightarrow \rrr$ such that $g\circ e=h\circ f$.
\end{defi}

This property is enough to ensure that if $\bfk$ is the \fraisse
limit of $\kkk$, then ${\rm Aut}(\bfk)$ has a neighborhood basis at
the identity consisting of clopen subgroups with a dense conjugacy
class.

\begin{thm}\label{CJEP} Let $\kkk$ be a \fraisse class with
\fraisse limit $\bfk$ and suppose that $\kkk_p$ satisfies the {\rm
CJEP}. Then for any finite substructure $\bfa\subseteq \bfk$,
there is a finite substructure $\bfa \subseteq\bfb\subseteq\bfk$
such that ${\rm Aut}(\bfk)_{(\bfb)}=\{g\in {\rm Aut}(\bfk)\colon
g\restriction_\bfb=id_\bfb\}$ has a dense conjugacy
class.\end{thm}

\begin{proof} Let $\bfa\subseteq\bfk$ be given and find by the CJEP some
$\bfa\leq\bfb\in\kkk$ satisfying the conditions of the definition.
By ultrahomogeneity of $\bfk$, we can suppose that
$\bfa\subseteq\bfb\subseteq\bfk$. Let $\ov b=\langle
b_1,\ldots,b_n\rangle$ be new names for the elements of $\bfb$ and
let $(\bfk,\ov b)$ be the corresponding expansion of $\bfk$.
Notice that ${\rm Aut}(\bfk,\ov b)={\rm Aut}(\bf K)_{(\bfb)}$.
Similarly, we let $(\kkk,\ov b)$ be the class of expanded structures
$(\bfd,\ov b)$, where $\bfb\leq \bfd\in \kkk$ and $b_1,\ldots,b_n$
are the names for some fixed copy of $\bfb$ in $\bfd$. We claim
that (i) $(\kkk,\ov b)$ is a \fraisse class and (ii) $(\kkk,\ov b)_p$
has the JEP.

First (i): The HP for $(\kkk,\ov b)$ follows at once from the HP of
$\kkk$. Now, if $(\bfc,\ov b)$, $(\bfd,\ov b)$,$(\bfe, \ov b)\in
(\kkk,\ov b)$ and $e,f$ are embeddings of $(\bfc, \ov b)$ into
$(\bfd, \ov b)$ and $(\bfc, \ov b)$ into $(\bfe ,\ov b)$ resp.,
then $e,f$ embed $\bfc$ into $\bfd$ and $\bfc$ into $\bfe$ resp.
So by the AP for $\kkk$ there is a structure $\bff\in\kkk$ and
embeddings $g\colon\bfd\rightarrow \bff$,
$h\colon\bfe\rightarrow\bff$ such that $g\circ e=h\circ f$. In
particular, $g\circ e(b_i^\bfc)=h\circ f(b_i^\bfc)$ for
$i=1,\ldots,n$ and so we can expand $\bff$ by setting
$b_i^\bff=g\circ e(b_i^\bfc)$. Then $g$ and $h$ embed $(\bfd,\ov
b)$ and $(\bfe,\ov b)$ into $(\bff,\ov b)$ and $g\circ e=h\circ
f$. So $(\kkk,\ov b)$ has the AP. Finally, the JEP follows for
$(\kkk,\ov b)$ from the AP, since all structures have a common
substructure, namely the one generated by $b_1,\dots,b_n$.

Now for (ii): Suppose $\s,\ttt\in (\kkk, \ov b)_p$. Then for some
$\bfd, \bfe\subseteq\bfc$, $\bfg,\bfh\subseteq\bff$,
$$
\s=\langle (\bfc,\ov b), \psi\colon(\bfd,\ov
b)\rightarrow(\bfe,\ov b)\rangle
$$
and
$$
\ttt=\langle (\bff,\ov b), \phi\colon(\bfg,\ov
b)\rightarrow(\bfh,\ov b)\rangle.
$$
But as $\psi$ and $\phi$ both preserve $b_1,\ldots,b_n$ this means
that $\jjj_\bfb$ embeds into both $\tilde\s=\langle \bfc,
\psi\colon\bfd\rightarrow\bfe\rangle$ and $\tilde\ttt=\langle \bff,
\phi\colon\bfg\rightarrow\bfh\rangle$. Therefore there is an
amalgamation of $\tilde\s$ and $\tilde\ttt$, which can be expanded
to a common super-system of $\s$ and $\ttt$.

We now only need to show that the \fraisse limit of $(\kkk,\ov b)$
is $(\bfk,\ov b)$. But for this it is enough to notice that
$(\bfk,\ov b)$ is still ultrahomogeneous as $\bfk$ is, and ${\rm
Age}(\bfk,\ov b)=(\kkk,\ov b)$ as ${\rm Age}(\bfk)=\kkk$ and $\bfk$
is ultrahomogeneous. So by \fraisse's theorem on the uniqueness of
the \fraisse limit we have the result. For now we can just apply
Theorem \ref{criterion for density} to $(\kkk,\ov b)$ with limit
$(\bfk,\ov b)$.\end{proof}
\medskip

In any natural instance it is certainly easy to verify the CJEP,
for example, it is true in any of the cases considered in this
paper. However, we have no general theorem saying that it should
follow from the existence of a comeager conjugacy class, and
probably it does not.

\section{Generic automorphisms}\label{comeagre}

We will now consider the question of when Aut$(\bfk )$ has a
comeager conjugacy class, when $\bfk$ is the \fraisse limit of
some \fraisse class $\kkk$.

Now it is known, by a theorem due to Effros, Marker and Sami (see,
e.g., Becker-Kechris \cite{becker}), that any non-meager orbit of a Polish
group acting continuously on a Polish space is in fact $G_\delta$.
So we are actually looking for criteria for when there is a dense
$G_\delta$ orbit.

Recall now some basic facts about Hjorth's notion of turbulence
(see Hjorth \cite{hjorth} or Kechris \cite{kechris2}).

Suppose a Polish group $G$ acts continuously on a Polish space
$X$.  A point $x\in X$ is said to be {\it turbulent} if for every
open neighborhood $U$ of $x$ and every symmetric open neighborhood
$V$ of the identity $e\in G$, the {\it local orbit}, $\ooo (x,U,V)$,
is somewhere dense, i.e.,
\begin{equation}
{\rm Int} (\overline{\ooo (x,U,V)})\neq\emptyset.
\end{equation}
Here
\begin{equation}\begin{split}
\ooo (x,U,V)=&\{y\in X:\exists g_0,\dots g_k\in V\;\forall i\leq k\\
&(g_ig_{i-1}\dots g_0\cdot x\in U\text{ and }g_kg_{k-1}\dots
g_0\cdot x=y)\}.
\end{split}
\end{equation}
Notice that the property of being turbulent is $G$-invariant (see
Kechris \cite{kechris2} 8.3), so we can talk about {\it turbulent orbits}.
We also have the following fact (see Kechris \cite{kechris2} 8.4):

\begin{prop}\label{turbulence} Let a Polish group $G$ act
continuously on a Polish space $X$ and $x\in X$.  The following
are equivalent:

(i) $x$ is turbulent.

(ii) For every open neighborhood $U$ of $x$ and every symmetric
open neighborhood $V$~of $e\in G, x\in {\rm Int}(\overline{\ooo
(x,U,V)})$.\end{prop}

Notice that if $V$ is an open subgroup of $G$, then $\ooo (x,U,V)=U
\cap V\cdot x$.  So if $G$ is a closed subgroup of $S_\infty$,
then, as $G$ has an open neighborhood basis at the identity
consisting of open subgroups, we see that $x\in X$ is turbulent
iff $x \in {\rm Int} (\overline{V \cdot x})$ for all open
subgroups $V \leq G$. Therefore, the following provide equivalent
conditions for turbulence:

\begin{prop}\label{conditions turbulence} Let $G$ be a closed subgroup of $S_\infty$
acting continuously on a Polish space $X$ and let $x\in X$.  Then
the following are equivalent:

(i) The orbit $G\cdot x$ is non-meager.

(ii) For each open subgroup $V\leq G, V\cdot x$ is non-meager.

(iii) For each open subgroup $V\leq G, V\cdot x$ is somewhere dense.

(iv) For each open subgroup $V\leq G,x\in{\rm Int} (\overline{V\cdot x})$.

(v) The point $x$ is turbulent.\end{prop}

\begin{proof}  (i) $\Rightarrow$ (ii): Suppose $G\cdot x$ is non-meager and
$V\leq G$ is an open subgroup of $G$.  Then we can find $g_n\in G$
such that $G=\bigcup g_nV$, so some $g_nV\cdot x$ is non-meager
and therefore $V\cdot x$ is non-meager.

(ii) $\Rightarrow$ (iii):  Trivial.

(iii) $\Rightarrow$ (iv):  Suppose $V$ is an open subgroup of $G$
such that $V\cdot x$ is dense in some open set $U\neq\emptyset$.
Take $g\in V$ such that $g\cdot x\in U$.  Then $g^{-1}V\cdot
x=V\cdot x$ is dense in $g^{-1}\cdot U\ni x$ and $x\in{\rm Int}
(\overline{V\cdot x})$.

(iv) $\Rightarrow$ (i):  Suppose $F_n\subseteq X$ are closed
nowhere dense such that $G\cdot x\subseteq\bigcup_nF_n$.  Then
$K_n=\{ g\in G:g\cdot x\in F_n\}$ are closed and, as they cover
$G$, some $K_n$ has non-$\emptyset$ interior.  So suppose
$gV\subseteq K_n$ for some open subgroup $V\leq G$.  Then $gV\cdot
x\subseteq F_n$ and both $gV\cdot x$ and $V\cdot x$ are nowhere
dense.

The equivalence of (iv), (v) is clear using the remarks following
\ref{turbulence}.\end{proof}

It was asked in Truss \cite{truss} whether the existence of a generic
automorphism of the limit of a \fraisse class $\kkk$ could be
expressed in terms of AP and JEP for $\kkk_p$.  He gave a partial
answer to this showing that the existence of a generic
automorphism implied JEP for $\kkk_p$, and, moreover, if $\kkk_p$
satisfies JEP and the so-called {\it cofinal} AP (CAP), then there
is indeed a generic automorphism. Here we say that $\kkk_p$
satisfies the CAP if it has a cofinal subclass, under
embeddability, that has the AP.

This is clearly equivalent to saying that there is a subclass
$\llll\subseteq \kkk_p$, cofinal under embeddability, such that for any
$\s\in\llll , \ttt,\rrr\in\kkk_p$ and embeddings $f:\s\rightarrow\ttt
,g:\s\rightarrow\rrr$, there is $\h\in\kkk_p$ and embeddings
$h:\ttt\rightarrow\h ,i:\rrr\rightarrow\h$ such that $h\circ f=i\circ
g$.  We mention that for any \fraisse class $\kkk$, the class $\llll$
of $\s =\langle\bfa ,\psi :\bfb\rightarrow\bfc\rangle\in \kkk_p$,
such that $\bfa$ is generated by $\bfb$ and $\bfc$, is cofinal
under embeddability.  This follows from the ultrahomogeneity of
the \fraisse limit of $\kkk$.

\begin{defi} A class $\kkk$ of finite structures
satisfies the weak amalgamation property {\rm (WAP)} if for every
$\bfa\in\kkk$ there is a $\bfb\in\kkk$ and an embedding
$e:\bfa\rightarrow\bfb$ such that for all $\bfc ,\bfd\in\kkk$ and
embeddings $i:\bfb\rightarrow\bfc ,j:\bfb \rightarrow\bfd$ there
is an $\bfe\in\kkk$ and embeddings $l :\bfc\rightarrow \bfe
,k:\bfd\rightarrow\bfe$ amalgamating $\bfc$ and $\bfd$ over
$\bfa$, i.e., such that $l\circ i\circ e=k\circ j\circ
e$.\end{defi}

A similar definition applies to the class $\kkk_p$.

We can now use turbulence concepts to provide the following answer
to Truss' question.  We have recently found out that Ivanov \cite{ivanov}
had also proved a similar result, by other techniques, in a
somewhat different context, namely that of $\aleph_0$-categorical
structures.

\begin{thm}\label{criterion comeager} Let $\kkk$ be a \fraisse class and $\bfk$
its \fraisse limit.  Then the following are equivalent:

(i) $\bfk$ has a generic automorphism.

(ii) $\kkk_p$ satisfies the {\rm WAP} and the {\rm JEP}.\end{thm}

\begin{proof} (i) $\Rightarrow$ (ii):  We know that if $\bfk$ has a generic
automorphism, then some $f\in$ Aut$(\bfk )$ is turbulent and has
dense conjugacy class.  So by Theorem 1.1, $\kkk_p$ has the JEP.

Given now $\s =\langle\bfa ,\psi
:\bfb\rightarrow\bfc\rangle\in\kkk_p$, we can assume that
$\bfa\subseteq\bfk$.  Moreover, as $[\psi :\bfb\rightarrow\bfc ]$
is open nonempty and $f$ has dense conjugacy class, we can suppose
that $\psi\subseteq f$.  Let $V_\bfa =[id:\bfa\rightarrow\bfa ]$,
which is a clopen subgroup of Aut$(\bfk )$.  By the turbulence of
$f$, we know that $c(V_\bfa ,f)=\{g^{-1}fg: g\in V_\bfa\}$ is
dense in some open neighborhood $[\phi :\bfe\rightarrow\bff ]$ of
$f$.

Now, since $f\supseteq\psi$, we can suppose that
$\psi\subseteq\phi$. Let $\bfd$ be the substructure of $\bfk$
generated by $\bfa ,\bfe$ and $\bff$, and put $\ttt =\langle\bfd
,\phi :\bfe\rightarrow\bff\rangle$.  So $\s$ is a subsystem of
$\ttt$ and denote by $e$ the inclusion mapping of $\s$ into $\ttt$.

Assume now that $i:\ttt\rightarrow\f =\langle\bfh ,\chi
:\bfm\rightarrow\bfn \rangle$ and $j:\ttt\rightarrow\g =\langle\bfp
,\xi :\bfq\rightarrow\bfr \rangle$ are embeddings.  Then we wish
to amalgamate $\f$ and $\g$ over $\s$.

By using the extension property of $\bfk$ we can actually assume
that $\bfh$ and $\bfp$ are substructures of $\bfk$ and
$\bfd\subseteq\bfh ,\bfp$, with $\chi$ and $\xi$ extending $\phi$,
and $i,j$ the inclusion mappings. By the density of $c(V_\bfa ,f)$
in $[\phi :\bfe\rightarrow\bff ]$ there are $g,k\in V_\bfa =[id
:\bfa\rightarrow\bfa ]$ such that $g^{-1}fg\supseteq\chi$ and
$h^{-1}fh\supseteq\xi$.  Let $\bfs$ be the substructure of $\bfk$
generated by $g''\bfh$ and $h''\bfp ,\ \bft$ the substructure
generated by $g''\bfm$ and $h''\bfq$ and $\bfu$ the
substructure generated by $g''\bfn$ and $h'' \bfr$.  Finally put
$\theta =f\restriction_\bft :\bft\rightarrow\bfu$ and $\e
=\langle\bfs ,\theta :\bft\rightarrow\bfu\rangle$.

As $f$ restricts to an isomorphism of $g''\bfm$ with $g''\bfn$ and
an isomorphism of $h''\bfq$ with $h''\bfr$, it also restricts to
an isomorphism of $\bft$ with $\bfu$.  So $\theta$ is well
defined.  Moreover, $g$ and $h$ obviously embed $\f$ and $\g$ into
$\e$, as $g^{-1}fg\supseteq\chi$ and $h^{-1}fh\supseteq\xi$.  And
finally, as $g,h\in [id:\bfa\rightarrow\bfa ]$, we have that
$g\circ i\circ e=h\circ j\circ e=id_\bfa$.  So $\e$ is indeed an
amalgam over $\s$.  Therefore $\kkk_p$ satisfies WAP.

(ii) $\Rightarrow$ (i):  Now suppose $\kkk_p$ satisfies WAP and JEP.
We will construct an $f\in$ Aut$(\bfk )$ with turbulent and dense
conjugacy class. Notice that this will be enough to insure that
the conjugacy class of $f$ is comeager. For as the conjugacy
action is continuous, if the orbit of $f$ is non-meager, then it
will be comeager in its closure, i.e., comeager in the whole
group.

As before for $\psi :\bfb\rightarrow\bfc$, an isomorphism between
finite substructures of $\bfk$, we let
\begin{equation}
D(\psi :\bfb\rightarrow\bfc )=\{f\in {\rm Aut}(\bfk ): \exists
g\in {\rm Aut}(\bfk )(g ^{-1}fg\supseteq\psi )\}.
\end{equation}
For $\bfe$ a finite substructure of $\bfk$, we let $V_\bfe =[id:
\bfe \rightarrow\bfe ]$, which is a clopen subgroup of Aut$(\bfk
)$.

We have seen in the proof of Theorem \ref{criterion for density} that the JEP for $\kkk_p$ implies that each $D(\psi
:\bfb \rightarrow\bfc )$ is open dense in Aut$(\bfk )$.  And
moreover, any element in their intersection has dense conjugacy
class.

Now for $\psi :\bfb \rightarrow\bfc$ let $\bfa$ be the
substructure generated by $\bfb ,\bfc$ and $\s =\langle\bfa ,\psi
:\bfb\rightarrow\bfc\rangle$.  Then by WAP for $\kkk_p$ and the
extension property of $\bfk$, there is an extension $\hat\s
=\langle\hat\bfa ,\hat\psi :\hat\bfb\rightarrow\hat\bfc \rangle$
of $\s$ such that any two extensions of $\hat\s$ can be
amalgamated over $\s$.  By extending $\hat\s$, we can actually
assume that $\hat\bfa$ is generated by $\hat\bfb ,\hat\bfc$.
Enumerate all such $\hat\psi$'s as $\hat\psi_1,\hat\psi_2,\dots$.
Moreover, list all the possible finite extensions $\theta
:\bfm\rightarrow\bfn$ of $\hat\psi_m :\hat\bfb\rightarrow\hat\bfc$
as $\theta^{\hat\psi_m}_1, \theta^{\hat\psi_m}_2,\dots\;$. Define
\begin{equation}\begin{split}
E(\psi :\bfb\rightarrow\bfc )=&\{f\in {\rm Aut}(\bfk ):
f\supseteq\psi \Rightarrow \\
&\text{ (for some }\phi\supseteq\psi \text{ and some } m \text{ we
have }f\supseteq \hat\phi_m )\}
\end{split}\end{equation}
and
\begin{equation}
F_{m,n}(\psi :\bfb\rightarrow\bfc )=\{f\in{\rm Aut}(\bfk ):
f\supseteq \hat\psi_m \Rightarrow (c(V_\bfb ,f)\cap
[\theta^{\hat\psi_m}_n]\neq\emptyset )\}
\end{equation}
Now obviously both $E(\psi :\bfb\rightarrow\bfc )$ and
$F_{m,n}(\psi :\bfb \rightarrow\bfc )$ are open and $E(\psi
:\bfb\rightarrow\bfc )$ is dense.

\begin{lemme} Suppose $f\in{\rm Aut}(\bfk )$ is in $E(\psi :\bfb
\rightarrow\bfc )$ and $F_{m,n}(\psi :\bfb\rightarrow\bfc )$ for
all $\psi : \bfb\rightarrow\bfc$ and $m,n\in\bbN$.  Then $f$ is
turbulent.\end{lemme}

\begin{proof}  Let a clopen subgroup $V_\bfe \leq$ Aut$(\bfk )$ be given. We
shall show that $c(V_\bfe ,f)$ is somewhere dense.

As $f\in E(f\restriction_\bfe :\bfe\rightarrow f''\bfe )$ there are
some $\phi :\bfb\rightarrow\bfc , \bfe\subseteq\bfb$ and $m$ such
that $f\supseteq \hat\phi_m$.  We claim that $c(V_\bfb ,f)$ and a
fortiori $c(V_\bfe ,f)$ is dense in $[\hat\phi_m
:\hat\bfb\rightarrow\hat\bfc ]$.  This is because any basic open
subset of $[\hat\phi_m :\hat\bfb\rightarrow\hat\bfc ]$ is of the
form $[\theta^{\hat\phi_m}_n]$ for some $n\in\bbN$ and we know that
$c(V_\bfb ,f)\cap [\theta^{\hat\phi_m}_n]\neq\emptyset$.\end{proof}
\medskip

So now we only need to show that each $F_{m,n}(\psi
:\bfb\rightarrow\bfc )$ is dense, as any element in the
intersection of the $D,E$ and $F_{m,n}$'s will do.

Given a basic open set $[\phi :\bfa\rightarrow\bfd ]$, where we
can suppose that $\phi\supseteq\hat\psi_m$, we need to show that
for some $f\supseteq\phi$, $c(V_\bfb ,f)\cap
[\theta^{\hat\psi_m}_n]\neq\emptyset$.  Now
$\theta^{\hat\psi_m}_n$ and $\phi$ are both extensions of
$\hat\psi_m$, so by WAP for $\kkk_p$, they can be amalgamated over
$\psi$.  It follows by the extension property of $\bfk$ that there
is some $g\in$ Aut$(\bfk )$ fixing $\bfb$, such that $g^{-1}\phi
g$ and $\theta^{\hat\psi_m}_n$ are compatible in Aut$(\bfk )$, i.e.,
for some $f\supseteq\phi ,g^{-1}fg\in [\theta^{\hat\psi_m}_n ]$.

This finishes the proof.\end{proof}

{\bf Remark}.  It is easy to verify that the intersection of the
$D,E$, and $F_{m,n}$'s is conjugacy invariant, so actually this
intersection is exactly equal to the set of generic automorphisms.
\medskip

As obviously the CAP implies the WAP, we have:

\begin{cor} {\rm (Truss \cite{truss})} If $\kkk_p$ has the cofinal
amalgamation property, then there is a generic automorphism of
$\bfk$.\end{cor}

{\bf Remark}.  In Hodges \cite{hodges1} the \fraisse theory is developed in
the context of a class $\kkk$ of {\it finitely generated} (but not
necessarily finite) structures satisfying (HP), (JEP), and (AP).
It is not hard to see that all the arguments and results of \S \ref{dense},
\S \ref{comeagre} carry over without difficulty to this more general context.
\medskip

Truss \cite{truss} calls an automorphism $f\in$ Aut$(\bfk ),\bfk =$
\fraisse limit of $\kkk$, $\kkk$ a \fraisse class, {\it locally
generic} if its conjugacy class is non-meager.

Given a \fraisse class $\kkk$, let us say that $\kkk_p$ satisfies the
{\it local} WAP if there is $\s =\langle\bfa ,\psi
:\bfb\rightarrow\bfc\rangle\in \kkk_p$ such that WAP holds for the subclass $\llll_p$ of $\kkk_p$ consisting of all
$\ttt\in\kkk_p$ into which $\s$ embeds.

Using again Proposition \ref{conditions turbulence}, we have:

\begin{thm} \label{locally generic}Let $\kkk$ be a \fraisse class and $\bfk$
its \fraisse limit.  Then the following are equivalent:

(i) $\bfk$ admits a locally generic automorphism.

(ii) $\kkk_p$ satisfies the local {\rm WAP}.\end{thm}

\begin{proof}  (i) $\Rightarrow$ (ii):  Suppose $\bfk$ admits a locally
generic automorphism $f$.  Then by Proposition \ref{conditions turbulence}, the conjugacy
class of $f$ is turbulent and of course dense in some open set
$U\subseteq$ Aut$(\bfk )$. So there is some $\s =\langle\bfa ,\psi
:\bfb\rightarrow\bfc\rangle\in\kkk_p$ such that $U\supseteq [\psi
:\bfb\rightarrow\bfc ]$.  Then, as in the proof of (i)
$\Rightarrow$ (ii) in Theorem \ref{criterion comeager}, we can see that WAP holds for
all $\ttt\in\kkk_p$ into which $\s$ can be embedded.

(ii) $\Rightarrow$ (i):  Suppose $\s =\langle\bfa ,\psi
:\bfb\rightarrow \bfc\rangle\in\kkk_p$ witnesses that $\kkk_p$ has
the local WAP.  Then we repeat the construction of (ii)
$\Rightarrow$ (i) in the proof of Theorem \ref{criterion comeager},
by taking as our first approximation $\psi$ (and without of course
using the sets $D(\phi :\bfd\rightarrow\bfe )$ whose density was
ensured by JEP).\end{proof}
\medskip

Similarly one can see that the existence of a conjugacy class that
is somewhere dense is equivalent to a local form of JEP defined in
an analogous way.

We can also characterize WAP in terms of local genericity.

\begin{thm}\label{dense locally generic} Let $\kkk$ be a \fraisse class and $\bfk$
its \fraisse limit.  Then the following are equivalent:

(i) $\bfk$ admits a dense set of locally generic automorphisms.

(ii) $\bfk$ admits a comeager set of locally generic
automorphisms.

(iii) $\kkk_p$ satisfies the {\rm WAP}.\end{thm}

\begin{proof}  (i) $\Rightarrow$ (iii) and (iii) $\Rightarrow$ (i) are
proved as in Theorems \ref{criterion comeager} and \ref{locally generic}, beginning the construction of a
turbulent point from any given $\psi :\bfb\rightarrow\bfc$.

(ii) $\Rightarrow$ (i) is of course trivial and (i) $\Rightarrow$
(ii) follows from the fact that if an orbit is non-meager then it is
comeager in its closure.\end{proof}

\begin{thm} Let $\kkk$ be a \fraisse class with \fraisse
limit $\bfk$ and suppose that $\kkk_p$ satisfies the {\rm WAP} and
the {\rm CJEP}. Then ${\rm Aut}(\bfk)$ has a neighborhood basis at
the identity consisting of clopen subgroups each having a comeager
conjugacy class.\end{thm}

\begin{proof} By Theorem \ref{dense locally generic}, ${\rm Aut}(\bfk)$ has a dense set of
non-meager conjugacy classes. Moreover, by Theorem \ref{CJEP}, ${\rm
Aut}(\bfk)$ has a neighborhood basis at the identity consisting of
clopen subgroups with dense conjugacy classes. Now assume that $G$
is a clopen subgroup of ${\rm Aut}(\bfk)$ with a dense conjugacy
class and that $g\in G$ has a non-meager conjugacy class in ${\rm
Aut}(\bfk)$. Then by Proposition \ref{conditions turbulence} (ii),
$g$ has a non-meager conjugacy class in $G$. But as the set of group
elements with dense conjugacy classes forms a $G_\delta$, which,
when non-empty, is dense, the conjugacy class of $g$ intersects this
set. So $g$'s conjugacy class is both dense and non-meager in $G$,
i.e., comeager in $G$.\end{proof}

{\bf Remark.} As is easily checked, all the structures considered
in this paper that are shown to have a comeager conjugacy class
actually also satisfy the above theorem. Let us just mention that
for the class of $\omega$-stable, $\aleph_0$-categorical
structures one just needs to notice that they stay
$\omega$-stable, $\aleph_0$-categorical after having been expanded
by a finite number of constants (see Hodges \cite{hodges1}).

\section{Normal form for isomorphisms of finite subalgebras of
clop${(2^{\pmb\bbN} )}$} In this section we will develop some facts
needed in the proof of the existence of generic homeomorphisms of
$2^\bbN$, which will be given in the next section.

Recall that the \fraisse limit of the class of finite boolean
algebras is the countable atomless boolean algebra $\bfb_\infty$,
which can be concretely realized as clop$(2^\bbN )$, the boolean
algebra of clopen subsets of Cantor space $2^\bbN$.

We will eventually prove that the class $\kkk_p$, where $\kkk$ is
the class of finite boolean algebras, has the CAP. Thus, we will
first have to describe the cofinal class $\llll$ of $\kkk_p$ over
which we can amalgamate.  First of all it is clear that if
$\phi\colon\bfb_\infty\til\bfb_\infty$ is such that for some $a\in
\bfb_\infty$, $\phi(a)<a$, then the orbit of $a$ under $\phi$ is
infinite. This shows that there are partial automorphisms of
finite boolean algebras that cannot be extended to full automorphisms
of bigger finite boolean algebras. This fact makes the situation a
lot messier than in the case of measured boolean algebras and, in
order to prove the amalgamation property, requires us to be able to
describe the local structure of a partial automorphism of a finite
boolean algebra.

\begin{defi}
Let $\bfb ,\bfc$ be finite subalgebras of ${\rm clop}(2^\bbN )$ and
$\psi :\bfb\rightarrow\bfc$ an isomorphism.  A {\em refinement} of
$\psi :\bfb\rightarrow\bfc$ consists of finite superalgebras
$\bfb\subseteq\bfb' ,\bfc\subseteq\bfc'$ and an isomorphism
$\psi':\bfb'\rightarrow\bfc'$ such that $\psi'\restriction_\bfb
=\psi$.
\end{defi}

Fix an isomorphism $\psi :\bfb\rightarrow\bfc$ of finite subalgebras
of clop$(2^\bbN )$, and let $\bfb\vee\bfc$ be the algebra generated
by $\bfb$ and $\bfc$.

\begin{defi} A {\em cyclic chain} is a sequence $a_1,\dots ,a_n$ of distinct
atoms of $\bfb\vee\bfc$ belonging to both $\bfb$ and $\bfc$ such
that $\psi (a_1)=a_2,\dots , \psi (a_{n-1})=a_n,\psi (a_n)=a_1$:
\begin{equation}\label{cyclic}
\begin{array}{ccccccc}
\bfb & &\;\;\;\; a_1 & \;\;\;\;a_2 && \;\;\;\;a_{n-1} & \;\;\;\;a_n\\
& \psi & \;\;\;\downarrow & \;\;\;\downarrow &\dots &
\;\;\;\downarrow &
\;\;\;\downarrow\\
\bfc & & \;\;\;\;a_2 & \;\;\;\;a_3 && \;\;\;\;a_n & \;\;\;\;a_1
\end{array}
\end{equation}
\end{defi}

\begin{defi} A {\em stable chain} is a sequence $a_0,\dots ,a_n$ of distinct
atoms of $\bfb\vee\bfc$ plus an element $c$, which we call its {\em
end}, such that one of the two following situations occur: 

\noindent{\bf (I)}
\begin{enumerate}
\item $a_0,\dots ,a_n$ belong to $\bfb$,
\item $a_1,\dots ,a_n, c$ belong to $\bfc$,
\item $\psi (a_0)=a_1,\dots ,\psi (a_{n-1})=a_n$,
\item $c=\psi (a_n)=a_0\vee b_1\vee \dots\vee b_k =a_0\vee x$, where $b_1,\dots ,b_k\neq a_0$ are atoms of
$\bfb\vee\bfc$ and $k\geq 1$.
\end{enumerate}
Diagramatically,
\begin{equation}\label{stable1}
\begin{array}{ccccccc}
\bfb &  & \;\;\;\;a_0   & \;\;\;\;a_1   &  & a_{n-1}
&\;\;\;\;a_n\\
& \psi & \;\;\;\downarrow & \;\;\;\downarrow &\dots&
\downarrow & \;\;\;\downarrow\\
\bfc && \;\;\;\;a_1 & \;\;\;\;a_2 &&a_n     &\;\;a_0\vee x
\end{array}
\end{equation}
where $x=b_1\vee\dots\vee b_k\neq 0$.
 
\noindent  {\bf (II)}
\begin{enumerate}
\item $a_0,\dots ,a_n$ belong to $\bfc$,
\item $a_1,\dots ,a_n, c$ belong to $\bfb$,
\item $\psi^{-1}(a_0)=a_1,\dots,\psi^{-1}(a_{n-1})=a_n$,
\item $c=\psi^{-1}(a_n)=a_0\vee
b_1\vee\dots\vee b_k=a_0\vee x$, where $b_1,\dots ,b_k\neq a_0$ are atoms of
$\bfb\vee\bfc$ and $k\geq 1$.
\end{enumerate}
Diagramatically, 
\begin{equation}\label{stable2}
\begin{array}{ccccccc}
\bfb && \;\;\;\;a_1 & \;\;\;\;a_2 &&a_n& \;\;a_0\vee x\\
& \psi & \;\;\;\downarrow &\;\;\;\downarrow &\dots &
\downarrow&\;\;\;\downarrow\\
\bfc && \;\;\;\;a_0 & \;\;\;\;a_1 &&a_{n-1}& \;\;\;\;a_n
\end{array}
\end{equation}
where $x=b_1\vee\dots\vee b_k\neq 0$.

For a stable chain as above we say that an atom $b$ in
$\bfb\vee\bfc$ such that $b\neq a_0$ and $b<\psi (a_n)$
(respectively, $b <\psi^{-1}(a_n)$) is {\em free}. These are the
elements $b_1, \ldots,b_k$ above. Moreover, the atom $a_0$ is said
to be the {\em beginning} of the stable chain.
\end{defi}

\begin{defi}A {\em linking chain} is a sequence $a_1,\dots ,a_n$ of distinct
atoms of $\bfb\vee\bfc$, such that
\begin{enumerate}
\item $a_1,\dots ,a_{n-1}$ belong to $\bfb$,
\item $a_2,\dots ,a_n$ belong to $\bfc$,
\item $\psi (a_1)=a_2,\dots,\psi (a_{n-1})=a_n$,
\item $a_1, a_n$ are free in some stable chains.
\end{enumerate}
\begin{equation}\label{linking}
\begin{array}{cccccccccc}
\bfb &&&& \;\;\;\;a_1 & \;\;\;\;a_2 && \;\;\;\;a_{n-1} && \;\;a_n \vee y\\
& \psi &&& \;\;\;\downarrow & \;\;\;\downarrow &\dots & \;\;\; \downarrow\\
\bfc &&& \;\;\;\;a_1\vee z & \;\;\;\;a_2 & \;\;\;\;a_3 &&
\;\;\;\;a_n
\end{array}
\end{equation}
\end{defi}

\begin{defi} An isomorphism $\psi
:\bfb\rightarrow\bfc$ is said to be {\rm normal} iff
\begin{enumerate}
\item every atom of $\bfb$ or of $\bfc$ that is not an atom in
$\bfb\vee\bfc$ is the end of a stable chain,
\item every atom of $\bfb\vee\bfc$ is a term in either a stable, linking or cyclic chain, or is free in some stable
chain.
\end{enumerate}
\end{defi}

For any finite subalgebras $\bfb$ and $\bfc$ of a common algebra, we
let $n(\bfb ,\bfc )$ be the number of atoms in $\bfb$ that are not
atoms in $\bfb\vee\bfc$ plus the number of atoms in $\bfc$ that are
not atoms in $\bfb\vee\bfc$.

\begin{lemme} For any isomorphism $\psi :\bfb\rightarrow\bfc$
between finite subalgebras $\bfb ,\bfc\subseteq{\rm clop}(2^\bbN )$,
there is a refinement $\psi':\bfb'\rightarrow\bfc'$ satisfying
condition (i) of normality.\end{lemme}

\begin{proof}  The proof is by induction on $n(\bfb,\bfc)$.
For the basis of the induction, if $n(\bfb,\bfc)=0$, then every atom
of $\bfb$ and every atom of $\bfc$ is an atom of $\bfb\vee\bfc$ and
there is noting to prove. In this case, we see that the structure of
$\psi$ is particularly simple, since it is then just an automorphism
of $\bfb=\bfc=\bfb\vee\bfc$ and thus one can easily split
$\bfb\vee\bfc$ into cyclic chains, namely the $\psi$-orbits of
atoms.

Now for the induction step: Suppose $x\in\bfb$ is some atom of
$\bfb$ that is not an atom of $\bfb\vee \bfc$ (the case when $x$ is
an atom of $\bfc$ that is not an atom of $\bfb\vee\bfc$ is of course
symmetric), and trace the longest chain of atoms of $\bfb\vee\bfc ,\
a_i\in\bfb ,b_i\in\bfc$, such that
\begin{equation}
x\stackrel{\psi}{\mapsto}b_0=a_0\stackrel{\psi}{\mapsto}b_1=a_1
\stackrel{\psi}{\mapsto}b_2=a_2\stackrel{\psi}{\mapsto}\dots
\end{equation}
If for some $i<j,\ a_i=a_j$, then also $b_i=b_j$, and by the
injectivity of $\psi ,\ a_{i-1}=a_{j-1}$, etc.  So in that case, by
induction, we have $b_0=b_{j-i}$, whence $a_{j-i-1}=x$,
contradicting that $x$ was not an atom in $\bfb\vee\bfc$.

So as $\bfb\vee\bfc$ is finite, the chain has to stop either (i) on
some $a_n$ or (ii) on some $b_n$ (if it stops with $x$ we let
$x=a_{-1}$ and treat it as in case (i)).

Case (i):  Let $y=\psi (a_n)$, which is an atom in $\bfc$ but not in
$\bfb\vee\bfc$.  Now, as $\bigvee^n_{i=1}a_i=\bigvee^n_{i=1} b_i$,
both $x$ and $y$ are disjoint from $\bigvee^n_{i=0}a_i$ and let
$x=c_1\vee\dots\vee c_p, y=d_1\vee\dots\vee d_q$ be the
decompositions into atoms of $\bfb\vee\bfc$. Suppose $p\leq q$ (the
case when $q\leq p$ is similar).  Split $a_i=b_i$ into non zero
elements $a^1_i,\dots ,a^p_i$ of clop$(2^\bbN )$ for each $i\leq n$.
Then we let $\bfb'$ be the smallest algebra containing $\bfb$ and
new atoms $a^l_i$, for $l\leq p, i\leq n$, and $c_1,\dots ,c_p$.
Also let $\bfc'$ be the smallest algebra containing $\bfc$ and new
elements $a^l_i$ for $l\leq p, i\leq n$, and $d_1,d_2,\dots
,d_{p-1},d_p\vee \dots\vee d_q$. Finally, let $\psi'$ be the unique
extension of $\psi$ satisfying $\psi'(c_l
)=a^l_0,\;\psi'(a^l_i)=a^l_{i+1}$, for $i<n,\;l\leq p$, and
$\psi'(a^l_n)=d_l$ for $l <p,\;\psi'(a^p_n)= d_p\vee\dots\vee d_q$.
We remark that $n(\bfb',\bfc')<n(\bfb,\bfc )$.

Case (ii):  Notice that
$x\wedge\bigvee^{n-1}_{i=1}b_i=x\wedge\bigvee^{n-1}_{i=1} a_i=0$,
and, as $b_n$ is an atom of $\bfb\vee\bfc$, either $b_n<x$ or $b_n
\wedge x=0$.  If $b_n<x$, then $x$ is the end of a stable chain. If
$b_n \wedge x=0$, then we proceed as follows:

Suppose $x=c_1\vee\dots\vee c_p$ is the decomposition into atoms of
$\bfb\vee\bfc$ and split $b_i\ (i\leq n)$ into $b^1_i,\dots ,b^p_i$
in clop$(2^\bbN )$.  Now let $\bfb'$ be the algebra generated by
$\bfb$ and $b^1_i,\dots ,b^p_i$, for $i<n$ and $c_1,\dots ,c_p$.
Let $\bfc'$ be the algebra generated by $\bfc$ and $b^ 1_i,\dots
,b^p_i$ for $i\leq n$, and let $\psi'$ be the unique extension of
$\psi$ satisfying $\psi'(c_l ) =b^l_0\ (l\leq p),\;\psi
(b^l_i)=b^l_{i+1}\ (l\leq p)$.  Again we remark that
$n(\bfb',\bfc')<n(\bfb ,\bfc )$.  So this concludes the induction
step.\end{proof}

\begin{prop}\label{normal} Any isomorphism $\psi
:\bfb\rightarrow\bfc$ between finite subalgebras $\bfb
,\bfc\subseteq{\rm clop}(2^\bbN )$ has a finite normal
refinement.\end{prop}

\begin{proof} By the lemma we can suppose that $\psi :\bfb\rightarrow\bfc$
satisfies condition (i) of normality.  We will then see that it
actually satisfies condition (ii) too.

So consider an atom $a_0$ of $\bfb\vee\bfc$.  Find an atom $c$ of
$\bfb$ such that $a_0\leq c$.  Now if $a_0<c$, then we know from
condition (i) that $c$ is the end of a stable chain and $a_0$ is
therefore either the beginning of a stable chain or is free.
Otherwise if $a_0=c\in\bfb$, then we find the maximal chain of atoms
of $\bfb\vee\bfc ,a_i\in\bfb , b_i\in\bfc$ such that
\begin{equation}
\dots\stackrel{\psi}{\mapsto}b_i=a_i\stackrel{\psi}{\mapsto}b_{i+1}=
a_{i+1}\stackrel{\psi}{\mapsto}b_{i+2}=a_{i+2}\stackrel{\psi}
{\mapsto}\dots
\end{equation}
where the indices run over an interval of $\bbZ$ containing 0.  We
have now various cases:

Case (i):  $a_i=a_j$ for some $i<j$.  Then obviously
$b_{i+1}=b_{j+1}$ and $b_{i-1}=b_{j-1}$ by the injectivity of $\psi$
(notice also that these elements are indeed defined, i.e., are atoms
in $\bfb\vee\bfc$ and not only in $\bfc$).  But then also
$a_{i+1}=a_{j+1}$ and $a_{i-1}= a_{j-1}$, etc. So the chain is
bi-infinite and periodic and we can write it as
\begin{equation}
\begin{array}{cccccc}
\bfb && \;\;\;\;a_0 & \;\;\;\;a_1 && \;\;\;\;a_n\\
& \psi & \;\;\;\downarrow & \;\;\;\downarrow &\dots & \;\;\;
\downarrow\\
\bfc && \;\;\;\;a_1 &\;\;\;\; a_2 && \;\;\;\;a_0
\end{array}
\end{equation}
So $a_0$ is term in a cyclic chain.

Case (ii):  The chain ends with some $a_k$.  Then $c=\psi (a_k)$ is
an atom of $\bfc$ that is not an atom of $\bfb\vee\bfc$ and must
therefore be the end of some stable chain.  So using the bijectivity
of $\psi$ one sees that all of the above chain are terms in a stable
chain and in particular $a_0$ is a term in a stable chain.

Case (iii):  The chain begins with some $b_k$.  Let
$c=\psi^{-1}(b_k)$. Then $c$ must be some atom of $\bfb$ that is not
an atom of $\bfb\vee\bfc$, so it must be the end of some stable
chain.  Now looking at $\psi^{-1}$ instead of $\psi$ and switching
the roles of $\bfb$ and $\bfc$, we see as before that $a_0$ is a
term in  some stable chain.

Case (iv):  The chain begins with some $a_n$ and ends with some
$b_k$. Then $a_n\neq b_k$ and there must be atoms $c$ of $\bfc$ and
$d$ of $\bfb$ that are ends of stable chains such that $a_n\wedge
c\neq 0$ and $b_k \wedge d\neq 0$.  But, as $a_n$ and $b_k$ are
atoms, this means that $a_n<c$ and $b_k<d$, so the chain is
linking.\end{proof}

\section{Generic homeomorphisms of Cantor space}

A {\it generic homeomorphism} of the Cantor space $2^\bbN$ is a
homeomorphism whose conjugacy class in $H(2^\bbN )$ is dense
$G_\delta$.  As $H(2^\bbN )$ (with the uniform topology) is
isomorphic as a topological group to Aut$(\bfb_\infty )$, where
$\bfb_\infty$ is the countable atomless boolean algebra, the
existence of a generic homeomorphism of $2^\bbN$ is equivalent to
the existence of a generic automorphism of $\bfb_\infty$.

\begin{thm} Let $\kkk$ be the class of finite boolean
algebras.  Then $\kkk_p$ has the {\rm CAP}, and therefore there is a
generic automorphism of the countable atomless boolean algebra and a
generic homeomorphism of the Cantor space.\end{thm}

Before we embark on the proof, let us first mention that one can
construct counter-examples to $\kkk_p$ having the AP, and that
therefore the added complications are necessary.

\begin{proof} We will show that the class $\llll$ of $\s =\langle\bfb\vee\bfc
, \psi :\bfb\rightarrow\bfc\rangle$, where $\psi$ is normal, has the
AP. Notice first that $\llll$ is cofinal in $\kkk_p$ by Proposition
\ref{normal}.

Suppose now $\s =\langle\bfa ,\psi
\colon\bfb\rightarrow\bfc\rangle\in\llll$.  Let $\s^l =\langle\bfa^l
,\psi^l \colon\bfb^l\rightarrow \bfc^l\rangle$ and
$\s^r=\langle\bfa^r,\psi^r:\bfb^r\rightarrow\bfc^r\rangle$ be two
refinements of $\s$, which we do not necessarily demand belong to
$\llll$.

We claim that they can be amalgamated over $\s$. We remark first
that it is trivially enough to amalgamate any two refinements of
$\s^l$ and $\s^r$ over $\s$. List the atoms of $\bfa$ as $a_1,\dots
,a_n$.  Then by refining $\s^l , \s^r$ we can suppose they have
atoms
\begin{equation}
a^l_1(1),\dots ,a^l_1(k),\dots ,a^l_n(1), \dots ,a^l_n(k)
\end{equation}
and
\begin{equation}
a^r_1(1),\dots ,a^r_1(k),\dots ,a^r_n(1),\dots , a^r_n(k)
\end{equation}
respectively, where $a^v_t(1),\ldots ,a^v_t(k)$ ($v=l,r$ and
$t=1,\ldots, n$) is a splitting of $a_t$.  So $a_t\mapsto
\bigvee_{j=1}^{k} a^l_t (j)$ and $a_t\mapsto\bigvee_{j=1}^{k}
a^r_t(j)$ are canonical embeddings of $\bfa$ into $\bfa^l$ and
$\bfa^r$ and we can furthermore suppose that both $\bfb^l ,\bfc^l$
and $\bfb^r,\bfc^r$ contain the image of $\bfa$ by these embeddings.
This can be done by further refining $\s^l$ and $\s^r$. This means
that any atom of $\bfb^v$ and $\bfc^v$ ($v=l,r$) will be on the form
$\bigvee_{i\in\Gamma}a^v_t(i)$, where $\Gamma\subseteq
\{1,\ldots,k\}$ and  $1\leq t\leq n$.

Take new formal atoms $a^l_m (i)\otimes a^r_m (j)$ for $m\leq n$ and
$i,j\leq k$.  Our amalgam
$\s^a=\langle\bfa^a,\psi^a:\bfb^a\rightarrow \bfc^a\rangle$ will be
such that the atoms in $\bfa^a$ will be a subset $E$ of these new
formal atoms and the embeddings $l :\bfa^l\rightarrow \bfa^a$ and
$r:\bfa^r\rightarrow\bfa^a$ will be defined by
\begin{equation}
l (a^l_t(i))=\bigvee\{a^l_t(i)\otimes a^r_t(j)\in E\},
\end{equation}
\begin{equation}
r(a^r_t(i))=\bigvee\{a^l_t(j)\otimes a^r_t(i)\in E\}
\end{equation}
The atoms of $\bfb^a$ will be the following
\begin{equation}
E^{\Gamma ,\Delta}_t=\bigvee\{a^l_t(i)\otimes a^r_t(j)\in E:
(i,j)\in\Gamma \times\Delta\},
\end{equation}
where $\bigvee_{i\in\Gamma}a^l_t(i)$ is an atom in $\bfb^l$ and
$\bigvee_{j\in\Delta}a^r_t(j)$ is an atom in $\bfb^r$. Similarly for
$\bfc^a$.

Now obviously it is enough to define $\psi^a$ between the atoms of
$\bfb^a$ and $\bfc^a$.

(1) Suppose $\Gamma ,\Delta ,\Theta ,\Lambda\subseteq\{1,\dots ,k\}$
are such that
\begin{equation}
\bigvee_{i\in\Gamma}a^l_t(i), \;\bigvee_{i\in\Delta}a^r_t(i),
\;\bigvee_{i\in\Theta}a^l_m(i),\;\bigvee_{i\in\Lambda}a^r_m(i)
\end{equation}
are atoms in $\bfb^l ,\bfb^r,\bfc^l ,\bfc^r$, respectively, such
that
\begin{equation}
\psi^l (\bigvee_\Gamma a^l_t(i))=\bigvee_\Theta a^l_m(i)
\;\;\textrm{ and }\;\; \psi^r(\bigvee_\Delta
a^r_t(i))=\bigvee_\Lambda a^r_m(i).
\end{equation}
Then we let $\psi^a (E^{\Gamma ,\Delta}_t )=E^{\Theta ,\Lambda}_m$.

(2) Now suppose that $\bigvee_\Gamma a^l_t(i)$, $\bigvee_\Delta
a^r_t(i),\bigvee_\Theta a^l_m(i)$, $\bigvee_\Lambda a^r_s(i) ,m\neq
s$, are atoms in $\bfb^l$, $\bfb^r$, $\bfc^l$, $\bfc^r$,
respectively, such that
\begin{equation}\label{problem}
\psi^l (\bigvee_\Gamma a^l_t(i))=\bigvee_\Theta a^l_m(i)
\;\;\textrm{ and }\;\; \psi^r(\bigvee_\Delta
a^r_t(i))=\bigvee_\Lambda a^r_s(i).
\end{equation}
We notice first that, since
$\s=\langle\bfa,\psi\colon\bfb\rightarrow\bfc\rangle\in\llll$ is
normal with $\bfa=\bfb\vee\bfc$, the only time $\psi$ sends an atom
of $\bfa $ to a non-atom of $\bfa$, or inversely sends a non-atom of
$\bfa$ to an atom of $\bfa$, is in the last steps of a stable chain
corresponding to $a_n\mapsto a_0\vee x$ in diagram (\ref{stable1}) or
$a_0\vee x\mapsto a_0$ in diagram (\ref{stable2}). Moreover, whenever
$x$ is an atom of $\bfb$, then either $x$ or $\psi(x)$ is an atom of
$\bfa$.

Using this, it follows from equations (\ref{problem})  and $m\neq s$
that, since $\psi^l$ and $\psi^r$ are refinements of $\psi$, $a_t$
must be an atom of $\bfb$, while $a_m\vee a_s$ must be below the end
of a stable chain in $\s$.

Now, in this  situation we cannot have $a^l_t(i)\otimes a^r_t(j)\in
E$, for any $(i,j)\in \Gamma\times\Delta$, or, in other words, we
must have $E^{\Gamma ,\Delta}_t=0$. This is because we have $\psi^l
(\bigvee_\Gamma a^l_t(i))\leq a_m$, while $\psi^r(\bigvee_\Delta
a^r_t(i))\leq a_s$, which would force
\begin{equation}
\psi^a(E^{\Gamma ,\Delta}_t)\leq\bigvee\{a^l_m(i)\otimes a^r_m(j)\in
E\}
\end{equation}
and similarly
\begin{equation}
\psi^a(E^{\Gamma ,\Delta}_t)\leq\bigvee\{a^l_s(i)\otimes a^r_s(j)\in
E\},
\end{equation}
which leaves only the possibility $\psi^a(E^{\Gamma ,\Delta}_t)=0$.
Or said in another way, it would force
\begin{equation}
\psi^a(E^{\Gamma ,\Delta}_t)=\psi^a(\bigvee\{a^l_t(i)\otimes
a^r_t(j)\in E: (i,j)\in\Gamma \times\Delta\})
\end{equation}
to be the join of elements on the form $a^l_m(p)\otimes a^r_s(q)$.
But we do not include any such elements in our amalgam.

There is also a dual version of this problem, namely when
\begin{equation}
\psi^l (\bigvee_\Gamma a^l_m(i))=\bigvee_\Theta a^l_t(i)
\;\;\textrm{ and }\;\; \psi^r(\bigvee_\Delta
a^r_s(i))=\bigvee_\Lambda a^r_t(i)
\end{equation}
for distinct $s,m$.

With (1) and (2) in mind, we can now formulate the necessary and
sufficient conditions on $E$ for this procedure to work out:

(a) For $l$ to be well-defined as an embedding from $\bfa^l$ to
$\bfa^a$:
$$
\forall t\leq n\;\forall i\leq k\;\exists j\leq k\ \;a^l_t(i)\otimes
a^r_t (j)\in E
$$

(b) For $r$ to be well-defined as an embedding from $\bfa^r$ to
$\bfa^a$:
\[
\forall t\leq n\;\forall j\leq k\;\exists i\leq k\ \;a^l_t(i)\otimes
a^r_t (j)\in E
\]

(c) If $\Gamma ,\Delta ,\Theta ,\Lambda$ and $t,m$ are as in (1),
then $E^{\Gamma ,\Delta}_t\neq 0\Leftrightarrow E^{\Theta
,\Lambda}_m\neq 0$.

(d) If $\Gamma ,\Delta ,\Theta ,\Lambda$ and $t,m,s$ are as in (2),
then $E^{\Gamma ,\Delta}_t=0$.

We will define $E$ separately for terms of stable, linking and
cyclic chains.

Suppose we are given a cyclic chain
\begin{equation}
\begin{array}{cccccc}
\bfb & & \;\;\;\;a_{t_1} & \;\;\;\;a_{t_2} & & \;\;\;\;a_{t_w}\\
& \psi & \;\;\;\downarrow & \;\;\;\downarrow & \dots & \;\;\;
\downarrow\\
\bfc & & \;\;\;\;a_{t_2} & \;\;\;\;a_{t_3} & & \;\;\;\;a_{t_1}
\end{array}
\end{equation}
Then we let $a^l_{t_v}(i)\otimes a^r_{t_v}(j)\in E$ for all $i,j\leq
k$ and $1\leq v\leq w$.

Given a linking chain
\begin{equation}
\begin{array}{cccccccc}
\bfb & & & \;\;\;\;a_{t_1} &\;\;\;\; a_{t_2} & & \;\;\;\;a_{t_{w-1}}
& \;\;a_{t_w}\vee y\\
& \psi && \;\;\;\downarrow & \;\;\;\downarrow & \dots & \;\;\;
\downarrow\\
\bfc && \;\;x\vee a_{t_1} & \;\;\;\;a_{t_2} & \;\;\;\;a_{t_3} & &
\;\;\;\;a_{t_w}
\end{array}
\end{equation}
we let $a^l_{t_v}(i)\otimes a^r_{t_v}(j)\in E$ for all $i,j\leq k$
and $1\leq v\leq w$.

So, since we have included all the new formal atoms in these two
cases, we easily see that (a), (b) and (c) are verified. Moreover,
condition (d) is only relevant for the stable chains as it only
pertains to (2).

Finally, suppose we have a stable chain
\begin{equation}
\begin{array}{ccccccl}
\bfb & & \;\;\;\;a_{t_1} & \;\;\;\;a_{t_2} & &a_{t_{w-1}}&
\;\;\;\;\;\;\;\;a_{t_w}\\

& \psi & \;\;\;\downarrow & \;\;\;\downarrow & \dots &\downarrow&
\;\;\;\;\;\;\;\;\;\downarrow\\

\bfc & &\;\;\;\; a_{t_2} & \;\;\;\;a_{t_3} & &a_{t_w}&
\;\;\;a_{t_1}\vee a_{t_{w+1}}\vee\dots\vee a_{t_{w+q}}
\end{array}
\end{equation}
The case when $a_{t_1}\in\bfc$ and $a_{t_1}\vee
a_{t_{w+1}}\vee\dots\vee a_{t_{w+q}}\in\bfb$ is symmetric to this
one and is taken care of in the same manner.

By choosing $k$ big enough and refining the partitions of
$a_{t_1},\dots , a_{t_w}$ in $\bfb^l$ and $\bfb^r$, the partitions
of $a_{t_2},\dots , a_{t_w},a_{t_1}\vee a_{t_{w+1}}\vee\dots\vee
a_{t_{w+q}}$ in $\bfc^l , \bfc^r$ and subsequently extending $\psi^l
,\psi^r$ to these refined partitions, we can suppose that there are
partitions of $\{1,\dots ,k\}$ as follows:
\begin{equation}\begin{split}
\{1,\dots,k\}=&\Gamma^l (\gamma ,1)\;\sqcup\;\dots\;\sqcup\;\Gamma^l(\gamma ,p)\; \sqcup\;\Delta^l(\gamma,1,1)\;\sqcup\;\dots\;\\
&\sqcup\;\Delta^l (\gamma ,1,p)\;
\sqcup\;\dots\;\sqcup\;\Delta^l (\gamma,q,1)\;\sqcup\;\dots\;\sqcup\;\Delta^l (\gamma ,q,p)
\end{split}\end{equation}
for $\gamma=1,\dots ,w$, and
\begin{equation}
\Lambda^l(1)\sqcup\dots\sqcup\Lambda^l (p)=\{1,\dots ,k\}
\end{equation}
\begin{equation}
\Theta^l (\beta ,1) \sqcup\dots\sqcup\Theta^l (\beta,p)=\{1,\dots k\}
\end{equation}
for $\beta =1,\dots ,q$, such that
\begin{itemize}

\item for $e=1,\dots ,p;\gamma =1,\dots ,w;\beta =1,\dots q$\;:

$\bigvee_{\Gamma^l (\gamma ,e)}a^l_{t_\gamma}(i),\bigvee_{\Delta^l
(\gamma ,\beta ,e)}a^l_{t_\gamma}(i)$ are atoms in $\bfb^l$

\item for $e=1,\dots p;\gamma =2,\dots ,w;\beta =1,\dots ,q$\;:

$\bigvee_{\Gamma^l (\gamma ,e)}a^l_{t_\gamma}(i),\bigvee_{\Delta^l
(\gamma ,\beta ,e)}a^l_{t_\gamma}(i)$ are atoms in $\bfc^l$,

\item for $e=1,\dots ,p\;:\ \bigvee_{\Lambda^l (e)}a^l_{t_1}(i)$ is
an atom in $\bfc^l$,

\item for $e=1,\dots ,p;\beta=1,\dots ,q$\;: $\bigvee_{\Theta^l
(\beta ,e)}a^l_{t_{w+\beta}}(i)$ is an atom in $\bfc^l$,

\item for $e=1,\dots ,p;\gamma =1,\dots ,w-1$\;:
\begin{equation}
\psi^l\left (\bigvee_{\Gamma^l (\gamma ,e)}a^l_{t_\gamma}(i) \right
)= \bigvee_{\Gamma^l (\gamma +1,e)}a^l_{t_{\gamma +1}}(i),
\end{equation}

\item for $e=1,\dots ,p;\gamma =1,\dots ,w-1;\beta =1,\dots ,q$\;:
\begin{equation}
\psi^l\left (\bigvee_{\Delta^l (\gamma ,\beta ,e)}a^l_{t_\gamma}(i)
\right ) =\bigvee_{\Delta^l (\gamma +1,\beta ,e)}a^l_{t_{\gamma
+1}}(i),
\end{equation}

\item for $e=1,\dots ,p$\;:
\begin{equation}
\psi^l\left (\bigvee_{\Gamma^l (w,e)}a^l_{t_w}(i)\right )=
\bigvee_{\Lambda^l (e)}a^l_{t_w}(i),
\end{equation}

\item for $e=1,\dots ,p;\beta =1,\dots q$\;:
\begin{equation}
\psi^l\left (\bigvee_{\Delta^l (w,\beta ,e)}a^l_{t_w}(i) \right )=
\bigvee_{\Theta^l (\beta ,e)}a^l_{t_{w+\beta}}(i).
\end{equation}
\end{itemize}
And the same for $r$, i.e., with the same constants $k,p$, though
different partitions.

Let us see how this can be obtained.  Suppose $\psi^l :\bfb^l
\rightarrow\bfc^l$ is given.  First refine the partition of
$a_{t_2}$ in $\bfb^l$ to be finer than the partition of $a_{t_2}$ in
$\bfc^l$. Then refine the partition of $a_{t_3}$ in $\bfc^l$ so that
we can extend $\psi^l$ to be an isomorphism of the refined
partitions of $a_{t_2}$ in $\bfb^l$ and $a_{t_3}$ in $\bfc^l$,
respectively.  Again we refine the partition of $a_{t_3}$ in
$\bfb^l$ to be finer than the partition of $a_{t_3}$ in $\bfc^l$.
And refine the partition of $a_{t_4}$ in $\bfc^l$ so that we can
extend $\psi^l$ to be an isomorphism of the refined partitions of
$a_{t_3}$ in $\bfb^l$ and $a_{t_4}$ in $\bfc^l$ respectively. We
continue like this until we have refined the partition of $a_{t_w}$
in $\bfb^l$. Extend now the partition of $a_{t_1}\vee
a_{t_{w+1}}\vee\dots\vee a_{t_{w+q}}$ in $\bfc^l$, so we can let
$\psi^l$ be an isomorphism of the refined partition of $a_{t_w}$ in
$\bfb^l$ with the refined partition of $a_{t_1}\vee
a_{t_{w+1}}\vee\dots\vee a_{t_{w+q}}$ in $\bfc^l$. We will now
refine the elements of the chain in the other direction.  We recall
first that the partition of $a_{t_1}\vee a_{t_{w+1}} \vee\dots\vee
a_{t_{w+q}}$ in $\bfc^l$ was supposed to be fine enough to contain
the elements $a_{t_1}, a_{t_{w+1}},\dots ,a_{t_{w+q}}$.

To begin, refine the new partition of $a_{t_w}$ in $\bfc^l$ to be
the same as the new partition of $a_{t_w}$ in $\bfb^l$. Refine now
the partition of $a_{t_{w-1}}$ in $\bfb^l$ so that $\psi^l$ extends
to an isomorphism of the refined partitions of $a_{t_{w-1}}$ in
$\bfb^l$ and $a_{t_w}$ in $\bfc^l$.  Again refine the partition of
$a_{t_{w-1}}$ in $\bfc^l$ to be the same as the partition of
$a_{t_{w-1}}$ in $\bfb^l$.  And refine the partition of
$a_{t_{w-2}}$ in $\bfb^l$ so that $\psi^l$ extends to an isomorphism
of the refined partitions of $a_{t_{w-2}}$ in $\bfb^l$ and
$a_{t_{w-1}}$ in $\bfc^l$. Continue in this fashion until $a_{t_1}$
has been refined in $\bfb^l$.

This shows that the above form can be obtained for $\psi^l :\bfb^l
\rightarrow\bfc^l$ individually.  Now we can of course also obtain
it for $\psi^r:\bfb^r\rightarrow\bfc^r$.  So the only thing lacking
is to obtain the same constants $k,p$ for both $\psi^l
:\bfb^l\rightarrow \bfc^l$ and $\psi^r:\bfb^r\rightarrow\bfc^r$. But
this is easily done by some further splitting (this last bit is only
to reduce complexity of notation).

We now define a derivation $D$ on $\p (\{1,\dots ,k\}^2)$ ( i.e., a
mapping satisfying $D(X)\subseteq X$ for all $X\subseteq \{1,\dots
,k\}^2$) as follows:
\begin{equation}
D(Y)=\{(i,j)\in Y: (i,j)\in\Gamma^l
(1,e)\times\Gamma^r(1,d)\Rightarrow (\Lambda^l
(e)\times\Lambda^r(d))\cap Y\neq\emptyset\}.
\end{equation}
Let
\begin{equation}
X_0=\bigcup_{e,d}\Gamma^l
(1,e)\times\Gamma^r(1,d)\cup\bigcup^q_{\beta =1}
\bigcup_{e,d}\Delta^l (1,\beta ,e)\times\Delta^r (1,\beta ,d).
\end{equation}
So $D^{k_\infty+1} (X_0)=D^{k_\infty}(X_0)$ for some minimal
$k_\infty\in\bbN$, and we notice first that
\begin{equation}
\bigcup^q_{\beta =1}\bigcup_{e,d}\Delta^l (1,\beta ,e)\times\Delta^r
(1,\beta ,d)\subseteq D^{k_\infty} (X_0)
\end{equation}

{\bf Claim}.  {\it Suppose that for some $e,d,\sigma$ and
$(i,j)\in\Gamma^l (1,e) \times\Gamma^r(1,d)$, we have $(i,j)\not\in
D^{\sigma +1}(X_0)$. Then $(\Lambda^l (e)\times\Lambda^r(d))\cap
D^\sigma (X_0)=\emptyset$.}

\begin{proof}  By definition of $X_0$, there is a $\tau\geq 0$ such that
$(i,j)\in D^\tau (X_0)\setminus D^{\tau +1}(X_0)$ (so $\tau\leq
\sigma$).  And by definition of $D$, we must have $(\Lambda^l
(e)\times\Lambda^r (d))\cap D^\tau (X_0)=\emptyset$, whence also
$(\Lambda^l (e)\times\Lambda^r (d))\cap D^\sigma
(X_0)=\emptyset$.\end{proof}

\begin{lemme}\label{derivation} For all $i\in\{1,\dots ,k\}$ we have
$\{i\}\times \{1,\dots ,k\}\cap D^{k_\infty}(X_0)\neq\emptyset$ and
$\{1,\dots ,k\}\times \{i\}\cap
D^{k_\infty}(X_0)\neq\emptyset$.\end{lemme}

\begin{proof}  Otherwise we can take $\tau$ minimal such that for some $i$ we
have, e.g., $\{i\}\times\{1,\dots ,k\}\cap D^\tau (X_0)=\emptyset$.
Clearly $\tau >0$, so, as $\bigcup^r_{\beta =1}\bigcup_{e,d}\Delta^l
(1,\beta ,e)\times\Delta^r (1,\beta ,d)\subseteq D^{k_\infty}(X_0)$,
we must have $i\in\Gamma^l (1,e)$ for some $e$.  Therefore, by the
claim, letting $\sigma +1=\tau$, we have $\Lambda^l
(e)\times\{1,\dots ,k\}\cap D^{k_\infty}(X_0)=\emptyset$. So, in
particular, for any $j\in\Lambda^l (e)$ we have
$\{j\}\times\{1,\dots ,k\}\cap D^\sigma (X_0)=\emptyset$,
contradicting the minimality of $\tau$.\end{proof}
\medskip

Now put $a^l_{t_1}(i)\otimes a^r_{t_1}(j)\in E$ for all $(i,j)\in
D^{k_\infty} (X_0)$.  And put $a^l_{t_\gamma}(i)\otimes
a^r_{t_\gamma}(j)\in E$ for all $(i,j)\in\Gamma^l (\gamma
,e)\times\Gamma^r(\gamma ,d)$ such that $(\Gamma^l
(1,e)\times\Gamma^r(1,d))\cap D^{k_\infty}(X_0)\neq\emptyset$ and
for all $(i,j)\in\Delta^l (\gamma ,\beta ,e)\times\Delta^r(\gamma ,
\beta ,d)\ (\gamma =2,\dots w; \beta =1,\dots ,q)$.

By Lemma \ref{derivation}, (a) and (b) are satisfied for $a_{t_1}$
and therefore also by construction for $a_{t_2},\dots ,a_{t_w}$.
Moreover, (c) is satisfied between $a_{t_1},a_{t_2}$, between
$a_{t_2},a_{t_3},\dots ,$ between $a_{t_{w-1}}, a_{t_w}$.
Furthermore, by definition of $X_0$, condition (d) is satisfied in
the only relevant place, namely between $a_{t_w}$ and $a_{t_1}\vee
a_{t_{w+1}}\vee\dots\vee a_{t_{w+q}}$.  So we only have to check
condition (c) between $a_{t_w}$ and $a_{t_1}\vee
a_{t_{w+1}}\vee\dots\vee a_{t_{w+q}}$.

Now obviously, since $a^l_{t_w}(i)\otimes a^r_{t_w}(j)\in E$, for
$(i,j) \in\Delta^l (w,\beta ,e)\times\Delta^r (w,\beta ,d)\ (\beta
=1,\dots ,q)$, we only have to check the condition for the products
$\Gamma^l (w,e) \times\Gamma^r (w,d)$.  So suppose
$a^l_{t_w}(i)\otimes a^r_{t_w}(j)\in E$ for some $(i,j)\in\Gamma^l
(w,e)\times\Gamma^r (w,d)$.  Then we know that $(\Gamma^l
(1,e)\times\Gamma^r (1,d))\cap D^{k_\infty} (X_0)\neq \emptyset$,
therefore, as $D^{k_\infty} (X_0)$ is $D$-stable, we must have
$(\Lambda^l (e)\times\Lambda^r(d))\cap D^{k_\infty}
(X_0)\neq\emptyset$. This means that, if
\begin{equation}
\bigvee\{a^l_{t_w}(i)\otimes a^r_{t_w}(j)\in E: (i,j)\in\Gamma^l
(w,e) \times\Gamma^r (w,d)\}\neq 0,
\end{equation}
then also
\begin{equation}
\bigvee\{a^l_{t_1}(i)\otimes a^r_{t_1}(j)\in E: (i,j)\in\Lambda^l
(e) \times\Lambda^r(d)\}\neq 0,
\end{equation}
confirming (c) in one direction.  And, conversely, if
\begin{equation}
\bigvee\{a^l_{t_1}\otimes a^l_{t_1}(j)\in E: (i,j)\in\Lambda^l (e)
\times\Lambda^r(d)\}\neq 0,
\end{equation}
then, in particular, $(\Lambda^l (e)\times\Lambda^r(d))\cap
D^{k_\infty} (X_0) \neq\emptyset$.  So, by the claim, also
$(\Gamma^l (1,e)\times\Gamma^r (1,d))\cap D^{k_\infty}
(X_0)\neq\emptyset$, and therefore
\begin{equation}
\bigvee\{a^l_{t_w}(i)\otimes a^r_{t_w}(j)\in E: (i,j)\in\Gamma^l
(w,e) \times\Gamma^r(w,d)\}\neq 0.
\end{equation}
\end{proof}
\medskip

One can also use well-known results (see e.g. Hjorth \cite{hjorth})
to see that there is a generic increasing homeomorphism of $[0,1]$.
\begin{thm} The group $H_+([0,1])$ of increasing
homeomorphisms of of the unit interval has a comeager conjugacy
class.\end{thm}

\begin{proof} Using the notation of Hjorth it is easy to verify that the set
of $\pi\in H_+([0,1])$, for which $<^{\mathcal M(\pi)}$ is
isomorphic to $\bbQ$, $P_=^{\mathcal M(\pi)}=\emptyset$ and
$P_+^{\mathcal M(\pi)}$, $P_-^{\mathcal M(\pi)}$ are dense and
unbounded in both directions in $<^{\mathcal M(\pi)}$, is comeager
and forms a single conjugacy class.\end{proof}
\medskip

Notice that this group is connected, so is not a topological
subgroup of $S_\infty$ (In fact, it has been recently proved in
Rosendal--Solecki \cite{rosendal} that it is not even an abstract
subgroup of $S_\infty$ either.)

\section{Ample generic automorphisms}

\subsection{The general concept}
Suppose a Polish group $G$ acts continuously on a Polish
space $X$.  Then there is a natural induced (diagonal) action of
$G$ on $X^n$, for any $n=1,2,\dots$, defined by $g\cdot
(x_1,x_2,\dots, x_n )=(g\cdot x_1,g\cdot x_2, \dots , g\cdot
x_n)$.  Now, by Kuratowski--Ulam, if there is a comeager $G$-orbit
in $X^n$, and $k\leq n$, then there is also a comeager orbit in
$X^k$.  However this is far from being true in the other
direction.  Let us reformulate the property for $n=2$.

Recall that $\forall^*xR(x)$ means that $\{ x:R(x)\}$ is comeager,
where $x$ varies over elements of a topological space $X$.

\begin{prop} Let a Polish group $G$ act
continuously on a Polish space $X$ and suppose $X$ has a comeager
orbit $\ooo$.  Then the following are equivalent:

(i) There is a comeager orbit in $X^2$,

(ii) $\forall x\in\ooo\;\forall^*y\in\ooo (G_x\cdot y$ is comeager in
$X$), where $G_x=\{g\in G:g\cdot x=x\}$ is the {\rm stabilizer} of
$x$,

(iii) $\forall x\in\ooo\;\forall^*y\in\ooo (G_xG_y$ is comeager in
$G$)

(iv) $\forall x\in\ooo\;\forall^*h\in G(G_xhG_x$ is comeager in
$G$),

(v) $\exists x,y\in\ooo (G_x\cdot y$ is comeager).\end{prop}

\begin{proof}  (i) $\Rightarrow$ (ii) Let $C$ be a comeager orbit in $X^2$.
Then by the Kuratowski--Ulam Theorem
$\forall^*x\;\forall^*y(x,y)\in C$.  Thus $\{x:\forall^*y(x,y)\in
C\}$ is comeager and clearly $G$-invariant, so for every $x\in\ooo
,\forall^*y(x,y)\in C$.  It follows that $\forall^*y\in\ooo
(G_x\cdot y$ is comeager).

(ii) $\Rightarrow$ (iii)  Fix $x\in\ooo$ and $y\in\ooo$ such that
$G_x\cdot y$ is comeager.  Since $\ooo$ is $G_\delta$ in $X$, by
Effros' Theorem (see, e.g., Becker-Kechris \cite{becker}) the map $\pi
:G\rightarrow\ooo$ given by $\pi (g)=g\cdot y$ is continuous and
open and therefore $\pi^{-1}(G_x\cdot y)=\{g: g\cdot y \in
G_x\cdot y\}=\{g:\exists h\in G_x(h^{-1}g\in G_y)\}=G_xG_y$ is
comeager in $G$.

(iii) $\Rightarrow$ (iv) We have for any $x\in \ooo,\forall^*y\in\ooo
(G_xG_y$ is comeager in $G$), so, by Effros' Theorem again,
applied this time to $\sigma (h)=h\cdot x$, we have
$\forall^*h(G_xG_{h\cdot x}$ is comeager in $G$).  But $G_{h\cdot
x}=hG_xh^{-1}$, so $\forall^*h(G_xhG_x$ is comeager).

(iv) $\Rightarrow$ (v)  Fix $x\in\ooo$ and $h\in G$ such that
$G_xhG_x$ is comeager, so that $G_xhG_xh^{-1}$ is comeager and
$G_xG_y$ is comeager, where $y=h \cdot x\in\ooo$.  Since the map
$\pi (g)=g\cdot y$ is continuous and open from $G$ to $\ooo$, it
follows that $\pi (G_xG_y)=G_x\cdot y$ is comeager.

(v) $\Rightarrow$ (i) Fix $x,y\in\ooo$ with $G_x\cdot y$ comeager.
If $z=g \cdot x$, then $G_zg\cdot y=G_{g\cdot x}g\cdot y=gG_x\cdot
g^{-1}g\cdot y=gG_x\cdot y$ is comeager.  But $\{z\}\times G_zg\cdot
y=\{g\cdot x\} \times G_zg\cdot y\subseteq G\cdot (x,y)$, so
$\forall z\in \ooo\forall^* u\in X(z,u)\in G\cdot (x,y)$, thus by
the Kuratowski--Ulam Theorem, $G\cdot (x,y)$ is a comeager orbit in
$X^2$.\end{proof}

In particular, if $G$ is abelian, there cannot be a comeager orbit
in $X^2$, unless $X$ is a singleton.

Notice also that if there is a comeager orbit  in $X^2$ and $\ooo$ is the comeager orbit in $X$, then, by (iv), for all $x\in\ooo
,\forall^*h\in G(G_xhG_x$ is comeager).  But since the product of
two comeager sets in $G$ is equal to $G$, we have
\begin{equation}
\forall x\in\ooo\;\forall^*h\in G\;(G_xG_{h\cdot x}G_x=G_xhG_x\cdot
G_x h^{-1}G_x=G).
\end{equation}
\textbf{Remark.} Given a Polish space $X$ and $n=1,2,\dots$, there
is a natural continuous action of $S_n$ (the group of permutations
of $n$) on $X^n$ given by
\begin{equation}
\sigma\cdot (x_1,x_2,\dots ,x_n)=(x_{\sigma^{-1}
(1)},x_{\sigma^{-1} (2)}, \dots ,x_{\sigma^{-1} (n)}).
\end{equation}

If there is a comeager orbit $A$ in $X^n$ under the action of $G$,
this orbit will be invariant by permutation of the coordinates,
i.e., invariant under the action by $S_n$. Because for each
$\sigma\in S_n ,\sigma\cdot A$ will be comeager, so $\sigma\cdot
A\cap A\neq\emptyset$. Take some $(y_1,y_2,\dots
,y_n)\in\sigma\cdot A\cap A$ and $h\in G$ such that $(y_{\sigma
(1)},y_{\sigma (2)},\dots ,y_{\sigma (n)})=(h \cdot y_1, \dots,
h\cdot y_n)$. Then for any
\begin{equation}
(x_1,x_2,\dots , x_n)=(k\cdot y_1, \dots, k\cdot y_n)\in A,
\end{equation}
we have
\begin{equation}
(x_{\sigma (1)},x_{\sigma (2)},\dots ,x_{\sigma (n)})=(k\cdot
y_{\sigma (1)}, \dots, k \cdot y_{\sigma(n)})= (kh \cdot y_1,
\dots, kh\cdot y_n)\in A.
\end{equation}
\textbf{Remark.} Naturally one would like to know if there is any
Polish group $G$ acting continuously on a Polish space $X$, with
card$(X)>1$, for which there is a comeager orbit in
$X^{\mathbb{N}}$. This however turns out not to be the case. For
suppose card$(X) >1$ and let $\emptyset \neq V \subseteq X$ be
open, but not dense. Let $C \subseteq X^{\mathbb{N}}$ be a comeager
orbit, towards a contradiction. Then $A= \{(x_n) \in X^\bbN:
\exists^\infty n (x_n \in V)\}$ is dense and $G_\delta$, so fix $(x_n) \in A
\cap C$. Then $\Omega=\{ n \in \mathbb{N}: x_n \in V\}$ is
infinite. Put $B=\{(y_n)\in X^{\mathbb{N}}: \{y_n\}_{n\in \Omega}
\ {\rm is\ dense\ in}\ X\}$, which is a dense $G_\delta$. So fix $(y_n)
\in B\cap C$. Then there is $h \in G$ such that $h\cdot
(x_n)=(y_n)$, so that $\{h\cdot x_n\}_{n\in \Omega}\subseteq h
\cdot V$ is dense, a contradiction.

\medskip
If $G$ acts continuously on $X$, we call $(x_1,\dots, x_n) \in
X^n$ \emph{generic} if its orbit is comeager. We say that the
action has \emph{ample generics} if for each $n$, there is a
generic element in $X^n$. In particular, we say that a Polish
group $G$ has {\it ample generic} elements if there is a comeager
orbit in $G^n$ (with the conjugacy action of $G$ on $G^n$), for
each finite $n$.  If $\bar{g}\in G^n$ has a comeager orbit, we
will refer to it as a {\it generic element} of $G^n$.  This is an
abuse of terminology as $G^n$ is itself a Polish group and so it
makes sense to talk about a generic element of $G^n$ viewed as a
group.  Note that a generic element in the sense we use here is
also generic in the group $G^n$, so our current notion is
stronger, and it will be the only one we will use in the rest of
the paper. Note that $G^n$ as a group has a generic element iff
$G$ has a generic element.

Suppose $\kkk$ is a \fraisse class, and $\bfk$ its \fraisse limit.
We would like to characterize as before when $\bfk$ has ample
generic automorphisms.

We know that Aut$(\bfk )^n$ has a comeager diagonal conjugacy class iff
some $\bar f=(f_1,\dots ,f_n)\in$Aut$(\bfk )^n$ is turbulent and
has dense diagonal conjugacy class.  Therefore we can prove exactly as
before:

\begin{thm}\label{diagonal comeager} Let $\kkk$ be a \fraisse class and $\bfk$
its \fraisse limit.  Then the following are equivalent:

(i) There is a comeager diagonal conjugacy class in ${\rm Aut}(\bfk )^n$,

(ii) $\kkk^n_p$ has the {\rm JEP} and the {\rm WAP}.\end{thm}

So $\bfk$ has ample generic automorphisms iff for every $n$,
$\kkk^n_p$ has the JEP and the WAP.

We recall the results of Hodges et al. \cite{hodges2} stating that many
$\omega$-stable, $\aleph_0$-categorical structures and the random
graph have ample generic automorphisms, but this also holds for
many other structures.

\subsection{The Hrushovski property}

\begin{defi}\label{hrushovski property}We say that a Fraiss\'{e} class $\mathcal K$ satisfies the
\emph{Hrushovski property} if any system $\s =\langle \bfa , \psi_1
: \bfb_1 \rightarrow \bfc_1, \dots, \psi_n: \bfb_n \rightarrow
\bfc_n \rangle$ in $\kkk^n_p$ can be extended to some $\ttt
=\langle\bfd, \varphi_1: \bfd \rightarrow \bfd, \dots, \varphi_n:
\bfd \rightarrow \bfd\rangle$ in $\kkk^n_p$, i.e., to a sequence of
automorphisms of the same finite structure. \end{defi}

Let us first reformulate the Hrushovski property in topological
terms as a condition on the automorphism group.

\begin{prop}\label{approximate compactness}Let $\mathcal K$ be a \fraisse class (of finite structures), $\bfk$ its \fraisse limit, and ${\rm Aut}(\bfk)$ the
automorphism group of $\bfk$. Then $\kkk$ has the Hrushovski
property if and only if there is a countable chain $C_0\leq C_1\leq
C_2\leq \ldots\leq {\rm Aut}(\bfk)$ of compact subgroups whose union
is dense in ${\rm Aut}(\bfk)$.\end{prop}

\begin{proof}Suppose that $\kkk$ has the Hrushovski property and $G={\rm
Aut}(\bfk)$. Then each of the sets
\begin{displaymath}\begin{split}
F_n&=\{(f_1,\ldots,f_n)\in G^n\colon
\forall x\in \bfk\; x\text{'s orbit under } \langle f_1,\ldots,f_n\rangle \text{ is finite}\} \\
&=\{(f_1,\ldots,f_n)\in G^n\colon \langle f_1,\ldots,f_n\rangle \text{ is relatively compact in } G\}.
\end{split}\end{displaymath}
is comeager, and hence the generic infinite sequence $(f_i)$ in $G$
will generate a dense subgroup all of whose finitely generated
subgroups are relatively compact. In particular, $G$ is the closure
of the union of a countable chain of compact subgroups.

Conversely, suppose $C_0\leq C_1\leq C_2\leq \ldots\leq G$ is a
chain of compact subgroups whose union is dense in $G$. Let
$\bfa\subseteq \bfk$ be a finite substructure of $\bfk$ and let
$p_1,\ldots,p_n$ be a sequence of partial automorphisms of $\bfa$.
By density, we can find some $m$ such that there are
$f_1,\ldots,f_n\in C_m$ with $f_i\supseteq p_i$ for each $i\leq n$.
But then $B=C_m\cdot \bfa$ is a finite set and if we let $\bfb$ be
the finite substructure of $\bfk$ generated by $B$ then, as $C_m$ is
a group, one can check that $\bfb$ is invariant under elements of
$C_m$. In particular, $\bfb$ is closed under the automorphisms
$f_i$, whence $(\bfb, f_1|\bfb,\ldots,f_n|\bfb)$ is the extension of
the system $(\bfa,p_1,\ldots,p_n)$ needed.
\end{proof}

Hrushovski \cite{hrushovski} originally proved the Hrushovski property for the class $\mathcal K$ of finite graphs and this was
used in Hodges et al. \cite{hodges2} to show that the automorphism group of
the random graph has ample generic automorphisms. In fact, the
Hrushovski property simply means that the class of all $\mathcal T$ as
above is cofinal and this is often enough to show that $\kkk^n_p$,
has the CAP for all $n$. This, combined with the JEP for $\kkk^n_p$,
which is usually not hard to verify, implies the existence of
ample generics. For example, this easily works for the class of
finite graphs. Another case is the following. Solecki \cite{solecki} and,
independently, Vershik have shown that the class $\kkk$ of finite
metric spaces with rational distances satisfies the Hrushovski
property. From this it easily follows that $\kkk^n_p$ has the CAP
for each $n$. Take for notational simplicity $n=1$.

Let $\llll\subseteq\kkk_p$ be the class of systems $\s =\langle\bfa
,\psi :\bfa\rightarrow\bfa\rangle$, $\bfa \in \kkk$. By the
Hrushovski property of $\kkk$ this class is cofinal under
embeddability in $\kkk_p$.  We claim that $\llll$ has the AP. Let $\psi
: \bfa\rightarrow\bfa ,\phi :\bfb\rightarrow\bfb , \chi
:\bfc\rightarrow\bfc$ be in $\kkk_p$ with $\bfa\subseteq\bfb
,\bfa\subseteq\bfc$ and $\psi\subseteq\phi , \psi\subseteq\chi$.
Let $\partial$ be the metric on $\bfb$ and $\rho$ the metric on
$\bfc$ and suppose without loss of generality that $\bfb \cap\bfc
=\bfa$.

We define the following metric $d$ on $\bfb\cup\bfc$:
\begin{itemize}
\item for $x,y\in\bfb$, let $d(x,y)=\partial (x,y)$
\item for $x,y\in\bfc$, let $d(x,y)=\rho (x,y)$
\item for $x\in\bfb ,y\in\bfc$ let $d(x,y)=\min (\partial
(x,z)+\rho (z,y): z\in\bfa)$
\end{itemize}
One easily checks that this satisfies the triangle inequality.  So
now we only need to see that $\theta =\phi\cup\chi$ is actually an
automorphism of $\bfb\cup\bfc$.  Let us first notice that it is
well defined as $\phi$ and $\chi$ agree on their common domain
$\bfa$.  Moreover, trivially $d(x,y)=d(\theta (x),\theta (y))$,
whenever both $x,y\in\bfb$ or both $x,y\in\bfc$.  So let $x\in\bfb
,y\in\bfc$ and find $z\in\bfa$ such that $d(x,y)=d(x,z)+d(z,y)$.
Then
$$
d(\theta (x),\theta (y))\leq d(\theta (x),\theta (z))+d(\theta
(z),\theta (y))=d(x,z)+d(z,y)=d(x,y),
$$
as $\theta (z)\in\bfa$. And reasoning with $\theta^{-1}$ we get
$d(x,y)\leq d(\theta (x),\theta (y))$, so $\theta$ is indeed an
isometry of $\bfb\cup\bfc$ with the metric $d$.  This can
therefore be taken to be our amalgam.

Since the argument in the proof of \ref{dense urysohn} also shows that $\kkk^n_p$
has the JEP, for all $n$, it follows from \ref{diagonal comeager} that $\bfu_0$ has
ample generic automorphisms.

For a further example of a structure with the Hrushovski property
and ample generic automorphisms take $\kkk =\mba_\bbQ$, with
\fraisse limit $(\bff ,\lambda )$.  If
$$
\s =\langle\bfa ,\psi_1:\bfb_1\rightarrow \bfc_1,\dots
,\psi_n:\bfb_n\rightarrow\bfc_n\rangle
$$
is in $\kkk^n_p$, then it can be extended to some
$$
\ttt =\langle\bfd ,\phi_1: \bfd\rightarrow\bfd ,\dots
,\phi_n:\bfd\rightarrow\bfd\rangle,
$$
i.e., to a sequence of automorphisms of the same structure.
Moreover, $\bfd$ can be found such that all the atoms of $\bfd$ have
the same measure.

To see this, notice that, as the measure on $\bfa$ only takes
rational values, we can refine $\bfa$ to some $\bfd$ such that all
its atoms have the same measure.  But then, as $\psi_i$ preserves
the measure, for any $b\in \bfb_i, b$ is composed of the same
number of $\bfd$ atoms as $\psi_i(b)$. So $\psi_i$ can easily be
extended to an automorphism $\phi_i$ of $\bfd$.

Let $\llll^n$ be the subclass of $\kkk^n_p$ consisting of systems of
the same form as $\ttt$. We claim that $\llll^n$ has the AP. For
suppose that $\s =\langle \bfa ,\psi_1,\dots ,\psi_n\rangle$, $\ttt
=\langle\bfb ,\phi_1,\dots ,\phi_n \rangle$, $\rrr =\langle\bfc
,\chi_1 ,\dots ,\chi_n\rangle$ are in $\llll^n$ with $\s$ being a
subsystem of $\ttt$ and $\rrr$. List the atoms of $\bfa$ as
$a_1,\dots, a_p$, the atoms of $\bfb$ as
$$
b_1(1),\dots,b_1(k_1),\break \dots, b_p(1),\dots ,b_p(k_p),
$$
and the atoms of $\bfc$ as
$$
c_1(1),\dots ,c_1(l_1),\dots ,c_p(1),\dots ,\break c_p (l_p),
$$
where
$$
a_i=b_i(1)\vee\dots\vee b_i(k_i)=c_i(1)\vee\dots\vee c_i(l_i)
$$
Then we can amalgamate $\ttt$ and $\rrr$ over $\s$ by taking atoms
$b_i(e)\otimes c_i(d)$ and sending $b_i(e)$ to $\bigvee_{d\leq
l_i} b_i(e)\otimes c_i(d),\;c_i(d)$ to $\bigvee_{e\leq
k_i}b_i(e)\otimes c_i(d)$ and letting
$$\nu (b_i(e)\otimes c_i(d))=\tfrac{\delta (b_i(e))\gamma (c_i
(d))}{\mu (a_i)}$$ where $\mu ,\delta$ and $\gamma$ are the
measures on $\bfa ,\bfb$ and $\bfc$ respectively.  Moreover, let
$$
\theta_j(b_i(e)\otimes
c_i(d))=\phi_j(b_i(e))\otimes\chi_j(c_i(d)).
$$

One easily checks that this is indeed an amalgam of $\ttt$ and $\rrr$
over $\s$. As, in this case, all $n$-systems have a common
subsystem in $\llll^n$, AP for $\llll^n$ also implies JEP for $\kkk^n$,
and hence:

\begin{thm} Let $(\bff ,\lambda )$ be the \fraisse
limit of $\mba_\bbQ$.  Then $(\bff ,\lambda )$ has ample generic
automorphisms.\end{thm}

It is easy to check that the above works for $\mba_{\bbQ_2}$ as
well, so the group Aut(clop$(2^\bbN ),\sigma$) has ample generic
elements and so does the group of measure preserving
homeomorphisms of $2^\bbN$.

Another property that has been studied in the context of
automorphism groups is the existence of a dense locally finite
subgroup. Bhattacharjee and Macpherson \cite{bhattacharjee} showed
that such a group exists in the automorphism group of the random
graph and it is not difficult to see that also $H(2^\bbN , \sigma )$
has one.  Let us just mention that if ${\rm Aut}(\bfm )$, for $\bfm$
a (locally finite) \fraisse structure, has a dense locally finite
subgroup $H$, then $\bfm$ has the Hrushovski property. This follows
directly from Proposition \ref{approximate compactness}.

In particular, there is no dense locally finite subgroup of $H(2^\bbN
)$, since $\bfb_\infty$ does not even have the Hrushovski
property. Vershik \cite{vershik} poses the question of whether Aut$(\bfu_0 )$ has a
locally finite dense subgroup.

\subsection{Two lemmas} We will now prove two technical lemmas, which generalize
and extend some results in Hodges et al. \cite{hodges2}. The second will be
repeatedly used later on.

\begin{lemme}\label{mainsublemma} Let $G$ be a Polish group acting
continuously on a Polish space $X$ with ample generics. Let $A, B
\subseteq X$ be such that $A$ is not meager and $B$ is not meager
in any non-$\emptyset$ open set. Then if $\bar{x} \in X^n$ is
generic and $V$ is an open nbhd of the identity of $G$, there are
$y_0 \in A, y_1 \in B, h \in V$ such that $(\bar{x} , y_0),
(\bar{x} , y_1) \in X^{n+1}$ are generic and $h \cdot (\bar{x} ,
y_0)=(\bar{x} , y_1)$.\end{lemme}

\begin{proof} Let $C\subseteq X^{n+1}$ be a comeager orbit. Then, by
Kuratowski-Ulam, $\{\bar{z} \in X^n:\forall^\ast y(\bar{z}, y)\in
C\}$ is comeager and clearly $G$-invariant, so it contains
$\bar{x}$ and thus $\forall^\ast y (\bar{x} , y)\in C$. If $y\in
C_{\bar{x}}=\{z:(\bar{x} , z)\in C\}$, then
$C_{\bar{x}}=G_{\bar{x}}\cdot y$, where $G_{\bar{x}}$ is the
stabilizer of $\bar{x}$ for the action of $G$ on $X^n$. Thus for
any $y \in C_{\bar{x}}, G_{\bar{x}}\cdot y$ is comeager. Fix $y_0
\in A \cap C_{\bar{x}}$. Consider now the action of $G_{\bar{x}}$
on $X$. Since $G_{\bar{x}}\cdot y_0$ is comeager, it is
$G_\delta$, so by Effros' Theorem, the map
\[
\begin{array}{c}
\pi: G_{\bar{x}} \rightarrow G_{\bar{x}} \cdot y_0\\
g\mapsto g\cdot y_0
\end{array}
\]
is continuous and open. Thus $\pi (G_{\bar{x}} \cap V) =
(G_{\bar{x}} \cap V) \cdot y_0$ is open in $G_{\bar{x}} \cdot y_0$,
so $(G_{\bar{x}} \cap V)\cdot y_0 \cap B \neq \emptyset$. Fix then
$y_1 \in(G_{\bar{x}} \cap V)\cdot y_0 \cap B$. Then for some $h \in
G_{\bar{x}} \cap V, h\cdot y_0 = y_1$ and clearly $h \cdot (\bar{x},
y_0) = h\cdot(\bar{x}, y_1)$.\end{proof}
\medskip

One can also avoid the use of Effros' Theorem in the above proof and
instead give an elementary proof along the lines of Proposition
\ref{conditions turbulence}.

\begin{lemme}\label{mainlemma} Let $G$ be a Polish group acting
continuously on a Polish space $X$ with ample generics. Let $A_n,
B_n \subseteq X$ be such that, for each $n$, $A_n$ is not meager
and $B_n$ is not meager in any non-$\emptyset$ open set. Then
there is a continuous map $a\mapsto h_a$ from $2^{\bbN}$ into $G$,
such that if $a|n=b|n, a(n)=0, b(n)=1$, we have $h_a\cdot A_n \cap
h_b\cdot B_n \neq \emptyset$.\end{lemme}

\begin{proof} Fix a complete metric $d$ on $G$. For $s\in 2^{< \bbN}$ define
$f_{s^\wedge 0}, f_{s^\wedge 1} \in G, x_{s^\wedge 0}, x_{s^\wedge
1} \in X$ such that if $\bar{x}_s=(x_{s|1}, x_{s|2}, \dots, x_s),
h_s = f_{s|1} \ldots f_s \ (s\neq \emptyset)$, we have

\renewcommand{\labelenumi}{(\arabic{enumi})}

\begin{enumerate}
\item $\bar{x}_s$ is generic,
\item $x_{s^\wedge 0} \in A_{|s|},\  x_{s^\wedge 1} \in B_{|s|},$
\item $f_{s^\wedge 0} = 1_G,$
\item $d(h_s, h_{s} f_{s^\wedge 1}) < 2^{-|s|},$
\item $f_{s^\wedge 1}\cdot \bar{x}_{s^\wedge 1}=\bar{x}_{s^\wedge 0}.$
\end{enumerate}

We begin by using \ref{mainsublemma} to find $x_0, x_1, f_0, f_1$
(taking $h_\emptyset =1$).

Suppose $f_s, x_s$ are given. By \ref{mainsublemma} again, we can find
$x_{s^\wedge 0} \in A_{|s|}$, so that $\bar{x}_{s^\wedge 0},
\bar{x}_{s^\wedge 1}$ are generic and $g_{s^\wedge 1} \in
V=\{f:d(h_s, h_s f)<2^{-|s|}\}^{-1}$ with $g_{s^\wedge 1}\cdot
\bar{x}_{s^\wedge 0}=\bar{x}_{s^\wedge 1}$. Let $f_{s^\wedge
1}=g^{-1}_{s^\wedge 1}$.

By (3) and (4), for any $a\in 2^{\bbN}$ the sequence $(h_{a|n})$
is Cauchy, so converges to some $h_a \in G$. Now consider some $a\in \ca$ and $m<n$. If $a(n-1)=0$, then by (3), $f_{a|n}=1_G$ and hence
$$
f_{a|n}\cdot x_{a|m}=x_{a|m}.
$$
On the other hand, if $a(n-1)=1$, then by (5) we have
$$
f_{a|n}\cdot\bar x_{(a|n-1)^\wedge 1}=f_{(a|n-1)^\wedge 1}\cdot\bar x_{(a|n-1)^\wedge 1}=\bar x_{(a|n-1)^\wedge 0},
$$
whence, in particular,
$$
f_{a|n}\cdot\bar x_{a|n-1}=\bar x_{a|n-1}.
$$
But since $m<n$, $x_{a|m}$ is a term in $\bar x_{a|n-1}$ and thus also
$$
f_{a|n}\cdot x_{a|m}=x_{a|m}.
$$

Fix now $a\in \ca$. Then, as $h_{a|n}\Lim{n\til \infty}h_a$, we have for any $m$, $h_{a|n}\cdot x_{a|m}\Lim{n\til \infty}h_a\cdot x_{a|m}$.
Thus, for any $m$,
\begin{equation}
\begin{split}
h_a\cdot x_{a|m}&=\lim\limits_n h_{a|n} \cdot x_{a|m}\\
&=\lim\limits_{n}f_{a|1}\ldots f_{a|m}  \ldots f_{a|n}
\cdot x_{a|m}\\
&= f_{a|1}\ldots f_{a|m} \cdot x_{a|m}\\
&= h_{a|m}\cdot x_{a|m}.
\end{split}
\end{equation}

\noindent So if $a|n = b|n=s, a(n)=0, b(n) =1$, we have
\[
h_a \cdot x_{s^\wedge 0}=h_{s^\wedge 0} \cdot x_{s^\wedge 0} =
h_sf_{s^\wedge 0} \cdot x_{s^\wedge 0}= h_s \cdot x_{s^\wedge 0},
\]
as $f_{s^\wedge 0} = 1_G$, and
\[
\begin{split}
h_b \cdot x_{s^\wedge 1} = h_{s^\wedge 1} \cdot x_{s^\wedge
1}&=h_s f_{s^\wedge 1} \cdot x_{s^\wedge 1}=h_s \cdot x_{s^\wedge
0},
\end{split}
\]
so $h_a \cdot x_{s^\wedge 0} = h_b \cdot x_{s^\wedge 1}$, and, since
$x_{s^\wedge 0} \in A_n, x_{s^\wedge 1} \in B_n$, we have $h_a \cdot
A_n \cap h_b \cdot B_n \neq\emptyset$.\end{proof}

\subsection{The small index property}
We will next discuss the
connection of ample generics with the so-called {\it small index
property} of a Polish group $G$, which asserts that any subgroup
of index $<2^{\aleph_0}$ is open.

\begin{lemme}{\rm (Hodges et al. \cite{hodges2})}\label{meager subgroup} Let $G$ be a Polish
group.  Then any meager subgroup has index $2^{\aleph_0}$ in
$G$.\end{lemme}

\begin{proof} (Solecki) Notice that $G$ is perfect, i.e., has no isolated points, as otherwise it would be
discrete.  Let $E=\{(g,h)\in G^2: gh^{-1}\in H\}$, where $H\leq G$
is a meager subgroup.  Then as $(g,h)\mapsto gh^{-1}$ is continuous
and open, $E$ must be a meager equivalence relation and therefore by
Mycielski's Theorem (see Kechris \cite{kechris1} (19.1)) have
$2^{\aleph_0}$ classes. \end{proof}

\begin{thm}\label{small index} Let $G$ be a Polish group with ample
generics. Then $G$ has the small index property.\end{thm}

\begin{proof} Suppose $H\leq G$ has index $< 2^{\aleph_0}$ but is not open.
Then $H$ is not meager by \ref{meager subgroup}. Also $G\setminus H$
is not meager in any non-$\emptyset$ open set, since otherwise $H$
would be comeager in some non-$\emptyset$ open set and thus, by
Pettis' Theorem (see Kechris \cite{kechris1} (9.9)), $H$ would be
open. We apply now Lemma \ref{mainlemma} to the action of $G$ on
itself by conjugation and $A_n=H$, $B_n=G\setminus H$. Then, if
$a\neq b \in 2^{\bbN}$, with say $a|n=b|n, a(n) = 0, b(n) =1, h_a H
h_a^{-1} \cap h_b (G\setminus H) h_b^{-1}\neq \emptyset$, so
$(h_b^{-1} h_a) H (h_a^{-1} h_b)\cap (G\setminus H)\neq \emptyset$,
therefore $h_b^{-1} h_a \notin H$, thus $h_a, h_b$ belong to
different cosets of $H$, a contradiction.\end{proof}

We will now apply  this to the group of measure preserving
homeomorphisms of the Cantor space.

\begin{lemme} Suppose $\bfa, \bfb$ are two finite
subalgebras of {\rm clop}$(2^\bbN )$ and $G={\rm Aut(clop}(2^\bbN
),\sigma$).  Then $\langle G_{(\bfa )}\cup G_{(\bfb )}\rangle
=G_{(\bfa\cap\bfb )}$, where $G_{(\bfa )}$ is the pointwise
stabilizer of $\bfa$.\end{lemme}

\begin{proof}  Since $G_{(\bfa )}$ is open, so is $\langle G_{(\bfa )}\cup
G_{(\bfb )}\rangle$, and it is therefore also closed. Moreover, it
is trivially contained in $G_{(\bfa\cap\bfb )}$, so it is enough to
show that $\langle G_{(\bfa )}\cup G_{(\bfb )}\rangle$ is dense in
$G_{(\bfa \cap\bfb )}$.  Suppose $\bfd$ is a finite subalgebra of
clop$(2^\bbN )$ with a set of atoms $X$ all of which have the same
measure.  Suppose moreover, that $\bfa ,\bfb\subseteq\bfd$. Let
$\{a_1,\dots ,a_n\},\{b_1,\dots ,b_m\}, \{c_1,\dots ,c_k\}$ be the
partitions of $X$ given by the atoms of the subalgebras $\bfa ,\bfb
,\bfa\cap\bfb =\bfc$ of $\bfd =\p (X)$. Since all of the elements of
$X$ have the same measure, it is enough to show that any permutation
$\rho$ of $X$, pointwise fixing $\bfc$, i.e., $\rho''
c_i=c_i,i=1,\dots ,k$, is in $\langle G_{(\bf{A})} \cup
G_{(\bf{B})}\rangle$. In fact it is enough to show this for any
transposition $\rho$.  Let $\sigma_{x,y}$ denote the transposition
switching $x$ and $y$ (we allow here $x=y$). Fix $x_0\in X$ and let
$V_{x_0}\subseteq X$ be the set of atoms $y\in X$ such that
$\sigma_{x_0,y}\in\langle G_{(\bfa )}\cup G_{(\bfb )}\rangle$. Then
$V_{x_0}\in\bfc$. For if, e.g., $a_i\cap V_{x_0}\neq \emptyset$,
then there is some $y\in a_i\cap V_{x_0}$. So for any $z\in a_i$ we
have $\sigma_{x_0,y},\sigma_{y,z}\in\langle G_{(\bfa )}\cup G_{(\bfb
)} \rangle$, and hence
$\sigma_{x_0,z}=\sigma_{x_0,y}\circ\sigma_{y,z}\circ\sigma_{x_0,
y}\in\langle G_{(\bfa )}\cup G_{(\bfb )}\rangle$. This shows that
$V_{x_0}\in \bfa$ and a similar argument shows that $V_{x_0}\in
\bfb$. Therefore, if $\sigma_{x,y}$ is a transposition in Sym$(X)$
pointwise fixing $\bfc$, then $x,y$ belong to the same atom of
$\bfc$ and we have $y\in V_x$, showing that $\sigma_{x,y}$ $\in
\langle G_{(\bf{A})} \cup G_{(\bf{B})} \rangle$. \end{proof}

Let $G$ be a closed subgroup of $S_\infty$.  $G$ is said to have
the {\it strong small index property} if whenever $H\leq G$ and
$[G:H]<2^{\aleph_0}$ there is a finite set $X\subseteq\bbN$ such
that $G_{(X)}\leq H\leq G_{\{ X\}}$, where $G_{\{X\}}$ is the {\it
setwise stabilizer} of $X$.

\begin{thm} Let $G={\rm Aut(clop}(2^\bbN ),\sigma )
\cong H(2^\bbN ,\sigma )$. Then

(i) $G$ has ample generics,

(ii) $G$ has the strong small index property.
\end{thm}

\begin{proof}  We know that $G$ has ample generics, so has the small index
property.

Suppose $[G:H]<2^{\aleph_0}$.  Then $H$ is open and therefore for
some finite $Y\subseteq$ clop$(2^\bbN )$, we have $G_{(Y)}\leq H$.
If $\bfa$ is the subalgebra generated by $Y$, then obviously
$G_{(\bfa )}=G_{(Y)} \leq H$.  We notice that if $h\in H$ then
$G_{(h'' \bfa )}=hG_{(\bfa )}h^{-1}\leq H$ and by the lemma
$G_{(\bfa \cap h'' \bfa )}=\langle G_{(\bfa )}\cup G_{(h'' \bfa
)}\rangle\leq H$. So, as $\bfa$ is finite, we have, for
$\bfb=\bigcap_{h\in H}h'' \bfa$, that  $G_{(\bfb )}\leq H$. But
trivially also $H\leq G_{\{\bfb\}}$. So $G$ has the strong small
index property.
\end{proof}

\subsection{Uncountable cofinality}
In the same manner as in Hodges
et al. \cite{hodges2}, we can also obtain the following result (analogous to
their Theorem 6.1).

\begin{thm}\label{nonopen subgroups} Let $G$ be a Polish group with ample
generic elements.  Then $G$ is not the union of a countable chain
of non-open subgroups.\end{thm}

\begin{proof} Let $G=\bigcup_n A_n$, where $A_0 \subseteq A_1 \subseteq
\dots$ are non-open subgroups. We can then assume that each $A_n$
is non-meager. Also, as in the proof of \ref{small index}, each $G\setminus A_n$
is non-meager in every non-$\emptyset$ open set. Apply then Lemma
\ref{mainlemma} to $A_n$ and $B_n=G\setminus A_n$ to find $h_a$, $a\in 2^{\bbN}$, in
$G$, so that if $a|n=b|n, a(n)=0, b(n)=1$, then
$$(h_a A_n h_a^{-1})\cap
(h_b (G\setminus A_n) h_b^{-1} )\neq \emptyset,
$$
whence $h_a \notin A_n$ or $h_b \notin A_n$. Find uncountable
$C\subseteq 2^{\bbN}$ and $m$ so that $h_a \in A_m$ for all $a\in
C$. Then find $a, b\in C$ such that $h_a, h_b\in A_m$ and $a|n=b|n$,
$a(n)=0$, $b(n)=1$ for some $n>m$. As $h_a, h_b \in A_n$ we have a
contradiction. \end{proof}

Let us also mention here the following result of Cameron (see
Hodges et al. \cite{hodges2}): If $G$ is an oligomorphic closed subgroup of
$S_\infty$, i.e., if $G$ is the automorphism group of some
$\aleph_0$-categorical structure or equivalently $G$ has finitely
many orbits on each $\bbN^n$, then any open subgroup of $G$ is
only contained in finitely many subgroups of $G$. Recall also that
any open subgroup of a Polish group is actually clopen. So we have
that if $G$ is either connected or an oligomorphic closed subgroup
of $S_\infty$, then $G$ is not the union of a chain of open proper
subgroups.  And therefore if $G$ has ample generic elements it
must have {\it uncountable cofinality}, i.e., $G$ is not the union
of a countable chain of proper subgroups.

Notice also that if a Polish group $G$ is topologically finitely
generated, then it cannot be the union of a countable chain of
proper open subgroups, so, if it also has ample generics, then it
has uncountable cofinality.

To summarize: If $G$ is a Polish group with ample generics and one
of the following holds:
\begin{itemize}
\item[(i)] $G$ is an oligomorphic closed subgroup of $S_\infty$,
\item[(ii)] $G$ is connected,
\item[(iii)] $G$ is topologically finitely generated,
\end{itemize}
\noindent then $G$ has uncountable cofinality.

\subsection{Coverings by cosets}
One can also use the main lemma to prove an analog of the following
well-known group theoretical result due to B.H. Neumann \cite{neumann}: If a group $G$ is covered by
finitely many cosets $\{g_i H_i\}_{i\leq n}$, then for some $i$,
$H_i$ has finite index.

Before we consider the main result in this setting, let us briefly look at a simple fact.

\begin{prop}Suppose $G$ is a Polish group with a comeagre conjugacy class. Then the smallest number of
cosets of proper subgroups needed to cover $G$ is $\aleph_0$.
\end{prop}

\begin{proof} Assume $G=g_1H_1\cup\ldots\cup g_nH_n$. Then by Neumann's
lemma, some $H_i$ is of finite index in $G$. But then $\langle
H_i\rangle$ must contain some subgroup $N$ which is normal and of
finite index in $G$. As a subgroup of finite index cannot be meagre,
$N$ must intersect and therefore contain the comeagre conjugacy
class, and thus be equal to $G$. So $G=N\subseteq H_i$.\end{proof}

In the same manner one can show that if $G$ is partitioned into finitely many pieces, one of these pieces generates $G$. We should
mention that this result does not generalize to partitions of a generating set. For example, Ab\'ert \cite{abert} has shown that $\ku S_\infty$
is generated by two abelian subgroups, but evidently $\ku S_\infty$ is not abelian itself.

\begin{thm} \label{b.h.neumann}Let $G$ be a Polish group with ample
generics. Then for any countable covering of $G$ by cosets $\{g_i
H_i\}_{i\in \bbN}$, there is an $i$ such that $H_i$ is open and thus
has countable index.\end{thm}

The theorem clearly follows from the following more general lemma.

\begin{lemme} Let $G$ be a Polish group with ample
generics. If $\{k_i A_i\}_{i \in \bbN}$ is any covering of $G$,
with $k_i\in G$ and $A_i\subseteq G$, then for some $i$
$$
A_i^{-1} A_i A_i^{-1} A_i^{-1} A_i A_i^{-1} A_i A_i A_i^{-1} A_i
$$
contains an open nbhd of the identity.\end{lemme}

\begin{proof} First enumerate in a list $\{B_n\}_{n \in \bbN}$ all
non-meager $A_i$'s, so that each one of them appears infinitely
often.  Then $\bigcup_{l, n \in \bbN} k_l B_n = B$ is clearly comeager. Notice now that if for some $n$ the set $B_n^{-1} B_n B_n
B_n^{-1} B_n$ is comeager in some non-$\emptyset$ open set, then,
by Pettis' Theorem, we are done, so we can assume that
$C_n=G\backslash (B_n^{-1} B_n B_n B_n^{-1} B_n$) is not meager in
any non-$\emptyset$ open set, so, by Lemma \ref{mainlemma}, there is a
continuous map $a\mapsto h_a$ from $2^{\bbN}$ into $G$ such that
if $a|n = b|n, a(n)=0, b(n) =1$, then $h_a B_n h_a^{-1} \cap h_b
C_n h_b^{-1} \neq \emptyset$. Since $B_n \cap C_n =\emptyset$ it
also follows that $a\mapsto h_a$ is injective, so $K=\{h_a: a\in
2^{\bbN}\}$ is homeomorphic to $2^{\bbN}$.

The map $(g, h)\in G\times K\mapsto g^{-1} h \in G$ is continuous
and open so, since
$$
\bigcup_{l, n \in \bbN} k_l B_n = B
$$
is comeager in $G$, we have, by Kuratowski--Ulam,
$$
\forall^\ast g\in G \;\forall^\ast h\in K\;(g^{-1} h\in B),
$$
so fix $g\in G$ with $\forall^\ast h\in K(h\in g B)$. Then fix $l,
n \in \bbN$ such that, letting $g_0=gk_l$, the set
$$
\{ a\in 2^{\bbN}: h_a \in g_0 B_n\}
$$
is non-meager, so dense in some
$$
N_t =\{ a\in 2^{\bbN}: t \subseteq a\},
$$
$t\in 2^{<\bbN}$. Let then $m>|t|$ be such that $B_n =B_m$, and
let $a, b$ be such that $a|m = b|m, a(m)=0, b(m)=1$ and $h_a, h_b
\in g_0 B_n = g_0 B_m$, say $h_a =g_0 h_1, h_b =g_0 h_2$, with
$h_1, h_2 \in B_m$. Then
\begin{equation}
\begin{split}
h_b^{-1} h_a B_m h_a^{-1} h_b
&=h_2^{-1} g_0^{-1} g_0 h_1 B_m h_1^{-1}g_0^{-1}g_0 h_2\\
&=h_2^{-1} h_1 B_m h_1^{-1} h_2 \\
&\subseteq B_m^{-1} B_m B_m
B_m^{-1} B_m,
\end{split}
\end{equation}
which contradicts that $h_b^{-1} h_a B_m h_a^{-1} h_b \cap C_m \neq
\emptyset$.\end{proof}

\begin{cor}\label{p.m.neumann}Suppose $G$ is a Polish group
with ample generics and that $G$ acts on a set $X$. Then the following are equivalent:

(i) All orbits are of size $2^{\aleph_0}$.

(ii) All orbits are uncountable.

(iii) For every countable set $A\subseteq X$ there is a $g\in X$
such that $g\cdot A\cap A=\tom$.
\end{cor}

\begin{proof} (i)$\saa$(ii) is trivial.

(ii)$\saa$ (iii): Suppose $A=\{a_0,a_1,\ldots\}$ and that (iii)
fails for $A$. Then clearly, $G=\bigcup_{i,j\in \N}G_{a_i,a_j}$,
where $G_{a_i,a_j}=\{g\in G\del g\cdot a_i=a_j\}$. Since each
$G_{a_i,a_j}$ is a coset of the subgroup $G_{a_i}$, this implies
by Theorem \ref{b.h.neumann} that for some $i$,
$[G:G_{a_i}]\leq \aleph_0$, and hence $G\cdot a_i$ is countable.

(iii)$\saa$(i): Suppose some orbit $\ku O$ has cardinality strictly
smaller than the continuum. Then for any $x\in \ku O$, $[G:G_x]=|\ku
O|<2^{\aleph_0}$, so by the small index property of $G$, $G_x$ is
open and hence of countable index in $G$. Thus $|\ku
O|\leq\aleph_0$, and hence $A=\ku O$ contradicts (iii). \end{proof}

\subsection{Finite generation of groups} We will now investigate a
strengthening of uncountable cofinality that concerns the finite
generation of permutation groups, a subject recently originated in
Bergman \cite{bergman} and also studied in Droste-G\"obel \cite{droste2}, Droste-Holland \cite{droste3}.

\begin{defi} A group $G$ is said to have the {\rm
Bergman property} iff for each exhaustive sequence of subsets
$W_0\subseteq W_1\subseteq W_2\subseteq \ldots\subseteq G$, there
are $n$ and $k$ such that $W_n^k=G$.

If $G$ has the stronger property that for some $k$ and each
exhaustive sequence of subsets $W_0\subseteq W_1\subseteq
W_2\subseteq \ldots\subseteq G$, there is $n$ such that $W_n^k=G$,
we say that $G$ is $k$-{\rm Bergman}.\end{defi}

Bergman \cite{bergman} proved the $k$-Bergman property for $S_\infty$ and subsequently Droste
and G\"obel \cite{droste2} found a sufficient condition for certain
permutation groups to have this property.

We will see that ample generics also provide an approach to this
problem.

\begin{prop}\label{bergmanlemma} Let $G$ be a Polish group with ample
generic elements and suppose $A_0\subseteq
A_1\subseteq\dots\subseteq G$ is an exhaustive sequence of subsets
of $G$. Then there is an $i$ such that $1\in {\rm
Int}(A_i^{10})$.\end{prop}

\begin{proof} Notice first that also $A_0\cap A_0\inv\subseteq A_1\cap
A_1\inv\subseteq\ldots$ exhausts $G$. For given $g\in G$ find $m$
such that $g,g\inv\in A_m$; then $g\in A_m\cap A_m\inv$. So we can
suppose that each $A_n$ is symmetric. As $\{A_n\}_{n\in\bbN}$ is a
covering of $G$, there is, by Lemma \ref{mainlemma}, some $i$ such
that 
\begin{equation}
A_i^{10}=A_i^{-1} A_i A_i^{-1} A_i^{-1} A_i A_i^{-1} A_i A_i
A_i^{-1} A_i
\end{equation}
contains an open neighborhood of the identity.
\end{proof}

In the case of oligomorphic groups we have the following, whose
proof generalizes some ideas of Cameron:

\begin{thm} Suppose $G$ is a closed oligomorphic
subgroup of $S_\infty$ with ample generic elements. Then $G$ is
$21$-Bergman.\end{thm}

\begin{proof} Suppose that $W_0\subseteq W_1\subseteq W_2\subseteq
\ldots\subseteq G$ is an exhaustive chain of subsets of $G$. Then by
Proposition \ref{bergmanlemma}, there is an $n$ such that $W_n^{10}$
contains an open neighborhood of the identity. Find some finite
sequence $\overline a\in \bbN^{<\bbN}$, $|\overline a|=m$, such that
$G_{(\overline a)}\subseteq W_n^{10}$. Then as $G$ is oligomorphic
there are only finitely many distinct orbits of $G_{(\overline a)}$
on $\bbN^m$. Choose representatives $\ov b_1,\ldots,\ov b_k$ for
each of the $G_{(\overline a)}$ orbits on $\bbN^m$ that intersect
$G\cdot\ov a$ and find $h_1, \ldots, h_k\in G$ such that $h_i\cdot
\ov a=\ov b_i$. In other words, for each $f\in G$ there are $i\leq
k$ and $g\in G_{(\ov a)}$ such that $f\cdot \ov a=g\cdot\ov
b_i=gh_i\cdot\ov a$. Now find $l\geq n$ sufficiently big such that
$h_1,\ldots,h_k\in W_l$. We claim that $W_l^{21}=G$. Let $f$ be any
element of $G$ and find $i\leq k$ and $g\in G_{(\ov a)}$ with
$f\cdot \ov a=g\cdot\ov b_i=gh_i\cdot\ov a$. Then
$h_i^{-1}g^{-1}f\in G_{(\ov a)}$ and thus $f\in G_{(\ov
a)}W_lG_{(\ov a)}\subseteq W_l^{21}$, i.e., $W_l^{21}=G$.\end{proof}

Let us now make the following trivial but useful remarks:  Suppose
$\bfm$ is some countable structure and $G=$ Aut$(\bfm )$ has ample
generics.  If furthermore for any finitely generated substructure
$\bfa\subseteq\bfm$ there is a $g\in G$ such that $G=\langle
G_{(\bfa )}G_{(g''\bfa)}\rangle$, then $G$ has uncountable
cofinality.  This follows easily from Theorem \ref{nonopen subgroups}, as $G_{(g''\bfa
)}=gG_{(\bfa )}g^{-1}$.

If moreover there is a finite $n$ such that $G=(G_{(\bfa
)}G_{(g''\bfa )})^n$, then $G$ has the Bergman property. For
suppose $W_0\subseteq W_1\subseteq W_2\subseteq \ldots\subseteq G$
is an exhaustive chain of subsets of $G$. Then by Proposition
\ref{bergmanlemma}, there is an $m$ such that $W_m^{10}$ contains an open
neighborhood of the identity.  Say it contains $G_{(\bfa )}$ for
some finitely generated substructure $\bfa\subseteq\bfm$.  Find
$g\in G$ and $n$ as above. Then for $k\geq m$ big enough,
$g,g^{-1}\in W_k$ and 
\begin{equation}
G=(G_{(\bfa )}G_{(g''\bfa )})^n=(G_{(\bfa )}gG_{(\bfa )} g^{-1})^n\subseteq
(W_m^{10}W_kW_m^{10}W_k)^n\subseteq W_k^{22n}.
\end{equation}

This situation is less rare than one might think. In fact:

\begin{thm} The group $H(2^\bbN ,\sigma )$ of measure
preserving homeomorphisms of the Cantor space is
$32$-Bergman.\end{thm}

We will, as always, identify $H(2^\bbN ,\sigma )$ and
Aut(clop$(2^\bbN ), \sigma )$.
\medskip

\begin{proof} Recall that $2^\bbZ$ is measure preservingly homeomorphic to
$2^\bbN$, so we will temporarily work on the former space.  Let
$g$ be the Bernoulli shift on $2^\bbZ$ seen as an element of
Aut(clop$(2^\bbZ ),\sigma )$.

Suppose $\bfa$ is a finite subalgebra of clop$(2^\bbZ )$.  By
refining $\bfa$, we can suppose that $\bfa$ has atoms $N_s=\{x\in
2^\bbZ :x \restriction_{[-k,k]}=s\}$ for $s\in 2^{2k+1}$.  Then
$(g^{2k+1})''\bfa$ is independent of $\bfa$.  Let $h=g^{2k+1}$. We
wish to show that $G_{(\bfa )}G_{(h''\bfa )}G_{(\bfa )}=G$. Since
$G_{(\bfa )}$ is an open neighborhood of the identity, it is
enough to show that $G_{(\bfa )}G_{(h''\bfa )}G_{(\bfa )}$ is
dense in $G$, and for this it suffices to show that whenever
$\bfd$ is a finite subalgebra and $f_0\in G_{\{\bfd\}}$, the
setwise stabilizer of $\bfd$, then there is some $f_1\in G_{(\bfa
)}G_{(h''\bfa )}G_{(\bfa )}$ agreeing with $f_0$ on $\bfd$, i.e.,
$f_0\!\!\restriction_\bfd =f_1\!\!\restriction_\bfd$. We notice
first that by refining $\bfd$, we can suppose that $\bfd$ is the
subalgebra of clop$(2^\bbZ )$ generated by $\bfa$ and some finite
subalgebra $\bfb\supseteq h''\bfa$ that is independent of $\bfa$
and all of whose atoms have the same measure.  For concreteness,
we can suppose that $\bfb$ has atoms $N_{t,r}= \{x\in 2^\bbZ
:x\restriction_{[-l ,-k[}=t$ and $x\restriction_{]k,l ]} =r\}$,
for $t,r\in 2^{l -k}$, where $l$ is some number $>k$.  List the
atoms of $\bfa$ as $a_1,\dots ,a_n$ and the atoms of $\bfb$ as
$b_1,\dots ,b_m$. Then the atoms of $\bfd$, namely $a_i\cap b_j$,
all have the same measure and can be identified with formal
elements $a_i\otimes b_j, i\leq n,j\leq m$. So an element of
$G_{\{\bfd\}}$ gives rise to an element of $S=\ $Sym$(\{a_i\otimes
b_j:i\leq n,j\leq m\})$, whereas elements of $G_{(\bfa )} \cap
G_{\{\bfd\}}$ and $G_{(\bfb )}\cap G_{\{\bfd\}}$ give rise to
elements of $S$ that preserve respectively the first and the
second coordinates.

We now need the following well-known lemma, see, for example,
Ab\'ert \cite{abert} for a proof:

\begin{lemme} Let $\Omega =\{a_i\otimes b_j:i\leq n,j\leq
m\}$ and let $F\leq {\rm Sym}(\Omega )$ be the subgroup that
preserves the first coordinates and $H\leq {\rm Sym}(\Omega )$ be
the subgroup that preserves the second coordinates.  Then {\rm
Sym}$(\Omega )={\rm FHF}$.\end{lemme}

This finishes the proof.\end{proof}

Some of the results of this section are related to independent research by A. Ivanov \cite{ivanov2}. His setup is slightly different from ours as he formulates his results in terms of amalgamation bases.

\subsection{Actions on trees}
Macpherson and Thomas \cite{macpherson2} have recently found a relationship between the existence of a comeager
conjugacy class in a Polish group and actions of the group on trees.  This is further connected with Serre's property (FA) that
we will verify for the group of (measure preserving) homeomorphisms of $2^\bbN$.

A {\it tree} is a graph $T=(V,E)$ that is uniquely path connected, i.e., $E$ is a symmetric irreflexive relation on the set of
vertices $V$ such that any two vertices are connected by a unique path.

A group $G$ is said to {\it act without inversions} on a tree $T$ if there is an action of $G$ by automorphisms on $T$ such that for
any $g\in G$ there are no two adjacent vertices $a,b$ on $T$ such that $g\cdot a=b$ and $g \cdot b=a$.  Now we can state Serre's
property (FA):

A group $G$ is said to have {\it property} (FA) if whenever $G$ acts without inversions on a tree $T=(V,E)$, there is a vertex
$a\in V$ such that for all $g\in G, g\cdot a=a$.

We say that a free product with amalgamation $G=G_1*_AG_2$ is {\it trivial} in case one of the $G_i$ is equal to $G$.

For $G$'s that are {\it not} countable we have the following characterization of property (FA), from Serre \cite{serre}:  $G$
has property (FA) iff

(i) $G$ is not a non-trivial product with amalgamation,

(ii) $\bbZ$ is not a homomorphic image of $G$,

(iii) $G$ is not the union of a countable chain of proper
subgroups.

\begin{thm} {\rm (Macpherson-Thomas \cite{macpherson2})} Let $G$ be a
Polish group with a comeager conjugacy class. Then $G$ cannot be
written as a free product with amalgamation.\end{thm}

Moreover, it is trivial to see that if a Polish group has a
comeager conjugacy class then also (ii) holds.  In fact every
element of $G$ is a commutator.  For suppose $C$ is the comeager
conjugacy class and $g$ is an arbitrary element of $G$.  Then both
$C$ and $gC$ are comeager, so $C\cap gC\neq\emptyset$.  Take some
$h,f\in C$ such that $h=gf$.  Then for some $k\in G$ we have
$gf=h=kfk^{-1}$ and $g=kfk^{-1}f^{-1}$.  So $G$ has only one
abelian quotient, namely $\{e\}$.  For if for some $H\unlhd G
,G/H$ is abelian, then $G=[G,G]\leq H$.

So to verify whether a Polish group with a comeager conjugacy
class has property (FA), we only need to show it has uncountable
cofinality.

Another way of approaching property (FA) is through the Bergman property, which one can see is actually very strong. For example, one can show that it implies that any action of the group by isometries on a metric space has bounded orbits. But in the case of isometric actions on
real Hilbert spaces or automorphisms of trees, having a bounded orbit is enough to ensure that there is a fixed point. So groups
with the Bergman property automatically have property (FA) and property (FH) (the latter says that any action by isometries on a
real Hilbert space has a fixed point). Thus in particular, ${\it H}(2^\N,\sigma)$ has both properties (FA) and (FH).

\subsection{Generic freeness of subgroups}\label{generic freeness}
Let us next mention another application of the existence of dense diagonal conjugacy classes in each $G^n$. Suppose $G$ is Polish and has a dense diagonal conjugacy class in $G^n$, for each $n\in\bbN$. Then there is a dense $G_\delta$ subset $C\subseteq G^\bbN$ such that any two sequences
$(f_n)$ and $(g_n)\in C$ generate isomorphic groups, i.e., the mapping $f_n\mapsto g_n$ extends to an isomorphism of $\langle
f_n\rangle$ and $\langle g_n\rangle$.

To see this, suppose $w(X_1,\dots ,X_n)$ is a reduced word.  Then either $w(g_1,g_2,\dots ,g_n)=e$, for all $g_1, \dots ,g_n$, or
$w(g_1,g_2,\dots ,g_n)\neq e$ for an open dense set of $(g_1,\dots ,g_n)\in G^n$.  (This is trivial as
$$
w(g_1,g_2,\dots ,g_n)=e\Leftrightarrow w(kg_1k^{-1},kg_2k^{-1},\dots ,kg_nk^{-1})=e,
$$
for all $k\in G$ and $(g_1,\dots ,g_n)\in G^n$.)  In any case, taking
the intersection over all reduced words, we have that there is a dense $G_\delta$ set $C\subseteq G^\bbN$ such that any two
sequences $(f_n)$ and $(g_n)$ in $C$ satisfy the same equations $w(X_1,\dots ,X_n)=e$, i.e., $f_n\mapsto g_n$ extends to an
isomorphism.

In particular, if for any nontrivial reduced word $w(X_1,\dots ,X_n)$ there are $g_1,\dots ,g_n\in G$ such that $w(g_1,\dots ,g_n)\neq e$,
then any sequence in $C$ freely generates a free group and, using the Kuratowski--Mycielski Theorem (see Kechris \cite{kechris1} (19.1)) the
generic compact subset of $G$ also freely generates a free group.

Macpherson \cite{macpherson1} shows that any oligomorphic closed subgroup of
$S_\infty$ contains a free subgroup of infinite rank. It is not
hard to see that it also holds for the groups ${\rm Aut}([0,1], \lambda
)$, $H(2^\bbN ,\sigma )$, ${\rm Aut}(\bfuo )$ and $\rm Iso(\bfu)$.
Moreover, $H(2^\bbN)$ and these latter groups also have dense
conjugacy classes in each dimension (for ${\rm Aut}(\bfuo)$ and
$H(2^\bbN)$ it is enough to use the multidimensional version of
Theorem \ref{criterion for density}, analogous to Theorem \ref{diagonal comeager}).

This property has been studied by a number of authors (see e.g.
Gartside and Knight \cite{gartside}) and it seems to be a fairly common
phenomenon in bigger Polish groups.

\subsection{Automatic continuity of homomorphisms}
We now come to the study of automatic continuity of homomorphisms from Polish
groups with ample generics. Notice that the small index property
can be seen as a phenomenon of automatic continuity. In fact, if a
topological group $G$ has the small index property, then any
homomorphism of $G$ into $S_\infty$ will be continuous. For the
inverse image of a basic open neighborhood of $1_{S_\infty}$ will
be a subgroup of $G$ with countable index and therefore open. But
of course this puts a strong condition on the target group, namely
that it should have a neighborhood basis at the identity
consisting of open subgroups. We would like to have some less
restrictive condition on the target group that still insures
automatic continuity of any homomorphism from a Polish group with
ample generics. Of course, some restriction is necessary, for the
identity function from a Polish group into itself, equipped with
the discrete topology, is never continuous unless the group is
countable.

\begin{lemme} Let $H$ be a topological group and $\kappa$
a cardinal number. Then the following are equivalent:

(i) For each open neighborhood $V$ of $1_H$, $H$ can be covered by
$<\kappa$ many right translates of $V$.

(ii) For each open neighborhood $V$ of $1_H$, there are not
$\kappa$ many disjoint right translates of $V$.\end{lemme}

\begin{proof} $(i)\Rightarrow(ii)$: Suppose $(ii)$ fails and $\{Vf_\xi
\}_{\xi<\kappa}$ is a family of $\kappa$ many disjoint right
translates of some open neighborhood $V$ of $1_H$. By choosing $V$
smaller, we can suppose that $V$ is actually symmetric. Then if
$\{Vg_\xi \}_{\xi<\lambda}$ covers $H$ for some $\lambda<\kappa$,
there are $f_\xi, f_\zeta$ and $g_\nu$, $\xi\neq \zeta$, such that
$f_\xi, f_\zeta\in Vg_\nu $, and hence $g_\nu\in Vf_\xi \cap
Vf_\zeta $, contradicting that $Vf_\xi \cap Vf_\zeta =\emptyset$.
So (i) fails.

$(ii)\Rightarrow(i)$: Suppose that (ii) holds and $V$ is some open
neighborhood of $1_H$. Find some open neighborhood $U\subseteq V$ of
$1_H$ such that $U^{-1}U\subseteq V$ and choose a maximal family of
disjoint right translates $\{Uf_\xi \}_{\xi<\lambda}$ of $U$. By
(ii), $\lambda<\kappa$. Suppose that $g\in H$. Then there is a
$\xi<\lambda$ such that $Ug \cap Uf_\xi \neq \emptyset$, whereby
$g\in U^{-1}Uf_\xi \subseteq Vf_\xi $. So $\{Vf_\xi
\}_{\xi<\lambda}$ covers $H$ and (i) holds. \end{proof}

Recall that the {\it Souslin number} of a topological space is the
least cardinal $\kappa$ such that there is no family of $\kappa$
many disjoint open subsets of the space. By analogy, for a
topological group $H$, let the {\it uniform Souslin number} be the
least cardinal $\kappa$ satisfying the equivalent conditions of
the above lemma.

Notice that the uniform Souslin number is at most the Souslin
number, which is again at most density$(H)^+$, where the density
of a topological space is the smallest cardinality of a dense
subset. In particular, if $H$ is a separable topological group,
then its (uniform) Souslin number is at most $\aleph_1\leq
2^{\aleph_0}$. We should mention that groups with uniform Souslin
number at most $\aleph_1$ have been studied extensively under the
name $\aleph_0$-bounded groups (see the survey article by
Tkachenko \cite{tkachenko}). A well-known result (see, e.g., Guran \cite{guran}) states that a
Hausdorff topological group is $\aleph_0$-bounded if and only if it (topologically) embeds
as a subgroup into a direct product of second countable groups. $\aleph_0$-bounded groups are easily seen to contain the
$\sigma$-compact groups. Moreover, the uniform Souslin number is productive in contradistinction to separability, i.e., any direct
product of groups with uniform Souslin number at most $\kappa$ has uniform Souslin number at most $\kappa$ (for any infinite
$\kappa$).

\begin{thm}\label{automatic} Suppose $G$ is a Polish group with ample
generic elements and $\pi:G\rightarrow H$ is a homomorphism into a
topological group with uniform Souslin number at most
$2^{\aleph_0}$ (in particular, if $H$ is separable). Then $\pi$ is
continuous.\end{thm}

\begin{proof} It is enough to show that $\pi$ is continuous at $1_G$. So let
$W$ be an open neighborhood of $1_H$. We need to show that
$\pi^{-1}(W)$ contains an open neighborhood of $1_G$. Pick a
symmetric open neighborhood  $V$ of $1_H$ such that
$V^{20}\subseteq W$ and put $A=\pi^{-1}(V^{-1}V)=\pi^{-1}(V^2)$.
\medskip

{\bf Claim 1.} {\it $A$ is non-meager.}

\begin{proof} Otherwise, as $(g,h)\mapsto gh^{-1}$ is open and continuous
from $G^2$ to $G$, there is by the Mycielski theorem a
Cantor set $C\subseteq G$, such that for any $g\neq h$ in $C$,
$gh^{-1}\notin A$. But
$$V\pi(g)\cap V\pi (h)=\emptyset\Leftrightarrow \pi(gh^{-1})\notin
V^{-1}V\Leftrightarrow gh^{-1}\notin A.$$ So this means that there
are continuum many disjoint translates of $V$, contradicting that
the uniform Souslin number of $H$ is at most $2^{\aleph_0}$.
\end{proof}

{\bf Claim 2.} {\it $A$ covers $G$ by $<2^{\aleph_0}$ many right
translates.}

\begin{proof} By the condition on the uniform Souslin number we can find a
covering $\{Vf_\xi \}_{\xi<\lambda}$ of $H$ by
$\lambda<2^{\aleph_0}$ many right translates of $V$. So for each
$Vf_\xi $ intersecting $\pi (G)$ take some $g_\xi\in G$ with $\pi
(g_\xi)\in Vf_\xi$. Then $Vf_\xi \subseteq VV^{-1}\pi(g_\xi)
=V^2\pi(g_\xi)$, so the latter cover $\pi (G)$. Now, if $g\in G$
find $\xi<\lambda$ such that $\pi(g)\in V^2\pi (g_\xi)$, whence $\pi
(gg_\xi^{-1})\in V^2$, i.e., $gg_\xi^{-1} \in A$ and $g\in Ag_\xi$.
So the $Ag_\xi$ cover $G$. \end{proof}

{\bf Claim 3.} {\it $A^5$ is comeager in some non-$\emptyset$ open
set.}

\begin{proof} Otherwise, by Lemma \ref{mainlemma}, we can find $h_a, a\in 2^{\bbN}$, in
$G$ so that if $a|n = b|n, a(n)=0, b(n)=1$,
$$h_a Ah_a^{-1} \cap h_b (G \setminus A^5)h_b^{-1} \neq \emptyset$$
\noindent or equivalently
$$h_b^{-1} h_a A h_a^{-1} h_b \cap (G \setminus A^5)\neq \emptyset.$$

Since $A$ covers $G$ by $<2^{\aleph_0}$ right translates, and thus
left translates (as $A$ is symmetric), there is uncountable $B
\subseteq 2^{\bbN}$ and $g\in G$ such that for $a\in B,  h_a \in
gA$. If $a, b\in B, a|n =b|n, a(n)= 0, b(n)=1$, then let $g_a, g_b
\in A$ be such that $h_a=g g_a, h_b=g g_b$. Then
$$
\begin{array}{l}
h_b^{-1} h_a A h_a^{-1} h_b=g_b^{-1} g^{-1} g g_a A g_a^{-1}
g^{-1} g g_b =g_b^{-1}g_a A g_a^{-1} g_b \subseteq A^5,
\end{array}
$$
a contradiction.\end{proof}

So, by Pettis' theorem, $A^{10}\subseteq \pi^{-1}(W)$ contains an
open neighborhood of $1_G$ and $\pi$ is continuous.\end{proof}

\begin{cor} Suppose $G$ is a Polish group with
ample generic elements. Then $G$ has a unique Polish group
topology.\end{cor}

\begin{proof} Suppose that $\tau$ and $\sigma$ are two Polish group
topologies on $G$ such that $G$ has ample generics with respect to
$\tau$. Then the identity mapping from $G, \tau$ to $G, \sigma$ is
continuous, i.e., $\sigma\subseteq \tau$. But then $\tau$ must be
included in the Borel algebra generated by $\sigma$ and in
particular the identity mapping is Baire measurable from $G,\sigma$
to $G, \tau$, so continuous. \end{proof}
\medskip

We see from the above proof, that the only thing we need is that
the two topologies be inter-definable, which is exactly what
Theorem \ref{automatic} gives us.

As one can easily show that $S_\infty$ has ample generics, our
result applies in particular to this group. So this implies that,
e.g., any unitary representation of $S_\infty$ on separable Hilbert space is actually a
continuous unitary representation. Moreover, whenever $S_\infty$
acts by homeomorphisms on some locally compact Polish space or by
isometries on some Polish metric space, then it does so
continuously. For these results it is enough to notice that the actions in question correspond to homomorphisms into the unitary group, resp. the homeomorphism and the isometry group, which are Polish when the spaces are separable, resp. locally compact.

In a beautiful paper Gaughan \cite{gaughan} proves that any Hausdorff group
topology on $S_\infty$ must extend its usual Polish topology. So
coupled with the above result this gives us the following rigidity
result for $S_\infty$ (we would like to thank V. Pestov for suggesting how to get rid of a Hausdorff condition in a previous version of the result):

\begin{thm} $S_\infty$ has exactly two separable
group topologies, namely the trivial one and the usual Polish
topology.\end{thm}

\begin{proof} Suppose that $\tau$ is a separable group topology on
$S_\infty$ and define
\begin{equation}
N=\bigcap\{U: U \;\textrm{ open } \;\&\; 1\in U\}
\end{equation}
We easily see that $N$ is conjugacy invariant, $N=N^{-1}$ and that
$N$ is closed under products. For if $x,y\in N$ and $W$ is any
open neighborhood of $xy$, then by the continuity of the group
operations, we can find open sets $x\in U$ and $y\in V$ such that
$UV\subseteq W$. But then as $1$ cannot be separated from $x$ by an open set, $x$ cannot be
separated from $1$ either. Thus $1\in U$ and similarly $1\in V$, whence
$1\in W$. So $xy$ cannot be separated from $1$ by an open set and
hence $1$ cannot be separated from $xy$ either. Thus $xy\in N$.

So $N$ is a normal subgroup of $S_\infty$ and hence equal to either
$\{1\}$, ${\rm Alt}$, ${\rm Fin}$ or $S_\infty$ itself. In the first
case, we see that $\tau$ is Hausdorff and thus that it extends the
Polish topology and in the last case that $\tau=\{\emptyset,
S_\infty\}$. Moreover, by the separability of $\tau$, we know that
$\tau$ is weaker than the Polish topology on $S_\infty$ and thus in
the first case we know that it must be exactly equal to the Polish
topology. We are therefore left with the two middle cases that we
claim cannot occur. The subgroups ${\rm Alt}$ and ${\rm Fin}$ are
dense in the Polish topology and thus also dense in $\tau$. If
$x\notin N$, then we can find some open neighborhood $V$ of $x$ not
containing $1$. But then $V\cap N=\emptyset$, because any element of
$V$ can be separated form $1$ and thus does not belong to $N$. This
shows that $N$ is closed and contains a dense subgroup, so
$N=S_\infty$, a contradiction. \end{proof}

We should also mention the following result that allows us to see the small index property as a special case of automatic continuity.

\begin{prop} \label{Cont-sip}Suppose $G$ is a Polish group such that any homomorphism from $G$ into a group with uniform Souslin
number $\leq 2^{\aleph_0}$ is continuous. Then $G$ has the small index property.\end{prop}

\begin{proof} Suppose $H$ is a subgroup of $G$ of small index. Then ${\rm Sym}(G/H)$, where $G/H$ is the set of left cosets of $H$, has
uniform Souslin number $\leq 2^{\aleph_0}$ and clearly the action of
left translation of $G$ on $G/H$ gives rise to a homomorphism of $G$
into ${\rm Sym}(G/H)$. By automatic continuity, this shows that the
pointwise stabiliser of the coset $H$ is open in $G$, i.e., $H$ is
open in $G$.\end{proof}

\subsection{Automorphisms of trees}
We will finally investigate the structure of the group of Lipschitz
homeomorphisms of the Baire space $\n$. This group is of course
canonically isomorphic to Aut$(\bbN^{<\bbN})$, where $\bbN^{< \bbN
}$ is seen as a copy of the uniformly countably splitting rooted
tree. Though Age$(\bbN^{<\bbN})$ is not a \fraisse class in its usual relational language, we can
still see $\bbN^{<\bbN}$ as being the generic limit of the class
of finite rooted trees and we will see that the theory goes
through in this context. Alternatively, one can change the language by replacing the tree relation by a unary function symbol that assigns to each node its predecessor in the tree ordering. In this way, Age$(\bbN^{<\bbN})$ outright becomes a Fra\"iss\'e class.

In the following, $\bbN^{<\bbN}$ is considered a tree with root
the empty string, $\emptyset$, such that the children of a vertex
$s\in \bbN^{<\bbN}$ are $s\:\hat{}\: n$, for all $n\in\bbN$.  A
subtree of $\bbN^{<\bbN}$ is a subset $T\subseteq\bbN^{<\bbN}$
closed under initial segments, i.e., if $\langle n_0,n_1,\dots
,n_k\rangle\in T$, then so are $\emptyset , \langle n_0\rangle
,\dots ,\langle n_0,n_1,\dots ,n_{k-1}\rangle$. By $m^{\leq m}$ we
denote the tree of sequences $\langle n_0,\dots , n_k\rangle$ such
that $n_i<m$ and $k<m$.  Moreover, if $s,t\in \bbN^{<\bbN }$, we
write $s\subseteq t$ to denote that $t$ extends $s$ as a sequence.

\begin{lemme}\label{hrushovski for trees} Suppose
$\phi :T\rightarrow S$ is an isomorphism between finite subtrees
of $\bbN^{<\bbN}$, say $T,S\subseteq m^{\leq m}$.  Then there is
an automorphism $\psi$ of $m^{\leq m}$ extending
$\phi$.\end{lemme}

\begin{proof}  Notice first that $\phi$ restricts to a bijection between two
subsets of $m^1=\{\langle 0\rangle ,\langle 1\rangle ,\dots ,\langle
m-1\rangle\}$, so can be extended to a permutation $\phi_1$ of
$m^1$. Then $\phi_1$ is an isomorphism of $T_1=T\cup m^1$ and
$S_1=S\cup m^1$.  Now suppose $\phi_1 (\langle 0\rangle )=\langle
j\rangle$. Then $\phi_1$ restricts to a bijection between two
subsets of $\{\langle 0,0\rangle ,\langle 0,1\rangle ,\dots ,\langle
0,m-1\rangle\}$ and $\{ \langle j,0\rangle ,\langle j,1\rangle
,\dots ,\langle j,m-1\rangle\}$ and can be extended as before to
some isomorphism of 
\begin{equation}
T_2=T_1\cup\{\langle 0,0\rangle ,\langle
0,1\rangle ,\dots ,\langle 0,m-1\rangle\}
\end{equation}
and 
\begin{equation}
S_2=S_1\cup\{\langle j,0\rangle ,\langle j,1\rangle ,\dots ,\langle j,
m-1\rangle\}
\end{equation}  
Now continue with $\phi_1(\langle 1\rangle )=\langle
\ell \rangle$, etc.  Eventually, we will obtain an automorphism
$\psi$ of $m^{\leq m}$ that extends $\phi$.\end{proof}

\begin{lemme} Suppose $m\leq n$ and $(\psi_1,\dots,\psi_k)$,
$(\phi_1,\dots ,\phi_k)$ are sequences of automorphisms of
$n^{\leq n}$ extending automorphisms $(\chi_1,\dots ,\chi_k)$ of
$m^{\leq m}$, i.e., $\chi_i\subseteq\psi_i,\phi_i$.  Then there is
an automorphism $\xi$ of $\ell^{\leq\ell}$, for $\ell =2n$, such
that $(\psi_1,\dots ,\psi_k)$ and $(\xi\circ\phi_1 \circ\xi^{-1},\dots
,\xi\circ\phi_k\circ\xi^{-1})$ can be extended to some common
sequence of automorphisms of $\ell^{\leq \ell}$. Moreover, $\xi$
fixes $m^{\leq m}$ pointwise.\end{lemme}

\begin{proof}  Let $\xi$ be an automorphism of $\ell^{\leq\ell}$ pointwise
fixing $m^{\leq m}$ and such that $\xi (n^{\leq n})\cap n^{\leq
n}=m^{\leq m}$.  Then for each $i\leq k,\psi_i$ and
$\xi\circ\phi_i\circ\xi^{-1}$ agree on their common domain, and the
same for $\psi^{-1}_i$ and $(\xi\circ\phi_i\circ\xi^{-1})^{-1}$. So
$\psi_i\cup\xi\circ\phi_i\circ\xi^{-1}$ is an isomorphism of finite
subtrees of $\ell^{\leq\ell}$ and the result follows from Lemma
\ref{hrushovski for trees} \end{proof}
\medskip

Though $\bbN^{<\bbN}$ is not ultrahomogeneous in its relational language, it is in the functional language, and the preceeding results show that its age in the functional language satisfies the conditions of Theorem \ref{diagonal comeager}. Therefore ${\rm Aut}(\bbN^{<\bbN})$ has ample generics.

\begin{lemme} \label{strong small index for trees} Suppose $S$ and
$T$ are finite subtrees of $\bbN^{<\bbN}$ and $G={\rm
Aut}(\bbN^{<\bbN})$.  Then $G_{(S\cap T)}
=G_{(T)}G_{(S)}G_{(T)}$.\end{lemme}

\begin{proof}  Since $G_{(T)}G_{(T)}G_{(S)}G_{(T)}=G_{(T)}G_{(S)}G_{(T)}$,
and $G_{(T)}$ is an open neighborhood of the identity, it is
enough to show that $G_{(T)}G_{(S)}G_{(T)}$ is dense in $G_{(S\cap
T)}$. By Lemma \ref{hrushovski for trees}, it is enough to show
that if $\phi$ is an automorphism of $m^{\leq m}$ for some $m$
such that $S,T\subseteq m^{\leq m}$ and $\phi$ pointwise fixes
$S\cap T$, then there is a $g\in G_{(T)}G_{(S)}G_{(T)}$ with
$g\supseteq\phi$.

For $s\in S\cap T$, let $A_s=\{n\in\bbN :s\:\hat{}\: n\not\in
S\cap T$ and $s\:\hat{}\: n\in S\}$ and fix a permutation
$\sigma_s$ of $\bbN$ pointwise fixing $\{0,1,2,\dots
,m-1\}\setminus A_s$, but such that $\sigma_s(A_s)\cap
A_s=\emptyset$. Define $f\in G_{(T)}$ as follows:  If $u\in S\cap
T$, let $f(u)=u$. Otherwise if 
\begin{equation}
u=s\:\hat{}\: n\:\hat{}\: t, s\in
S\cap T, s\:\hat{}\: n\not\in S\cap T, t\in\bbN^{<\bbN},
\end{equation}
let $f(u)=s\:\hat{}\: \sigma_s(n)\:\hat{}\: t$.

Now let $g\supseteq\phi$ be defined by:  If $u=s\:\hat{}\: t$, where
$s$ is the maximal initial segment such that $s\in m^{\leq m}$, let
$g(u)= \phi (s)\:\hat{}\: t$.  Then $f^{-1}gf\in G_{(S)}$.  For
suppose $u=s\: \hat{}\: n\:\hat{}\: t\in S, s\:\hat{}\: n\not\in T,
s\in S\cap T$. Then 
\begin{equation}
f^{-1}\circ g\circ f(u)=f^{-1}\circ g
(s\:\hat{}\: \sigma_s(n)\:\hat{}\: t)=f^{-1}(\phi (s)\:\hat{}\:
\sigma_s(n)\:\hat{}\: t)=f^{-1}(s\:\hat{}\: \sigma_s(n)\:\hat{}\:
t)=u,
\end{equation}
since $\phi$ fixes $S\cap T$ pointwise.  And if $u\in S\cap
T$, then $f^{-1}\circ g\circ f(u)=u$, as both $f,g\in G_{(S\cap
T)}$.  So 
\begin{equation}
g=ff^{-1}gff^{-1}\in G_{(T)}G_{(S)}G_{(T)}.
\end{equation} \end{proof}

\begin{thm} Let $G={\rm Aut}(\bbN^{<\bbN})$.  Then:

(i) $G$ has ample generic elements.

(ii) {\rm (R. M\"oller \cite{moller})} $G$ has the strong small index property.

(iii) $G$ is $32$-Bergman.

(iv) $G$ has a locally finite dense subgroup.
\end{thm}

\begin{proof}  (i) has been verified and so $G$ has the small index
property. Suppose $G_{(S)}\leq H\leq G$ for some finite subtree
$S\subseteq \bbN^{<\bbN}$ and open subgroup $H$.  Then, by Lemma
\ref{strong small index for trees}, we have, for any $h\in H$,
that $G_{(S\cap h''S)}=G_{(S)}G_{(h''S)}G_{(S)}\leq H$. So, as $S$
is finite, we have, if $T=\bigcap_{h\in H}h''S$, that $G_{(T)}
\leq H$ and also $H\leq G_{\{T\}}$.  So $G$ has the strong small
index property.

Now suppose $W_0\subseteq W_1\subseteq \ldots\subseteq G$ is an
exhaustive sequence of subsets. Then, by Proposition \ref{bergmanlemma}, there
is an $n\in\bbN$ such that $W_n^{10}$ contains an open subgroup
$G_{(S)}$, where $S\subseteq\bbN^{<\bbN}$ is some finite tree.
Take $g\in G$ such that $g''S\cap S=\{\emptyset\}$. Then for
$m\geq n$ sufficiently big we have $g,g^{-1}\in W_m$ and $G=G_{\{
\emptyset\}}=G_{(S)}G_{(g''S)}G_{(S)}=G_{(S)}gG_{(S)}g^{-1}G_{(S)}
\subseteq W_m^{32}$.  This verifies the Bergman property.

Finally, the following is easily seen to be a dense locally finite
subgroup of $G: K=\{g\in G:\exists m\exists\phi$ automorphism of
$m^{\leq m}$ such that if $u=s\:\hat{}\: t\in\bbN^{<\bbN}$, where
$s$ is the maximal initial segment of $u$ such that $s\in m^{\leq
m}$, then $g(u)=\phi (s)\:\hat{}\: t\}$.\end{proof}

We denote by $\bft$ the $\aleph_0$-regular tree on $\N$, i.e., the tree in which each vertex has valency $\aleph_0$. Notice that ${\rm
Aut}(\N^{<\N})\iso {\rm Aut}(\bft,a_0)={\rm Aut}(\bft)_{a_0}$, for any $a_0\in \bft$. Hence, ${\rm Aut}(\N^{<\N})$ is isomorphic to a clopen
subgroup of ${\rm Aut}(\bft)$.

\begin{cor}\label{SIP}${\rm Aut}(\bft)$ satisfies automatic
continuity and has the small index property.\end{cor}

This follows from the general fact that automatic continuity passes from an open subgroup to the whole group.

For the next couple of results, we need some more detailed information about the structure of group actions on trees. If $g$ is an automorphism of a tree $S$ that acts without inversion, then $g$ either has a fixed point, in which case $g$ is said to be {\em elliptic}, or $g$ acts by translation on some line in the tree, in which case $g$ is said to be {\em hyperbolic}. Serre's book \cite{serre} is a good reference for more information on these concepts.

\begin{lemme}\label{generation}Suppose $\ell_g=(a_i\del i\in \Z)\subseteq \bft$ is a line and $g$ a
hyperbolic element of ${\rm Aut}(\bft)$ acting by translation on $\ell_g$ with amplitude $1$, i.e. $g\cdot a_i=a_{i+1}$, $\forall i\in
\Z$. Then ${\rm Aut}(\bft)=\big\langle {\rm Aut}(\bft,a_0)\cup \{g\}\big\rangle$.\end{lemme}

\begin{proof} Suppose $t_0$ is any element of $\bft$. Then there is an
$h\in{\rm Aut}(\bft,a_0)$ such that $t_0\in h\cdot \ell_g$ and so
$hg^{d(t_0, a_0)}\cdot a_0=t_0$. Put $k=hg^{d(t_0, a_0)}$, then
\begin{equation}
{\rm Aut}(\bft,t_0)=k{\rm Aut}(\bft,a_0)k\inv\subseteq \langle {\rm
Aut}(\bft,a_0)\cup \{g\}\big\rangle.
\end{equation}
This shows that $\langle {\rm
Aut}(\bft,a_0)\cup \{g\}\big\rangle$ contains all elliptic elements
of ${\rm Aut}(\bft)$.

Now, suppose $k$ is any other hyperbolic element of ${\rm Aut}(\bft)$
acting by translation on a line $\ell_k$ with amplitude $m=\|k\|$. Let $\alpha=(b_n,\ldots,
b_0)$ be the geodesic from $\ell_k$ to $\ell_g$ and $a_i=b_0$ be
its endpoint.

There are two cases: Either $\ell_k\cap \ell_g\neq \tom$, in which
case it is easy to find some $h\in{\rm Aut}(\bft,a_i)$ such that
$k=hg^mh\inv\in\langle {\rm Aut}(\bft,a_0)\cup \{g\}\big\rangle$. Otherwise, take some $f\in{\rm Aut}(\bft,a_i)$
such that $f\cdot a_{i+j}=b_j$ for $j=0,\ldots, n$. So
\begin{equation}
b_n=f\cdot a_{i+n}\in (f\cdot\ell_g)\cap
\ell_k=\ell_{fgf\inv}\cap\ell_k.
\end{equation}
Replacing $g$ by $fgf\inv$ we can repeat the argument above to
find an $h\in{\rm Aut}(\bft,b_n)$ such that $k=h(fgf\inv)^mh\inv\in
\langle {\rm Aut}(\bft,a_0)\cup \{g\}\big\rangle$. Therefore,
$\langle {\rm Aut}(\bft,a_0)\cup \{g\}\big\rangle$ contains all
hyperbolic elements of ${\rm Aut}(\bft)$.

We now only have the inversions left. So suppose $k$ inverts an
edge $(a,b)$ of $\bft$. Find some hyperbolic $f$ with amplitude
$\norm{f}=1$ such that its characteristic subtree $\ell_f$ passes
through $a$ and $b$ with $f\cdot a=b$. Take also some elliptic
$h\in {\rm Aut}(\bft,b)$ such that $h\cdot (f\cdot b)=a$, then
clearly
\begin{equation}
hfk\in{\rm Aut}(\bft,a)\subseteq\langle {\rm Aut}(\bft,a_0)\cup
\{g\}\big\rangle.
\end{equation}
Therefore, ${\rm Aut}(\bft)=\langle {\rm Aut}(\bft,a_0)\cup
\{g\}\big\rangle$.\end{proof}

\begin{thm}${\rm Aut}(\bft)$ is of uncountable
cofinality.\end{thm}

\begin{proof} Suppose $H_0\leq H_1\leq \ldots\leq{\rm Aut}(\bft)$ exhausts
${\rm Aut}(\bft)$. Then evidently, $H_0\cap {\rm Aut}(\bft,t_0)\leq
H_1\cap{\rm Aut}(\bft,t_0)\leq \ldots\leq{\rm Aut}(\bft,t_0)$ also
exhausts ${\rm Aut}(\bft,t_0)$ and thus by the uncountable
cofinality of the latter, there is some $n$ such that ${\rm
Aut}(\bft,t_0)\leq H_n$. Now, fix a line $\ell_g=(a_i\del i\in
\Z)\subseteq \bft$ and a hyperbolic element $g\in{\rm Aut}(\bft)$
acting by translation on $\ell_g$ with amplitude $1$, i.e. $g\cdot
a_i=a_{i+1}$, $\forall i\in \Z$ such that $t_0=a_0$. We know by
Lemma \ref{generation} that ${\rm Aut}(\bft)=\big\langle {\rm
Aut}(\bft,t_0)\cup \{g\}\big\rangle$, whence if $m\geq n$ is such
that $g\in H_m$, also ${\rm Aut}(\bft)=H_m$.\end{proof}

\begin{rem} {\rm It is easy to see that ${\rm Aut}(\bft)$ fails
property (FA). For we can just add a vertex to every edge of $\bft$
and extend the action of ${\rm Aut}(\bft)$ to the tree obtained. We
now see that ${\rm Aut}(\bft)$ acts without inversion, but does not
fix a vertex. An exercise in Serre's book (\cite{serre}, p. 34)
implies that ${\rm Aut}(\bft)$ is actually a non-trivial free product
with amalgamation.}\end{rem}

We shall now see that ${\rm Aut}(\bft)$ also satisfies the analogue of Neumann's Lemma. This is a special case of the following general fact.
\begin{prop}
Suppose $G$ is a Polish group having an open subgroup $K$ with ample generics. Then $G$ satisfies Theorem \ref{b.h.neumann} and Corollary \ref{p.m.neumann}.
\end{prop}

\begin{proof}
Suppose towards a contradiction that $\{g_iH_i\}_{i\in \N}$ covers $G$ such that no $H_i$ is open. Then as $G$ has the small index property, every $H_i$ has index $\conti$ in $G$ and hence when $G$ acts by left translation on
$$
X=G/H_0\sqcup G/H_1\sqcup\ldots,
$$
any orbit is of size $\conti$. Since $K$ is open in $G$ there are
$f_i\in G$ such that $\{f_iK\}_{i\in \N}$ cover $G$. So for any
$x\in X$, $G\cdot x=\bigcup_i f_iK\cdot x$, whence every $K$-orbit
is of size $\conti$. Applying Theorem \ref{p.m.neumann} to $K$, we
find some $k\in K$ such that $k\cdot A\cap A=\tom$, where
$A=\{H_i,g_iH_i\}_{i\in \N}$. Thus for each $i$, $kH_i\neq g_iH_i$,
contradicting that  $\{g_iH_i\}_{i\in \N}$ covers $G$. Therefore,
some $H_i$ is open proving Theorem \ref{b.h.neumann} for $G$.
Theorem \ref{p.m.neumann} for $G$ now follows from the small index
property and Theorem \ref{b.h.neumann} for $G$. \end{proof}

\subsection{Concluding remarks}
It is interesting to see that there seem to be two different paths to the study of automorphism groups of homogeneous
structures.  On the one hand, there are the methods of moieties dating at least back to Anderson \cite{anderson} which have been developed and used by a great number of
authors, and, on the other hand, the use of genericity.

The main tenet of this last section is that although ample genericity can sometimes be quite non-trivial to verify, it is
nevertheless a sufficiently powerful tool for it to be worth the effort looking for.  In particular, it provides a uniform approach
to proving the small index property, uncountable cofinality, property (FA), the Bergman property and automatic continuity.

\subsection{Some questions}
\begin{enumerate}
\item[(1)] Are there any examples of Polish groups that are
not isomorphic to closed subgroups of $S_\infty$ but have ample
generic elements? That have the small index property? {\bf
Addendum.} It has now been verified by  S. Solecki and the second
author that the homeomorphism group of the reals indeed has the
small index property, but of course is not a topological subgroup of
$S_\infty$.

\item[(2)]
Does Aut$(\bfb_\infty )$ have ample generic elements?

\item[(3)]
Can a Polish locally compact group have a comeager conjugacy class?

\item[(4)]
Characterize the generic elements of Aut$(\bfb_\infty )$ and Aut$(\bff ,\lambda )$, where $(\bff ,\lambda)$ is as in Proposition
\ref{measure}. {\bf Added in proof.} We have recently received a preprint of Akin, Glasner, and Weiss \cite{akin2} in which they give another proof of the existence of a comeager conjugacy class in the homeomorphism group of the Cantor space, a proof that also gives an explicit characterization of the elements of this class.

\item[(5)]
Is the conjugacy action of Iso$(\bfu )$ turbulent? (Vershik) Does Iso($\bfu_0$) have a dense locally finite subgroup?

\item[(6)]Suppose a Polish group has ample generics and acts by homeomorphisms on a Polish space. Is the action necessarily
continuous?

\item[(7)] Is property (T) somehow related to the existence of dense or comeager conjugacy classes? Concretely, if $G$ is a
Polish group, which is not the union of a countable chain of proper open subgroups and such that there is a diagonally dense
conjugacy class in $G^\bbN$, does $G$ have property (T)?
\end{enumerate}

\end{document}